\documentclass[11pt,reqno]{amsart}
\usepackage{amssymb,mathrsfs,color,mathtools,empheq, verbatim}
\usepackage{pinlabel}
\mathtoolsset{showonlyrefs}

\usepackage{hyperref}

\usepackage{enumerate, bbm}
\usepackage{tensor}

\usepackage{graphicx}
\usepackage{xcolor} % A package to add color.
\usepackage{slashed}
\usepackage{cite}

\usepackage[shortlabels]{enumitem}

\usepackage{geometry}

\numberwithin{equation}{section}

\newtheorem{thm}{Theorem}[section]
\newtheorem{cor}[thm]{Corollary}
\newtheorem{lem}[thm]{Lemma}
\newtheorem{prop}[thm]{Proposition}

\theoremstyle{definition} 
\newtheorem{rem}[thm]{Remark}
\newtheorem{defn}[thm]{Definition}

\theoremstyle{remark}

\def\cD {\mathcal{D}}

\def\cR {\mathcal{R}}
\def\cT {\mathcal{T}}

\newcommand{\bs}[1]{\boldsymbol{#1}}

\newcommand{\supp}{\operatorname{supp}}

\newcommand{\dist}{\operatorname{dist}}

\newcommand{\iii}{\mathrm i}

\definecolor{deepgreen}{cmyk}{1,0,1,0.5}

\newcommand{\N}{\mathbb{N}}
\newcommand{\R}{\mathbb{R}}

\newcommand{\al}{\alpha}
\newcommand{\be}{\beta}
\newcommand{\ga}{\gamma}
\newcommand{\de}{\delta}

\newcommand{\lam}{\lambda}

\newcommand{\ta}{\tau}

\newcommand{\loc}{\operatorname{loc}}

\makeatletter

\newcommand{\Rmnum}[1]{\expandafter\@slowromancap\romannumeral #1@}
\makeatother

\newcommand{\ti}{\widetilde}

\newcommand{\EQ}[1]{\begin{equation}\begin{split} #1 \end{split}\end{equation}}

\newcommand{\Del}[1]{}

\definecolor{green}{rgb}{0,0.8,0} % Redefines the color green.
\newcommand{\ud}{\mathrm{d}}

\newcommand{\alp}{\alpha}

\newcommand{\veps}{\varepsilon}

\newcommand{\bfd}{{\bf d}}
\newcommand{\bfe}{{\bf e}}

\newcommand{\calC}{\mathcal C}

\newcommand{\calJ}{\mathcal J}

\newcommand{\calL}{\mathcal L}

\newcommand{\calS}{\mathcal S}

\begin{document}
	\parindent=0pt
	\title[Bubbling for the Nonlinear Heat Flow]{Continuous-in-time Bubbling for the Energy-Critical Heat Flow}
	\author{Shrey Aryan}
\begin{abstract}
We show that any finite energy solution of the energy-critical nonlinear heat flow in dimensions $d\geq 3$ asymptotically resolves into a sum of possibly time-dependent solitons, a weak limit when the flow exists for finite time, and an error term that vanishes in the energy space. As a consequence, under the additional assumption that the initial data is nonnegative, we settle the Soliton Resolution Conjecture.
\end{abstract}
\maketitle
\section{Introduction}
\subsection{Problem Setting}
In this work, we study the long-time behavior of solutions to the energy-critical nonlinear heat flow in dimensions $d \geq 3$:
\begin{align}\label{eqn:NLH}
\partial_t u &= \Delta u + |u|^{p-1} u,\\
u(0, x) &= u_0(x),\quad u_0 \in \dot{H}^1(\mathbb{R}^d),
\end{align}
where $p:=\frac{d+2}{d-2}$. From a variational perspective, equation \eqref{eqn:NLH} is the negative gradient flow of the nonlinear energy functional
\begin{align}\label{eqn:energy}
E(u):=\frac12\int_{\R^d}|\nabla u|^2 \ud x -\frac{1}{p+1}\int_{\R^d}|u(x)|^{p+1} \ud x,
\end{align}
which arises in the study of extremizers of the Sobolev inequality and, more generally, is related to the Yamabe problem on the sphere via stereographic projection. The local well-posedness of \eqref{eqn:NLH} in $\dot{H}^1$-norm is classical and was initiated by Weissler in \cite{W1,W2}, with further contributions by Giga \cite{giga-3}, Ni–Sacks \cite{ni-sacks}, and Brezis-Cazenave \cite{brezis-cazenave}. Observe that the solutions of \eqref{eqn:NLH} are invariant under translations and the following parabolic scaling,
\begin{align*}
u(t, x) \mapsto u_{\lambda}(t, x):=\lambda^{-\frac{d-2}{2}} u\left(t / \lambda^{2}, x / \lambda\right), \quad \lambda>0.   
\end{align*}
Since $E(u_\lambda(t))=E(u(t/\lambda^2))$, the equation \eqref{eqn:NLH} is energy-critical. Testing \eqref{eqn:NLH} against $\partial_{t} u$ and integrating by parts, we observe the formal energy identity
\begin{align}\label{eqn:energy-identity-intro}
E(u(T))+\int_0^T\|\partial_tu(t)\|_{L^2}^2\,\ud t=E(u(0)),
\end{align}
for each $0<T<T_+$, where $T_+>0$ is the maximal time of existence. In particular, this implies that the nonlinear energy is nonincreasing along the flow. A \emph{bubble} or \emph{soliton} will mean a nonzero finite-energy stationary solution $W\in\dot H^1(\R^d)$ of
\begin{align}\label{eq:yamabe}
\Delta W+|W|^{p-1}W=0.
\end{align}
\subsection{Statement of the Main Results}\label{subsection:statement of main result}
We now state our first main result.
\begin{thm}[Bubbling along any sequence]\label{thm:main1}
Let $u(t)$ be a solution of~\eqref{eqn:NLH} with initial data $u_0 \in \dot{H}^1$. Let $T_+ = T_+(u_0) \in (0, \infty]$ denote its maximal time of existence and assume that $u(t)$ has uniformly bounded energy, i.e., $\sup_{t\in [0,T_+)}\|\nabla u(t)\|_{L^2}<\infty.$ Then the following hold:
\newline 
$\mathrm{(i)}$ If $T_+ < \infty$, then there exist a map $u^*\in\dot H^1$,  an integer $K \ge 1$, and distinct points $\{x^i\}_{i =1}^K \subset \R^d$ such that the following holds. Let $t_n \to T_+$ be any time sequence. After passing to a subsequence (still denoted by $t_n$) we can associate with each $i\in\{1,\dots,K\}$ an integer $J_i\in\N$ with $J_i\geq1$, sequences $a_{j, n}^{i} \in \R^d$ and $\lam_{ j, n}^{i} \in (0, \infty)$ for each $j \in \{1, \dots, J_i\}$,  with 
\begin{align*}
\lim_{n\to \infty} \frac{|a_{j,n}^i-x^i|+\lam_{ j, n}^{i}}{\sqrt{T_+- t_n}}  =0,
\end{align*}
and nonzero bubbles $W_{1}^{i},  \dots, W_{J_i}^{i}$ such that 
\EQ{ \label{eq:ao-finite} 
\lim_{n \to \infty} \bigg( \frac{\lam_{ j, n}^{i}}{\lam_{ k, n}^{i}} + \frac{\lam_{k, n}^{i}}{\lam_{ j, n}^{i}} + \frac{ | a_{j, n}^{i} - a_{k, n}^{i}|}{ \lam_{j, n}^{i}}\bigg)  = \infty \quad\textrm{for all }j \neq k, 
}
and 
\EQ{\label{eqn:decomp-finite}
u(t_n) = u^* + \sum_{i = 1}^{K}  \sum_{j =1}^{J_i} W^i_{j} [a^i_{j,n},\lam^{i}_{j,n}]  + o_{\dot{H}^1}(1),
}
where $W^i_j[a^i_{j,n},\lam^i_{j,n}](x)=(\lam^i_{j,n})^{-(d-2)/2}W^i_j((x-a^i_{j,n})/\lam^i_{j,n})$ and $o_{\dot H^1}(1)$ denotes a term converging strongly to zero in $\dot H^1$.
\newline 
$\mathrm{(ii)}$ If $T_+ = \infty$, then let $t_n \to \infty$ be any time sequence. After passing to a subsequence we can find an integer $K \ge 0$, and when $K\geq 1$, sequences $a_{j, n} \in \R^d$ and $\lam_{ j, n} \in (0, \infty)$  for each $j \in \{1, \dots, K\}$, with
\EQ{
 \lim_{n \to \infty} \frac{|a_{j, n}| + \lam_{ j, n}}{\sqrt{t_n}} = 0
}
and nonzero bubbles $W_1, \dots, W_{K}$, so that 
\EQ{\label{eqn:ao-global}
\lim_{n \to \infty} \bigg( \frac{\lam_{j, n}}{\lam_{k, n}} + \frac{\lam_{k, n}}{\lam_{j, n}} + \frac{ | a_{ j, n} - a_{k, n}|}{ \lam_{j, n}} \bigg) = \infty \quad\textrm{for all }j \neq k,
}
and
\EQ{\label{eqn:decomp-global}
u(t_n)  = \sum_{j =1}^{K} W_{j} [a_{j,n},\lam_{j,n}]  +o_{\dot{H}^1}(1)
}
where $W_{j} [a_{j,n},\lam_{j,n}](x)=(\lam_{j,n})^{-(d-2)/2} W_j((x-a_{j,n})/\lam_{j,n})$ and the error term $o_{\dot{H}^1}(1)\to 0$ strongly in $\dot{H}^1.$ When $K=0$, the sum is empty and $u(t_n)\to0$ in $\dot H^1$.
\end{thm}
\begin{rem}
Since our arguments do not use modulation techniques, we do not obtain continuous families of scaling and translation parameters of the kind available in the radial case (cf. Theorem 1.1 in \cite{nlh-radial}). Theorem~\ref{thm:main1} gives a bubble-tree decomposition along a subsequence of any prescribed sequence of times, whereas Struwe \cite{struwe-global} previously obtained such a decomposition only along a well-chosen (or Palais--Smale) sequence. Our next main result, Theorem~\ref{thm:main}, shows that the distance to the set of multi-bubble configurations localized to suitable parabolic regions vanishes continuously in time.
\end{rem}
The bubbles obtained in the above decomposition may depend on the sequence of times. This issue also arises for the harmonic map heat flow (cf. Remark 1.8 in \cite{lawrie-harmonic-map-nonradial}). For \eqref{eqn:NLH}, the additional assumption of nonnegative finite energy initial data allows us to choose the bubble profiles independently of the time sequence, yielding the full Soliton Resolution Conjecture. Before stating the result, recall that the classification theorem of Caffarelli--Gidas--Spruck \cite{CGS} implies that every nontrivial, nonnegative finite-energy solution of \eqref{eq:yamabe} is, up to scaling and translation, the Aubin--Talenti bubble
\begin{align}\label{eqn:talenti}
\widetilde W(x):=\left(1+\frac{|x|^2}{d(d-2)}\right)^{-\frac{d-2}{2}}.
\end{align}
This elliptic rigidity identifies the bubbles in Theorem~\ref{thm:main1} up to the symmetries of the equation \eqref{eqn:NLH}. Consequently, we obtain
\begin{cor}\label{cor:src}
Let $u(t)$ be a solution to~\eqref{eqn:NLH} with nonnegative initial data $u_0\geq 0$ and $u_0 \in \dot{H}^1$. Let $T_+ = T_+(u_0) \in (0, \infty]$ denote its maximal time of existence and assume that $u(t)$ has finite energy, i.e., $\sup_{t\in [0,T_+)}\|\nabla u(t)\|_{L^2}<\infty.$ Then, after modifying the translation and scaling parameters, every bubble in \eqref{eqn:decomp-finite} and \eqref{eqn:decomp-global} may be taken to be $\widetilde W$. Moreover, the number of bubbles denoted by $K$ in the global in time case and by $\sum_{i=1}^KJ_i$ in the finite time case are independent of the chosen sequence.
\end{cor}
\begin{rem}
Corollary~\ref{cor:src} does not assert the uniqueness of the scale and translation parameters. In general, the scaling parameter can oscillate as demonstrated, for instance, by Theorem 3 in \cite{harada2025oscillatory}.
\end{rem}
To state the continuous-in-time bubbling result, we first introduce the definitions of \emph{scale} and~\emph{center} for a nonzero stationary solution. Let $S=S(d)>0$ be the sharp constant in the Sobolev inequality on $\R^d$,
\begin{align}
\|u\|_{L^{p+1}}\leq S\|\nabla u\|_{L^2}
\end{align}
for all $u\in \dot{H}^1.$ Then observe that for any nonzero stationary solution $W\in \dot{H}^1$ we have
\begin{align}
\int_{\R^d} |\nabla {W}|^2  \ud x = \int_{\R^d} |W|^{p+1}  \ud x. 
\end{align}
By the variational characterization of the Sobolev inequality, the best constant (or equality) is attained by a positive stationary solution $\ti{W}$, i.e.,
\EQ{\label{eqn:minimizer-sobolev}
\int_{\R^d} |\ti{W}|^{p+1} \ud x = S^{p+1}\left(\int_{\R^d} |\nabla \ti{W}|^2 \ud x\right)^{(p+1)/2}= \int_{\R^d} |\nabla \ti{W}|^2  \ud x
}
implying that for any positive stationary solution, we have
\EQ{
\int_{\R^d} |\nabla \ti{W}|^2 \ud x = S^{-d}.
}
If $W$ changes sign, testing \eqref{eq:yamabe} with $W^+=\max\{W,0\}$ and $W^-=\max\{-W,0\}$ and applying the Sobolev inequality gives
\begin{align}
    \int_{\R^d} |\nabla W|^{2} \ud x = \int_{\R^d} |\nabla W^+|^2 \ud x + \int_{\R^d} |\nabla W^-|^2 \ud x \geq 2 S^{-d}.
\end{align}
In particular, we deduce that any nonzero stationary solution $W\in \dot H^1$ satisfies $\|\nabla W\|_{L^2}^2\geq S^{-d}.$ Denote this minimal Dirichlet energy by $\bar E_*:=S^{-d}$. Define $\bar{E}(u;\Omega):=\|\nabla u\|_{L^2(\Omega)}^2$ for any $\Omega\subset \R^d$ and $u\in \dot{H}^1$ with the convention that $\bar{E}(u):=\bar{E}(u;\R^d).$ Since $E(W)=\bar E(W)/d$ for every stationary solution $W$, the corresponding minimal nonlinear energy is $\bar E_*/d$. Let $B(x,r)\subset \R^d$ denote the ball centered at $x\in \R^d$ with radius $r>0.$ Given any nonzero stationary solution $W:\R^d\to \R$, we define its scale and center as follows:
\begin{defn}[Scale of a stationary solution] 
\label{def:scale}
Let $\gamma_0 \in (0, \bar{E}_*/2)$. Then the scale associated to a nontrivial stationary solution $W$, denoted by $\lam( W; \gamma_0)$, is defined by 
\EQ{
\lam( W; \gamma_0):= \inf\{ \lam \in (0, \infty) \mid \exists\,a \in \R^d\textrm{ such that } \bar{E}( W; B(a, \lam)) \ge \bar{E}(W) - \gamma_0\}.
}
\end{defn}
\begin{defn}[Center of a stationary solution]
\label{def:center}
Let $\gamma_0\in(0,\bar E_*/2)$ and let $\lam(W;\gamma_0)$ be the scale
of a nonzero stationary solution $W$. A center of $W$, denoted by
$a(W;\gamma_0)\in\R^d$, is any point satisfying
\EQ{\label{eq:cent def}
\bar E\bigl(W;B(a(W;\gamma_0),\lam(W;\gamma_0))\bigr)
\geq
\bar E(W)-\gamma_0.
}
\end{defn}
These quantities are well defined, as shown in Lemma~\ref{lem:well-defn-scale-center}. Since our main result says that finite energy solutions of \eqref{eqn:NLH} eventually approach a sum of stationary solutions, it will be convenient to define their sum, which we will often refer to as a multi-bubble configuration.
\begin{defn}[Multi-bubble configuration] 
Let $K \in \{0,1, 2, \dots\}$. A $K$-multi-bubble configuration is the sum
\EQ{
\mathbf{W}(x) = \sum_{j=1}^K W_j(x), 
}
where $W_j\in\dot H^1(\R^d)$ are nonzero stationary solutions. By convention, if $K=0$ then $\mathbf{W}\equiv 0$. To emphasize the dependence of $\mathbf{W}$ on the collection $\{W_j\}_{j=1}^K$, we will occasionally write $\mathbf{W}=\mathbf{W}(\vec{W})$, where $\vec{W}=(W_1,\ldots, W_K).$
\end{defn} 
For each nonzero stationary $W$ and $\gamma_0\in(0,\bar E_*/2)$, fix one center $a(W;\gamma_0)$ satisfying Definition~\ref{def:center}. When $\gamma_0$ is fixed, write $\lambda(W):=\lambda(W;\gamma_0)$ and $a(W):=a(W;\gamma_0)$. Next, we quantify the distance of a function to some multi-bubble configuration.
\begin{defn}[Localized distance to a multi-bubble configuration]  \label{def:d}
Given
\begin{enumerate}
    \item radii $\xi, \rho, \nu \in (0, \infty)$, such that $\xi \le \rho\le \nu$;
    \item a point $y\in\R^d$, a function $u\in \dot{H}^1(B(y,\nu))$, and $\gamma_0\in(0,\bar E_*/2)$;
    \item a nonnegative integer $K\geq 0$ and nonzero stationary solutions $W_1, \dots, W_K$ with centers $a(W_j;\gamma_0)\in B(y, \xi)$ and scales $\lambda(W_{j};\gamma_0) \in (0, \infty)$ for each $j \in \{1, \dots, K\}$;
    \item a collection of radii $ \vec \nu = (\nu,\nu_1, \dots, \nu_K) \in (0, \infty)^{K+1}$ such that $B(a(W_j), \nu_j) \subset B(y, \xi)$ and smaller scales $\vec \xi = ( \xi,\xi_1, \dots,\xi_K) \in (0, \infty)^{K+1}$ such that $\xi_j < \lam(W_j;\gamma_0)$ for each $j \in \{1, \dots, K\}$. Denote $I_j=\{k\in \{1,\ldots,K\}\setminus\{j\}:B(a(W_k),\xi_j)\subset B(a(W_j),\nu_j)\}.$
    % \item excised balls 
    % \EQ{
    %  B_{j}^*:=B( a(W_j), \nu_j) \setminus \cup_{k\in \calI_j }B(a(W_k), \xi_j),
    % }
    % where $\calI_j:= \{ k \neq j \mid  B(a(W_k), \xi_j) \subset B(a(W_j), \nu_j)\}$;
    % \item and cutoff functions $\phi,\psi\in C^\infty_c(\R^d)$ supported on $B(y,2\rho)$ and on the annulus $B(y,2\nu)\setminus B(y,\xi/2)$ respectively.
\end{enumerate}
Then the localized distance is defined as 
 \EQ{
\bfd_{\gamma_0}(u,\mathbf W;B(y,\rho);\vec\nu,\vec\xi)
&:= \bar E\bigl(u-\mathbf W(\vec W);B(y,\rho)\bigr)
+\|u-\mathbf W(\vec W)\|_{L^{p+1}(B(y,\rho))}^{p+1} \notag\\
&\quad+\bar E(u;B(y,\nu)\setminus B(y,\xi))
+\|u\|_{L^{p+1}(B(y,\nu)\setminus B(y,\xi))}^{p+1} \notag\\
&\quad+\frac{\xi}{\rho}+\frac{\rho}{\nu}
+\sum_{j\neq k}
\left(
\frac{\lambda(W_j)}{\lambda(W_k)}
+\frac{\lambda(W_k)}{\lambda(W_j)}
+\frac{|a(W_j)-a(W_k)|}{\lambda(W_j)}
\right)^{-\frac{d-2}{2}} \notag\\
&\quad+\sum_j
\left(
\frac{\lambda(W_j)}{\dist(a(W_j),\partial B(y,\xi))}
+\frac{\lambda(W_j)}{\nu_j}
+\frac{\xi_j}{\lambda(W_j)}
\right) \notag\\
&\quad+\sum_j\sum_{k\in I_j}
\frac{\xi_j}{\dist(a(W_k),\partial B(a(W_j),\nu_j))}.
  }
\end{defn} 
Minimizing over all parameters in the preceding definition yields the
following proximity function.
\begin{defn}[Localized multi-bubble proximity function]\label{defn:delta}
Given $T_+\in(0,\infty]$, $u:[0,T_+)\to\dot H^1(\R^d)$, $t\in[0,T_+)$, $y\in\R^d$, $\rho>0$, and $\gamma_0\in(0,\bar E_*/2)$, define
\EQ{
\bs\de_{\gamma_0}( u(t); B(y, \rho)) := \inf_{\mathbf{W}, \vec \nu, \vec \xi}  \bfd_{\gamma_0}( u(t),\mathbf{W};  B(y, \rho); \vec \nu, \vec \xi) 
}
where the infimum is taken over all nonnegative integers $K$, all $K$-multi-bubble configurations, and all parameters $\vec\nu\in(0,\infty)^{K+1}$ and $\vec\xi\in(0,\infty)^{K+1}$ as in Definition~\ref{def:d}.
\end{defn} 
We next state the localized form of the bubbling theorem.
\begin{thm}[Continuous-in-time bubbling] \label{thm:main} 
Let $u(t)$ be a solution of~\eqref{eqn:NLH} with initial data $u_0 \in \dot{H}^1$. Let $T_+ = T_+(u_0) \in (0, \infty]$ denote its maximal time of existence and assume that $u(t)$ has finite energy, i.e., $\sup_{t\in [0,T_+)}\bar{E}(u(t))<\infty.$  Then there exists $\gamma_0 = \gamma_0(\bar{E}_*,d)>0$ with $\gamma_0\in(0,\bar E_*/2)$ such that the following holds:

$\mathrm{(i)}$ If $T_+< \infty$, then for any $y \in \R^d$,  
\EQ{
 \lim_{t \to T_+} \bs \de_{\gamma_0}\big( u(t); B(y, \sqrt{T_+-t})\big)  = 0. 
 }
Moreover, let $t_n \to T_+$ be any sequence and let  $B(y_n, \rho_n)$ be any sequence of balls for $y_n\in \R^d$, $\rho_n>0$ such that $B(y_n, R_n \rho_n) \subset B(y, \sqrt{T_+-t_n})$ for some sequence $R_n>0$ such that $R_n \to \infty$. Suppose $\al_n,\be_n\in(0,\infty)$ are sequences with $\al_n\to0$, $\be_n\to\infty$, $\lim_{n \to \infty} \be_n R_n^{-1} = 0$, and
  \EQ{
 \lim_{n \to \infty}   \left[\bar{E}\big(u(t_n); B(y_n, \be_n \rho_n) \setminus B( y_n, \al_n \rho_n)\big) + \|u(t_n)\|^{p+1}_{L^{p+1}(B(y_n, \be_n \rho_n) \setminus B( y_n, \al_n \rho_n))}\right]  = 0. 
 }
 Then, 
 \EQ{
 \lim_{ n \to \infty} \bs \de_{\gamma_0}\big( u(t_n); B( y_n, \rho_n)\big) = 0. 
 }

$\mathrm{(ii)}$ If $T_+ = \infty$, then for every $y \in \R^d$, 
\EQ{
 \lim_{t \to \infty} \bs \de_{\gamma_0}\big( u(t); B(y, \sqrt{t})\big)  = 0. 
 }
 Moreover, let $t_n \to \infty$ be any sequence and let  $B(y_n, \rho_n)$ be any sequence of balls for $y_n\in \R^d$, $\rho_n>0$ such that $B(y_n, R_n \rho_n) \subset B(y, \sqrt{t_n})$ for some sequence $R_n>0$ such that $R_n \to \infty$. Suppose $\al_n,\be_n\in(0,\infty)$ are sequences with $\al_n\to0$, $\be_n\to\infty$, $\lim_{n \to \infty} \be_n R_n^{-1} = 0$ and
  \EQ{
 \lim_{n \to \infty} \left[ \bar{E}\big(u(t_n); B(y_n, \be_n \rho_n) \setminus B( y_n, \al_n \rho_n)\big) + \|u(t_n)\|^{p+1}_{L^{p+1}(B(y_n, \be_n \rho_n) \setminus B( y_n, \al_n \rho_n))} \right]= 0. 
 }
 Then, 
 \EQ{
 \lim_{ n \to \infty} \bs \de_{\gamma_0}\big( u(t_n); B( y_n, \rho_n)\big) = 0. 
 }
\end{thm} 
\subsection{Background and Motivation}
A fundamental problem in the analysis of nonlinear partial differential equations (PDEs) is to describe the long-time behavior of their solutions. The {Soliton Resolution Conjecture} asserts that any finite-energy solution to a nonlinear PDE asymptotically decomposes into a sum of decoupled solitons that are stationary solutions of the underlying equation, a radiation term that behaves like a solution to the linear flow, and an error term that vanishes in the natural energy norm. This conjecture arose from the numerical experiments of Fermi–Pasta–Ulam–Tsingou \cite{pasta} and Zabusky–Kruskal \cite{zabusky}, which provided evidence that it holds for the Korteweg--de Vries (KdV) equation. Since then, the problem has been extensively studied for the KdV equation as well as for several other integrable models.

Beyond integrable systems, analogues of the Soliton Resolution Conjecture have emerged across various areas of mathematics. In general relativity, the {Final State Conjecture} (cf. \cite{Klainerman}) predicts that generic solutions to Einstein's field equations asymptotically approach a finite number of stationary solutions, namely Kerr black holes, moving apart from each other. In geometric analysis, Soliton Resolution arises naturally in the study of gradient flows associated with conformally invariant variational problems. For example, pioneering works of Struwe \cite{struwe,struwe1994yang}, Qing \cite{qing}, Qing–Tian \cite{qing-tian}, and Hong–Tian \cite{hong-tian} have established bubbling or Soliton Resolution along a well-chosen sequence of times for the harmonic map and Yang–Mills heat flows. 

Motivated by these parabolic works, we study the energy-critical nonlinear heat flow in dimension $d\geq 3.$ Our main results, Theorem \ref{thm:main1} and Theorem \ref{thm:main}, establish a continuous-in-time bubble-tree decomposition for all $\dot{H}^1$ bounded solutions of \eqref{eqn:NLH}. More precisely, any solution with uniformly bounded $\dot{H}^1$-norm decomposes into a sum of solitons that may vary along different time sequences, a radiation term that is asymptotically trivial or captured by a weak limit in $\dot{H}^1$, and an error term that vanishes in the energy space. Moreover, when the initial data is nonnegative, Corollary \ref{cor:src} shows that Theorem \ref{thm:main1} implies the Soliton Resolution Conjecture, since positive solitons have been classified and are unique up to the symmetries of the equation due to \cite{obata1972conformal,CGS}. Theorem \ref{thm:main1} extends Struwe’s classical compactness result \cite{struwe-global} by obtaining a similar decomposition along a subsequence of any prescribed sequence of times, whereas Struwe's result uses a well-chosen sequence. It also strengthens Theorem 1.4 in \cite{ishiwata2018potential} in the global-in-time setting by upgrading the stated remainder convergence from $L^{p+1}$ to $\dot H^1$. Corollary~\ref{cor:src} establishes Soliton Resolution for a non-integrable PDE beyond radial symmetry and without restrictions on the size of the initial data. Ishiwata \cite{ishiwata2018potential} obtained a related result in the global setting, while our work treats both finite-time blowup and global-in-time solutions.

To explain the significance of our result, we now review some key developments in the literature. In the integrable setting, where tools such as the inverse scattering transform are available, the conjecture is well understood for models including the KdV equation \cite{eckhaus1983emergence}, the modified KdV equation \cite{schuur2006asymptotic}, the one-dimensional cubic nonlinear Schrödinger equation (NLS) \cite{borghese2018long}, the derivative NLS \cite{jenkins2019derivative}, and, more recently, the Calogero–Moser derivative NLS \cite{kim2024soliton}.

For non-integrable equations with radial symmetry, where the solitons do not move in space, the conjecture has been settled for the nonlinear wave equation \cite{Duyckaerts-wave:d=3,Duyckaerts-wave-2,Duyckaerts-wave-5,Duyckaerts-wave-6,jia-kenig,collot2022soliton,lawrie-wave,lawrie-wavemaps}, damped Klein–Gordon equation \cite{burq2017long,gu2023soliton}, equivariant self-dual Chern–Simons–Schrödinger equation \cite{oh-chern-simons}, equivariant harmonic map heat flow \cite{lawrie-harmonic-map}, and energy-critical nonlinear heat flow \cite{nlh-radial}. The last two works led to the first complete classification of the bubble-tree term dynamics for the equivariant harmonic map heat flow and the radial energy-critical nonlinear heat flow by Kim and Merle in \cite{kim2024classification}.

For non-integrable equations without radial symmetry, Soliton Resolution is
known in one dimension for the damped Klein--Gordon equation
\cite{cote2021long}; near a finite collection of solitons for the
energy-critical nonlinear heat flow and the damped Klein--Gordon equation
\cite{collot,ishizuka2023global,ishizuka2025long}; continuously in time for
the harmonic map heat flow \cite{lawrie-harmonic-map-nonradial}; and along
a sequence of times, in dimensions $3\leq d\leq5$, for the energy-critical
nonlinear wave equation \cite{jia-dkm}.

Regarding \eqref{eqn:NLH} in the energy-critical setting, Wang--Wei \cite{wang-wei} proved that in dimensions $d\ge 7$ and for bounded nonnegative initial data $u_0\ge 0$, any finite-time blowup must be of Type-I, thereby ruling out finite-time Type-II blowup in higher dimensions. When $3\leq d \leq 6$, examples of Type-II blowup have been constructed, for instance by del~Pino--Musso--Wei--Zhang--Zhang \cite{del2020type}, Schweyer \cite{schweyer2012type}, del~Pino--Musso--Wei \cite{del2019type}, and Harada \cite{harada2020type}. Regarding
positive solutions, Wei--Zhang--Zhou \cite{wei2024fila} and Harada \cite{harada2025oscillatory} have constructed positive solutions near nonnegative bubbles at infinite time when $d=4$ and $d=6$ respectively.

Our argument treats every dimension $d\geq3$. In this setting, solitons have only weak decay, need not be radial, may translate in space, and may exhibit pathological behavior (cf. \cite{ding1986conformally,del2011large,del2013torus}). We also impose no restrictions on the size of the initial data, so the nonlinear energy is, in general, non-coercive, unlike in \cite{kenig-merle,roxanas}. We address these difficulties using profile decompositions. The resulting arguments seem robust and appear to be adaptable to other nonlinear parabolic flows.

\subsection{Proof Sketch}\label{sec:proof-sketch}
Let $u$ satisfy the hypotheses of Theorem~\ref{thm:main1} and let $t_n\to T_+$ be any sequence of times. Suppose first that $T_+<\infty$, let $u^*$ denote the weak limit and let $\calS$ denote the singular set as defined in Theorem~\ref{thm:lwp}. Set $v_n:=u(t_n)-u^*$. After passing to a subsequence, we apply Proposition~\ref{prop:profile-decomp} to $\{v_n\}_{n\in\N}$ and obtain a profile decomposition. The strong convergence of $u(t)$ to $u^*$ away from $\calS$, together with Lemmas~\ref{lem:finite-time-no-escape}, \ref{lem:ss-bu}, and~\ref{lem:ss-bu-lp}, shows that if $(a_n^j,\lambda_n^j)\in \R^d\times (0,\infty)$ are the translation and scaling parameters of a nonzero profile indexed by an integer $j\in \N$, then there exists a point depending on $j$ denoted by $x^{i(j)}\in\calS$ such that
\begin{align}
        \frac{|a_n^j-x^{i(j)}|+\lambda_n^j}
        {\sqrt{T_+-t_n}}
        \to0,
\end{align}
as $n\to \infty.$ If $T_+=\infty$, after passing to a subsequence, we instead apply Proposition~\ref{prop:profile-decomp} to $v_n:=u(t_n)$ to obtain a profile decomposition and denote the profile parameters by $(a_n^j,\lambda_n^j)\in \R^d\times (0,\infty)$ indexed by an integer $j\in \N$. In this case, Lemmas~\ref{lem:ss-global} and~\ref{lem:ss-global-weaker} give
\begin{align}
        \frac{|a_n^j|+\lambda_n^j}{\sqrt{t_n}}
        \to0,
\end{align}
as $n\to \infty$ for every nonzero profile. Rescaling the solution by the parameters of
each profile and using the finite dissipation in
\eqref{eq:tension-L2}, we show that every nonzero profile is a stationary
solution of \eqref{eq:yamabe} (cf. Lemma~\ref{lem:arbitrary-stationary-profiles}).

By the Pythagorean expansion and the lower bound $\bar E(W)\geq\bar E_*$ for every nonzero stationary profile, only finitely many nonzero profiles occur. Thus, for some integer $J\geq0$ (with $J\geq1$ if $T_+<\infty$), the profile decomposition gives
\begin{align}\label{eqn:upgraded-profile-decomp}
        v_n
        =
        \sum_{j=1}^J
        W^j[a_n^j,\lambda_n^j]
        +
        r_n,
        \qquad
        \|r_n\|_{L^{p+1}}\to0,
\end{align}
as $n\to \infty.$ The main point is to prove that
$\|\nabla r_n\|_{L^2}\to0$, as $n \to \infty.$ By the energy identity for the gradient in
Proposition~\ref{prop:profile-decomp} and the fact that
$\|r_n\|_{L^{p+1}}\to0$ as $n\to \infty$ we have
\begin{align}\label{eqn:r_n-identity}
  -\int_{\R^d}
        \nabla u(t_n)\cdot\nabla r_n\,\ud x
        +
        \int_{\R^d}
        |u(t_n)|^{p-1}u(t_n)r_n\,\ud x
        =
        -\|\nabla r_n\|_{L^2}^2+o_n(1).
\end{align}
By \eqref{eqn:NLH}, the expression on the left-hand side is obtained by testing $\partial_tu(t_n)$ against $r_n$. In the finite-time case, the $u^*$-term arising from the first integral vanishes since
\begin{align}\label{eqn:weak-rn-term}
\int_{\R^d}\nabla u^*\cdot\nabla r_n\,\ud x=o_n(1).
\end{align}
This follows from the fact that $v_n\rightharpoonup 0$ as $n\to \infty$ in $\dot{H}^1$ and every profile in \eqref{eqn:upgraded-profile-decomp} converges weakly to zero in
$\dot H^1$. This implies that $r_n\rightharpoonup0$ as $n\to \infty$ in $\dot{H}^1$, which gives \eqref{eqn:weak-rn-term}.

It therefore remains to show that the expression on the left-hand side of \eqref{eqn:r_n-identity} converges to zero. Suppose first that $d\geq5$. For $0<t<T<T_+$, we show the following main estimate
\begin{align}
        \int_0^{T-t}
        \|e^{s\Delta}\partial_tu(t)\|_{L^2}^2\,\ud s
        \lesssim
        \int_t^T
        \|\partial_tu(s)\|_{L^2}^2\,\ud s.
\end{align}
We prove this estimate by analyzing the PDE satisfied by the difference
quotients $h^{-1}[u(s+h)-u(s)]$, rather than $\partial_tu$, because of regularity issues. In the global case, letting $T\to\infty$ in the above estimate and using finite total dissipation gives $\|\partial_tu(t)\|_{\dot H^{-1}}\to0$ as $t\to\infty$. Therefore, by \eqref{eqn:NLH} and the definition of the $\dot H^{-1}$-norm, this implies that the left-hand side of
\eqref{eqn:r_n-identity} tends to zero.

In the finite-time case, set $\rho_n:=\sqrt{T_+-t_n}$. The same estimate gives
\begin{align}
        \|P_{>\rho_n^{-1}}\partial_tu(t_n)\|_{\dot H^{-1}}
        \to0,
\end{align}
as $n\to \infty.$ Here $P_{>\rho_n^{-1}}$ denotes the usual Fourier projection defined in the following subsection. On the other hand, the remainder can be localized to balls of radius
$\alpha_n\rho_n$, where $\alpha_n\to0$ as $n\to \infty.$ If we denote $f_n:=\chi_nr_n$, where $\chi_n$ is a suitable cutoff localizing to these balls, then
\begin{align}
        \|r_n-f_n\|_{\dot H^1}\to0,
        \qquad
        \|\nabla P_{\leq\rho_n^{-1}}f_n\|_{L^2}\to0,
\end{align}
as $n\to \infty.$ By \eqref{eqn:NLH} and the Sobolev inequality,
\EQ{
\|\partial_tu(t_n)\|_{\dot H^{-1}}\lesssim\|u(t_n)\|_{\dot H^1}+\|u(t_n)\|_{\dot H^1}^p\lesssim1.
} Together with the uniform bound for
$\|\partial_tu(t_n)\|_{\dot H^{-1}}$, these estimates show that the
left-hand side of \eqref{eqn:r_n-identity} tends to zero. Hence $\|\nabla r_n\|_{L^2}\to0$ as $n\to\infty$ when $d\geq5$.

It remains to consider $d\in\{3,4\}$. In these dimensions, the nonlinear
energy is analytic. Lemma~\ref{lem:no-ghost-LS} shows that the nonlinear
energy is locally constant on the set of stationary solutions. Since
$\dot H^1(\R^d)$ is second countable, the set of stationary nonlinear
energy values is countable. Moreover, every stationary solution $W$
satisfies
\begin{align}
        \bar E(W)
        =
        \|\nabla W\|_{L^2}^2
        =
        \|W\|_{L^{p+1}}^{p+1},
        \qquad
        E(W)=\frac1d\bar E(W),
\end{align}
so the set of stationary Dirichlet energy values is also countable. The energy identity and the uniform $\dot H^1$ bound imply that $L_E:=\lim_{t\to T_+}E(u(t))$ exists and is finite. Next, set $b(t):=\|u(t)\|_{L^{p+1}}^{p+1}$. The set of limit points of $b(t)$ as $t\to T_+$ is a nonempty compact interval because $b$ is bounded and continuous. On the other hand, by Lemma~\ref{lem:arbitrary-stationary-profiles} and, in finite time, \eqref{eq:no-ghost-finite-BL-splitting}, every limit point is a finite sum of stationary Dirichlet energies, together with $\|u^*\|_{L^{p+1}}^{p+1}$ in finite time. We use this to deduce that the limit $L_b:=\lim_{t\to T_+}b(t)$ exists.

Finally, finite total dissipation allows us to choose a sequence of times for which the corresponding localized time slices form a Palais--Smale sequence. Lemma~\ref{lem:whole-space-PS} gives a decomposition along this sequence
with a remainder converging to zero in $\dot H^1$. Along an arbitrary
sequence, after passing to a further subsequence, suppose that
$\|\nabla r_n\|_{L^2}^2\to\gamma$ as $n\to \infty$ for some $\gamma\in [0,\infty).$ We aim to show that $\gamma=0.$ The Pythagorean identities give
\begin{align}
    L_E&=E(u^*)+\frac1d\left(L_b-\|u^*\|_{L^{p+1}}^{p+1}\right)+\frac12\gamma,\text{ if } T_+<\infty \\
    L_E&=\frac1dL_b+\frac12\gamma, \text{ if } T_+=\infty.
\end{align}
Along the Palais--Smale sequence the same identities hold with $\gamma=0$. Hence $\gamma=0$ along the arbitrary sequence. This proves that the remainder term vanishes in the stronger $\dot{H}^1$ norm (cf. Lemma~\ref{lem:arbitrary-time-no-ghost}).

Theorem~\ref{thm:main1} follows from the preceding arguments and from Lemma~\ref{lem:arbitrary-singular-coverage}, which shows that in the finite-time case each point in the singular set is associated to at least one profile. To prove Theorem~\ref{thm:main}, we also use Lemmas~\ref{lem:ss-bu}, \ref{lem:ss-bu-lp}, \ref{lem:ss-global}, and~\ref{lem:ss-global-weaker}, which rule out Dirichlet and $L^{p+1}$ energy concentration at the self-similar parabolic scale, together with Lemma~\ref{lem:localize-strong-bubbles}, which localizes the bubble tree decomposition to the self-similar parabolic region.

\subsection{Notation and Conventions} 
We use the following conventions in this paper:
\begin{itemize}
\item We denote mixed Lebesgue spaces by $L_t^rL_x^q$, where the subscripts indicate the variables of integration.
\item We define $\dot H^1(\R^d)$ as the completion of $C_c^\infty(\R^d)$ under the norm $\|\nabla f\|_{L^2}$, for $f\in C^\infty_c(\R^d)$ and set
\begin{align}
\dot H^{-1}(\R^d)
:=
\bigl\{
F\in\cD'(\R^d):
|F[\varphi]| \leq C\|\nabla\varphi\|_{L^2}, \forall \varphi \in C_c^\infty(\R^d) \textrm{ and some  }C>0
\bigr\},
\end{align}
with $\|F\|_{\dot H^{-1}} = \sup \{|F(\varphi)|: \varphi\in C_c^\infty(\R^d),
\|\nabla\varphi\|_{L^2}\leq1\}$, where $\cD'(\R^d)$ denotes the space of distributions and
$F[\varphi]$ denotes the action of the distribution $F$ on the test function $\varphi$. 
\item Some constants that will occur frequently include $p:=\frac{d+2}{d-2}$, for $d\geq3$, and $\bar E_*:=S^{-d}$, which is the minimum Dirichlet energy of a nonzero stationary solution of
\eqref{eq:yamabe}. Furthermore, the inequality $A\lesssim B$ means that $A\leq CB$ for some constant $C>0$, while $A\simeq B$ means that $A\lesssim B$ and $B\lesssim A$. Unless stated otherwise, implicit constants may depend on $d$ and the finite energy bound of the solution of \eqref{eqn:NLH}.
\item An open ball is denoted by $B(x,r):=\{z\in\R^d:|z-x|<r\}$. A backward parabolic cylinder is denoted
by $Q_r(x,t):=B(x,r)\times(t-r^2,t)$ for $x\in\R^d$, $t\in\R$, and $r>0.$ For convenience, we write
$Q_1:=Q_1(0,0)=B(0,1)\times(-1,0)$.
\item We define the nonlinear energy and the Dirichlet energy densities by
\begin{align}
\bfe(u)
:=
\frac12|\nabla u|^2-\frac1{p+1}|u|^{p+1},
\qquad
\bar{\bfe}(u):=|\nabla u|^2.
\end{align}
Given a measurable set $A\subset\R^d$, we define
\begin{align}
E(u;A):=\int_A\bfe(u(x))\ud x,\quad 
\bar E(u;A):=\int_A\bar{\bfe}(u(x))\ud x.
\end{align}
If $\phi\in C^\infty(\R^d)\cap W^{1,\infty}(\R^d)$ satisfies
$0\leq\phi\leq1$, we also define
\begin{align}
E_\phi(u;A)
:=
\int_A\bfe(u(x))\phi^2(x)\ud x,
\quad
\bar E_\phi(u;A)
:=
\int_A\bar{\bfe}(u(x))\phi^2(x)\ud x.
\end{align}
We use the abbreviations
$E_\phi(u):=E_\phi(u;\R^d)$ and
$\bar E_\phi(u):=\bar E_\phi(u;\R^d)$.
\item A standard cutoff function is a radial function
$\chi\in C_c^\infty(\R^d)$ satisfying $0\leq\chi\leq1$, $\chi\equiv1$ on $B(0,1)$, and $\chi\equiv0$ on $\R^d\setminus B(0,2)$. For $R>0$, we set $\chi_R(x):=\chi(x/R)$.
\item Given $\lambda>0$, $z \in \mathbb{R}^d$ and a function $W:\R^d\to \R$, we define the rescaled function as
\EQ{
W[z,\lam](x):=\frac{1}{\lambda^{\frac{d-2}{2}}}W\left(\frac{x-z}{\lam}\right).
}
\item For $f\in C_c^\infty(\R^d)$, define
\begin{align}
        \widehat f(\xi)
        :=
        \frac{1}{(2\pi)^{d/2}}
        \int_{\R^d}
        e^{-\iii x\cdot\xi}f(x)\,\ud x.
\end{align}
By a density argument, the above definition extends to functions $f\in\dot H^1(\R^d)$. For $N>0$, define
\begin{align}
        \widehat{P_{\leq N}f}(\xi)
        &:=
        \mathbf{1}_{\{|\xi|\leq N\}}\widehat f(\xi),
        \quad
        P_{>N}f
        :=
        f-P_{\leq N}f,
\end{align}
where $\mathbf{1}_{\{|\xi|\leq N\}}$ is the indicator function of $\{|\xi|\leq N\}$. We use the same notation for the Fourier projections on $\dot H^{-1}(\R^d)$. For $f\in\dot H^{-1}(\R^d)$, define
\begin{align}
        (P_{\leq N}f)[g]
        :=
        f[P_{\leq N}g],
        \qquad
        g\in\dot H^1(\R^d),
\end{align}
and set $P_{>N}f:=f-P_{\leq N}f$.
\end{itemize}
\subsection{Acknowledgments}
The author is grateful to Andrew Lawrie for proposing the problem and for many valuable discussions, Tobias Colding for his constant encouragement and invaluable advice, and Yvan Martel for insightful conversations. This work was partially supported by NSF Grant DMS-2405393 and the Simons Dissertation Fellowship in Mathematics.
\section{Preliminaries}
\subsection{Properties of Stationary Solutions}\label{sec:stationary-sol}
In this section, we recall some standard properties of nonzero finite energy solutions to \eqref{eq:yamabe}. We first prove that the scale and center in Definitions~\ref{def:scale} and~\ref{def:center} are well-defined and then record some consequences of these definitions.

\begin{lem}[Center and scale]  \label{lem:well-defn-scale-center} 
Let $\gamma_0\in(0,\bar E_*/2)$ and let $W\in\dot H^1(\R^d)$ be a non-zero stationary solution. Set
$\lambda(W):=\lambda(W;\gamma_0)$. Then $0<\lambda(W)<\infty$, and there exists a center
$a(W):=a(W;\gamma_0)$ satisfying Definition~\ref{def:center}. Moreover, for every $(b,\mu)\in\R^d\times(0,\infty)$ and any choices of centers for $W$ and
$W[b,\mu]$, we have
\EQ{\label{eq:centscale}
\lambda\left(W[b,\mu]\right) = \lambda(W)\mu, \text{ and }   \left| a\left( W[b,\mu]\right) - b-a(W)\mu \right| \le 2\lambda(W)\mu. 
}
\end{lem} 
\begin{proof}
Since $\bar E(W;B(0,R))\to\bar E(W)$ as $R\to\infty$, the set in Definition~\ref{def:scale} is non-empty, so $\lambda(W)<\infty$. Suppose that $\lambda(W)=0$. For every $n\geq1$, choose $0<r_n\leq1/n$ and $a_n\in\R^d$ such that
\EQ{\label{eq:scale contra}
\bar E(W;B(a_n,r_n))\geq\bar E(W)-\gamma_0.
}
The balls $B(a_n,r_n)$ pairwise intersect, since otherwise $\bar E(W)\geq2\bar E(W)-2\gamma_0$ would imply $\bar E(W)\leq2\gamma_0<\bar E_*$. Thus
$|a_n-a_m|\leq r_n+r_m$, so $a_n\to a_\infty$, as $n\to \infty.$ For all large $n\geq 1$,
$B(a_n,r_n)\subset B(a_\infty,2/n)$, and hence
\[
0<\bar E(W)-\gamma_0
\leq\bar E(W;B(a_\infty,2/n))\to0,
\]
as $n\to \infty$ which is a contradiction. Choose sequences $\lambda_n>0$ and $a_n\in \R^d$ such that
\[
\lambda(W)\leq\lambda_n<\lambda(W)+1/n,
\qquad
\bar E(W;B(a_n,\lambda_n))\geq\bar E(W)-\gamma_0.
\]
By the same argument as above, the balls $B(a_n,\lambda_n)$ pairwise intersect. Thus, the sequence $\{a_n\}_{n\in \N}\subset \R^d$ is bounded, and, after passing to a subsequence, $a_n\to a_\infty$, as $n\to \infty.$ By dominated convergence,
\[
\bar E(W;B(a_\infty,\lambda(W)))
=
\lim_{n\to\infty}\bar E(W;B(a_n,\lambda_n))
\geq\bar E(W)-\gamma_0,
\]
so $a_\infty$ is a center. A change of variables gives
$\lambda(W[b,\mu])=\mu\lambda(W)$. Any two centers of $W$ have intersecting balls of radius $\lambda(W)$, so their distance is at most $2\lambda(W)$. Since ${(a(W[b,\mu])-b)}/{\mu}$ is a center of $W$, \eqref{eq:centscale} follows.
\end{proof}
\begin{lem}[Uniform exterior decay for bounded stationary solutions]
\label{lem:uniform-decay-stationary}
There exists $\gamma_*=\gamma_*(d)>0$ with the following property. Assume that the fixed number $\gamma_0$ used in Definitions~\ref{def:scale} and~\ref{def:center} satisfies
$0<\gamma_0<\min\{\bar E_*/2,\gamma_*\}$. Let $0<C<\infty$. Then there exists a function $\omega_C:[1,\infty)\to[0,\infty)$ with $\omega_C(R)\to0$ as $R\to\infty$ such that the following holds. If $W$ is a nontrivial solution of \eqref{eq:yamabe} with $ \bar E(W)\leq C,$ then
\begin{align}
        \bar E\bigl(W;\R^d\setminus B(a(W),R\lambda(W))\bigr)
        +
        \|W\|_{L^{p+1}(\R^d\setminus B(a(W),R\lambda(W)))}^{p+1}
        \leq \omega_C(R)
\end{align}
for every $R\geq1$.
\end{lem}

\begin{proof}
We write $a(W)=a(W;\gamma_0)$ and $\lambda(W)=\lambda(W;\gamma_0)$. Define
\begin{align}\label{eq:uniform-decay-normalized-profile}
        V(y):=\lambda(W)^{\frac{d-2}{2}}W(a(W)+\lambda(W)y).
\end{align}
Then $V$ solves \eqref{eq:yamabe}, and by scale invariance,
\begin{align}\label{eq:uniform-decay-normalized-bounds}
        \bar E(V)=\bar E(W)\leq C,
        \quad
        \|V\|_{L^{p+1}}=\|W\|_{L^{p+1}}.
\end{align}
Moreover, by the definition of $a(W)$ and $\lambda(W)$,
\begin{align}\label{eq:uniform-decay-normalized-small-tail}
        \bar E(V;\R^d\setminus B(0,1))\leq \gamma_0 .
\end{align}
We will use the exterior Hardy inequality
\begin{align}\label{eq:uniform-decay-exterior-Hardy}
\int_{\R^d\setminus B(0,r)}
\frac{|f|^2}{|x|^2}\,\ud x
\lesssim
\bar E(f;\R^d\setminus B(0,r)),
\qquad
f\in\dot H^1(\R^d),\quad r>0.
\end{align}
To see this, first note that for any $f\in C_c^\infty(\R^d)$, integration by parts gives
\begin{align}
(d-2)\int_{\R^d\setminus B(0,r)}
\frac{|f|^2}{|x|^2}\,\ud x
\leq
2\left(
\int_{\R^d\setminus B(0,r)}
\frac{|f|^2}{|x|^2}\,\ud x
\right)^{1/2}
\left(
\int_{\R^d\setminus B(0,r)}
|\nabla f|^2\,\ud x
\right)^{1/2},
\end{align}
and \eqref{eq:uniform-decay-exterior-Hardy} follows by density. See also proof of Lemma 2.10 in \cite{nlh-radial} for more details. It is enough to prove the desired estimate for $V$, since
\begin{align}\label{eq:uniform-decay-rescale-energy}
        \bar E(W;\R^d\setminus B(a(W),R\lambda(W)))
        =
        \bar E(V;\R^d\setminus B(0,R))
\end{align}
and
\begin{align}\label{eq:uniform-decay-rescale-Lp}
        \|W\|_{L^{p+1}(\R^d\setminus B(a(W),R\lambda(W)))}^{p+1}
        =
        \|V\|_{L^{p+1}(\R^d\setminus B(0,R))}^{p+1}.
\end{align}
Let $\zeta\in C^\infty(\R^d)$ satisfy $0\leq\zeta\leq1$, $\zeta=0$ on $B(0,1)$, $\zeta=1$ on $\R^d\setminus B(0,2)$, and $|\nabla\zeta(x)|\lesssim |x|^{-1}$. By Sobolev's inequality, \eqref{eq:uniform-decay-exterior-Hardy}, and
\eqref{eq:uniform-decay-normalized-small-tail},
\begin{align}\label{eq:uniform-decay-fixed-Lp-tail}
        \|V\|_{L^{p+1}(\R^d\setminus B(0,2))}^{p+1}
        &\lesssim
        \|\nabla(\zeta V)\|_{L^2}^{p+1} \lesssim
        \left(
        \int_{\R^d\setminus B(0,1)}|\nabla V|^2\,\ud x
        +
        \int_{\R^d\setminus B(0,1)}\frac{|V|^2}{|x|^2}\,\ud x
        \right)^{\frac{p+1}{2}} \lesssim
        \gamma_0^{\frac{p+1}{2}}.
\end{align}
Hence
\begin{align}\label{eq:uniform-decay-fixed-Lp-tail-pminus}
        \|V\|_{L^{p+1}(\R^d\setminus B(0,2))}^{p-1}
        \lesssim
        \gamma_0^{\frac{p-1}{2}} .
\end{align}

Let $R>4$ and let $\phi_R\in C^\infty(\R^d)$ be a radial logarithmic cutoff satisfying $0\leq\phi_R\leq1$, $\phi_R=0$ on $B(0,2)$, $\phi_R=1$ on $\R^d\setminus B(0,R)$, and $|\nabla\phi_R(x)| \lesssim \frac{1}{|x|\log R}$. Thus, $\|\nabla\phi_R\|_{L^d} \lesssim (\log R)^{-\frac{d-1}{d}}.$ By Hölder's inequality, $V\phi_R^2\in\dot H^1(\R^d)$. Thus, testing \eqref{eq:yamabe} by a density argument gives
\begin{align}\label{eq:uniform-decay-tested}
        \int_{\R^d}|\nabla V|^2\phi_R^2\,\ud x
        =
        \int_{\R^d}|V|^{p+1}\phi_R^2\,\ud x
        -
        2\int_{\R^d}V\phi_R\nabla V\cdot\nabla\phi_R\,\ud x .
\end{align}
Set
\begin{align}
        X_R:=\int_{\R^d}|\nabla V|^2\phi_R^2\,\ud x,
        \quad
        Y_R:=\int_{\R^d}|V|^2|\nabla\phi_R|^2\,\ud x .
\end{align}
By \eqref{eq:uniform-decay-normalized-small-tail}, $X_R\leq\gamma_0$. Also, since $p+1=\frac{2d}{d-2}$, Hölder's and Sobolev's inequalities, \eqref{eq:uniform-decay-normalized-bounds}, and $\|\nabla\phi_R\|_{L^d} \lesssim (\log R)^{-\frac{d-1}{d}}$ imply
\begin{align}\label{eq:uniform-decay-YR}
        Y_R
        \leq
        \|V\|_{L^{p+1}}^2\|\nabla\phi_R\|_{L^d}^2
        \lesssim_C
        (\log R)^{-\frac{2(d-1)}d}.
\end{align}
For the nonlinear term, Hölder's and Sobolev's inequalities, and \eqref{eq:uniform-decay-fixed-Lp-tail-pminus} give
\begin{align}\label{eq:uniform-decay-nonlinear-term}
        \int_{\R^d}|V|^{p+1}\phi_R^2\,\ud x
        &\leq 
        \|V\phi_R\|_{L^{p+1}}^2
        \|V\|_{L^{p+1}(\R^d\setminus B(0,2))}^{p-1} \lesssim
        \gamma_0^{\frac{p-1}{2}}\|\nabla(V\phi_R)\|_{L^2}^2 \lesssim
        \gamma_0^{\frac{p-1}{2}}(X_R+Y_R).\qquad 
\end{align}
Moreover, by Cauchy--Schwarz and Young's inequality, for every $\veps>0$,
\begin{align}\label{eq:uniform-decay-cross-term}
        \left|
        2\int_{\R^d}V\phi_R\nabla V\cdot\nabla\phi_R\,\ud x
        \right|
        \leq
        \veps X_R+C_\veps Y_R .
\end{align}
Combining \eqref{eq:uniform-decay-tested}, \eqref{eq:uniform-decay-nonlinear-term}, and \eqref{eq:uniform-decay-cross-term}, we obtain
\begin{align}
        X_R
        \lesssim
        \gamma_0^{\frac{p-1}{2}}(X_R+Y_R)
        +
        \veps X_R+C_\veps Y_R .
\end{align}
Choose $\gamma_*=\gamma_*(d)>0$ and then $\veps=\veps(d)>0$ small enough that the terms involving $X_R$ on the right-hand side can be absorbed into the left-hand side for every $0<\gamma_0<\gamma_*$. Hence
\begin{align}\label{eq:uniform-decay-XR}
        X_R
        \lesssim_C
        Y_R
        \lesssim_C
        (\log R)^{-\frac{2(d-1)}d}.
\end{align}
Since $\phi_R\equiv1$ on $\R^d\setminus B(0,R)$, this gives
\begin{align}\label{eq:uniform-decay-energy-normalized}
        \bar E(V;\R^d\setminus B(0,R))
        \lesssim_C
        (\log R)^{-\frac{2(d-1)}d}
        \qquad\textrm{for }R>4.
\end{align}
It remains to obtain an analogous estimate for the $L^{p+1}$-norm. Let $R>8$ and choose $\eta_R\in C^\infty(\R^d)$ satisfying $0\leq\eta_R\leq1$, $\eta_R=0$ on $B(0,R/2)$, $\eta_R=1$ on $\R^d\setminus B(0,R)$, and $|\nabla\eta_R(x)|\lesssim |x|^{-1}$. By Sobolev's inequality, \eqref{eq:uniform-decay-exterior-Hardy}, and
\eqref{eq:uniform-decay-energy-normalized},
\begin{align}\label{eq:uniform-decay-Lp-normalized}
        \|V\|_{L^{p+1}(\R^d\setminus B(0,R))}^{p+1} \lesssim
        \left(
        \bar E(V;\R^d\setminus B(0,R/2))
        +
        \int_{|x|>R/2}\frac{|V|^2}{|x|^2}\,\ud x
        \right)^{\frac{p+1}{2}}
        \lesssim_C
        (\log R)^{-\frac{(d-1)(p+1)}d}.\qquad 
\end{align}
Combining \eqref{eq:uniform-decay-energy-normalized} and \eqref{eq:uniform-decay-Lp-normalized}, and using $p+1>2$, we obtain
\begin{align}\label{eq:uniform-decay-combined-normalized}
        \bar E(V;\R^d\setminus B(0,R))
        +
        \|V\|_{L^{p+1}(\R^d\setminus B(0,R))}^{p+1}
        \lesssim_C
        (\log R)^{-\frac{2(d-1)}d}
        \qquad\textrm{for }R>8.
\end{align}
Rescaling back by \eqref{eq:uniform-decay-rescale-energy} and \eqref{eq:uniform-decay-rescale-Lp} proves the estimate for $R>8$. For $1\leq R\leq8$, Sobolev's inequality and $\bar E(W)\leq C$ give the uniform bound
\begin{align}
        \bar E(W)
        +
        \|W\|_{L^{p+1}}^{p+1}
        \lesssim_C 1 .
\end{align}
Thus, after increasing the constant, we may take
\begin{align}
        \omega_C(R):=
        \begin{cases}
        C_1, & 1\leq R\leq 8,\\
        C_1(\log R)^{-\frac{2(d-1)}d}, & R>8.
        \end{cases}
\end{align}
Here $C_1>0$ depends on $C$, and $\omega_C(R)\to0$ as $R\to\infty$.
\end{proof}

\subsection{Local well-posedness of the nonlinear heat flow}
We first recall the local well-posedness theory for the nonlinear heat equation \eqref{eqn:NLH} and then establish localized energy inequalities that will be used to propagate the proximity of solutions of \eqref{eqn:NLH} to stationary solutions of \eqref{eq:yamabe} for a short time. Local well-posedness is standard, but the short-time propagation results require new arguments because the nonlinear energy \eqref{eqn:energy} is not coercive, so standard energy estimates do not apply directly. We use profile decompositions to overcome this difficulty.
\begin{thm}[Local well-posedness and finite-time concentration]\label{thm:lwp}
Let $d\geq3$ and let $u_0\in\dot H^1(\R^d)$. Then there exists a maximal time
$T_+=T_+(u_0)\in(0,\infty]$ and a unique solution
\begin{align}\label{eq:lwp-solution-class}
        u\in C([0,T_+);\dot H^1(\R^d))
\end{align}
to \eqref{eqn:NLH} with $u(0)=u_0$. If $T_+<\infty$ and $u(t)\to u^*$ strongly in $\dot H^1(\R^d)$ as $t\to T_+$ for some $u^*\in \dot{H}^1(\R^d)$, then $u$ extends past $T_+$. Moreover, for every
$0<T_1<T_2<T_+$,
\begin{align}\label{eq:lwp-positive-time-regularity}
        u\in L^\infty(\R^d\times(T_1,T_2))\cap C^{2,1}(\R^d\times(T_1,T_2)).
\end{align}
Furthermore,
\begin{align}\label{eq:lwp-positive-time-H2}
        \partial_tu\in C((0,T_+);H^1(\R^d)),
        \quad
        D_x^\alpha u\in C((0,T_+);L^2(\R^d))
        \quad\text{for }2\leq |\alpha|\leq3 .
\end{align}
The nonlinear energy $E(u(t))$ is continuous and nonincreasing on
$[0,T_+)$ and in particular for every $0< t_1\leq t_2<T_+$,
\begin{align}\label{eq:energy-identity}
        E(u(t_2))
        +
        \int_{t_1}^{t_2}\|\partial_tu(t)\|_{L^2}^2\,\ud t
        =
        E(u(t_1)).
\end{align}
Assume in addition that
\begin{align}\label{eqn:limsup bound}
        \sup_{0\leq t<T_+}\bar E(u(t))\leq M<\infty .
\end{align}
Then
\begin{align}\label{eq:tension-L2}
        \int_0^{T_+}\|\partial_tu(t)\|_{L^2}^2\,\ud t
        \lesssim
        \sup_{0\leq t<T_+}
        \left(\bar E(u(t))+\bar E(u(t))^{\frac{p+1}{2}}\right)
        <\infty .
\end{align}
If $T_+<\infty$, then there exist $\veps_0>0$, $R_0>0$, an integer $L\geq 0$ and a finite set $\calS\subset\R^d$ with cardinality $L$ such that, for every $x\in \calS$ and every $0<R\leq R_0$,
\begin{align}\label{defn:singular-points}
        \limsup_{t\to T_+}
        \left[
        \bar E(u(t);B(x,R))
        +
        \|u(t)\|_{L^{p+1}(B(x,R))}^{p+1}
        \right]
        \geq \veps_0 .
\end{align}
Moreover, there exists $u^*\in\dot H^1(\R^d)$ such that $u(t)\rightharpoonup u^*$ weakly in $\dot{H}^1(\R^d)$ as $t\to T_+$ and $u^*$ has a $C^1$ representative on $\R^d\setminus\calS$ such that, for every compact set $K\Subset\R^d\setminus\calS$,
\begin{align}
\|u(t)-u^*\|_{C^1(K)}\to0
\qquad\textrm{as }t\to T_+.
\end{align}
\end{thm}

\begin{proof}
The local well-posedness theory for \eqref{eqn:NLH} is due to \cite{brezis-cazenave}. For our purposes, we use Proposition~2.2 of \cite{harada2024dynamics}, which also establishes \eqref{eq:lwp-positive-time-regularity} and \eqref{eq:lwp-positive-time-H2}. To prove the continuation statement, suppose that $T_+<\infty$ and that $u(t)\to u^*$ strongly in $\dot H^1(\R^d)$ as $t\to T_+$. Fix $0<\delta<T_+(u^*)$. By Proposition~2.1(iii) in \cite{ikeda2024global}, $T_+(u(t))>\delta$ for all $t>0$ sufficiently close to $T_+$. Choose such a time $t>0$ with $T_+-t<\delta$. By time translation, the solution with initial data $u(t)$ exists on $[t,t+\delta)$, and by uniqueness it agrees with $u$ on $[t,T_+)$. Since $t+\delta>T_+$, this extends $u$ past $T_+$.

The identity \eqref{eq:energy-identity} follows by differentiating the nonlinear energy $E$ and integrating by parts (cf. Proposition~2.2 in \cite{ikeda2024global} for details). The estimate \eqref{eq:tension-L2} follows from the Sobolev inequality and \eqref{eqn:limsup bound}.

Assume now that $T_+<\infty$. Since $u$ is a classical solution for positive times, it is a suitable weak solution in the sense of \cite{wang-wei}. Let $\veps_1,\delta_*,c_*,\theta_*>0$ be the constants in Theorem~3.9 of \cite{wang-wei}. Choose $C_1>1+\theta_*+\delta_*^2$ and then $R_0>0$ so that $T_+-C_1R_0^2>0$. Define
\begin{align}\label{eq:lwp-calS-def}
        \calS:=
        \{
        x\in\R^d:
        \int_{T_+-C_1r^2}^{T_+}\int_{B(x,r)}
        \left(|\nabla u|^2+|u|^{p+1}\right)\ud y\ud s
        \geq
        \veps_1r^2
        \quad\textrm{for every }0<r\leq R_0
        \}.\quad 
\end{align}
If $x\notin\calS$, then the reverse inequality in \eqref{eq:lwp-calS-def} holds for some $0<r\leq R_0$. Apply Theorem~3.9 of \cite{wang-wei} at $(x,T_+-\delta_*^2r^2)$. The hypothesis of Theorem~3.9 in \cite{wang-wei} at $(x,T_+-\delta_*^2r^2)$ holds because
\begin{align}
B(x,r)\times
\left(T_+-(1+\theta_*+\delta_*^2)r^2,
T_+-(\theta_*+\delta_*^2)r^2\right) \subset
B(x,r)\times(T_+-C_1r^2,T_+).
\end{align}
Hence, using the same notation as in \cite{wang-wei} we obtain
\begin{align}
\sup_{B(x,\delta_*r)\times
(T_+-2\delta_*^2r^2,T_+)}
|u|
\leq c_*r^{-\frac{d-2}{2}}.
\end{align}
Then standard interior estimates give uniform spatial $C^{2,\beta}$ bounds, for some $\beta\in(0,1)$, on a smaller parabolic cylinder. 

Now we show \eqref{defn:singular-points}. Fix $x\in\calS$ and $0<R\leq R_0$. Choose $r_j \to 0$ as $j\to \infty$ with $r_j<R$ for all $j\in \N$. By \eqref{eq:lwp-calS-def}, for each $j\in \N$ there exists $s_j\in(T_+-C_1r_j^2,T_+)$ such that
\begin{align}\label{eq:lwp-time-slice-lower}
        \int_{B(x,r_j)}
        \left(|\nabla u(s_j)|^2+|u(s_j)|^{p+1}\right)\ud y
        \geq
        \frac{\veps_1}{C_1}.
\end{align}
Since $s_j\to T_+$ as $j\to \infty$ and $B(x,r_j)\subset B(x,R)$, we obtain
\begin{align}\label{eq:lwp-terminal-limsup-lower}
        \limsup_{t\to T_+}
        \left[
        \bar E(u(t);B(x,R))
        +
        \|u(t)\|_{L^{p+1}(B(x,R))}^{p+1}
        \right]
        \geq
        \frac{\veps_1}{C_1}.
\end{align}
Thus \eqref{defn:singular-points} holds after setting $\veps_0:=\veps_1/C_1$.

We next prove that $\calS$ is finite. Let $x_1,\ldots,x_N\in\calS$ be distinct points for $N\geq 1$. Choose $0<r\leq R_0$ such that the balls $B(x_i,r)$, $1\leq i\leq N$, are pairwise disjoint. Summing \eqref{eq:lwp-calS-def} gives
\begin{align}\label{eq:lwp-finiteness}
        N\veps_1r^2
        &\leq
        \int_{T_+-C_1r^2}^{T_+}
        \int_{\cup_{i=1}^N B(x_i,r)}
        \left(|\nabla u|^2+|u|^{p+1}\right)\ud y\ud s \notag\\
        &\leq
        C_1r^2
        \sup_{0\leq t<T_+}
        \left(\bar E(u(t))+\|u(t)\|_{L^{p+1}}^{p+1}\right) \notag\\
        &\lesssim
        r^2\left(M+M^{\frac{p+1}{2}}\right),
\end{align}
where we used Sobolev's inequality and \eqref{eqn:limsup bound}. This proves that $\calS$ is finite.

Finally, we prove the existence of the weak limit. Since $u(t)$ is bounded in $\dot H^1(\R^d)$, every sequence $t_n\to T_+$ as $n\to \infty$ has a subsequence converging weakly in $\dot H^1$. To see that the weak limit is unique, let $\zeta\in C_c^\infty(\R^d)$. For $0<t_1<t_2<T_+$,
$$\left|\int_{\R^d}(u(t_2)-u(t_1))\zeta\,\ud x\right|
\leq \|\zeta\|_{L^2}\int_{t_1}^{t_2}\|\partial_tu(t)\|_{L^2}\,\ud t.$$
By \eqref{eq:tension-L2}, the right-hand side tends to $0$ as $t_1,t_2\to T_+$. Hence $\int_{\R^d} u(t)\zeta\,\ud x$ has a limit for every $\zeta\in C_c^\infty(\R^d)$. This distributional limit agrees with every weak $\dot H^1$ subsequential limit, and therefore $u(t)\rightharpoonup u^*$ weakly in $\dot H^1(\R^d)$ as $t\to T_+$. The preceding bound and standard interior estimates give uniform spatial $C^{2,\beta}$ bounds, for some $\beta\in(0,1)$, on a smaller parabolic neighborhood of every point in $\R^d\setminus\calS$. A finite covering of $K$, Arzelà--Ascoli, and the distributional convergence to $u^*$ now give convergence in $C^1(K)$.
\end{proof}
\begin{rem}
The non-emptiness of $\calS$ when $T_+<\infty$ and $\sup_{t<T_+}\bar E(u(t))<\infty$ will be proved in Lemma~\ref{lem:finite-time-no-escape}.
\end{rem}
We next collect localized energy inequalities that will be needed later in our arguments.
\begin{lem}[Localized energy identities]\label{lem:localized-energy-density}
Let $u(t)$ be the solution in Theorem~\ref{thm:lwp}. Let $\phi\in C_c^\infty(\R^d)$ and
let $0<t_1<t_2<T_+$. Then
\begin{align}
        E_{\phi}(u(t_2))-E_\phi(u(t_1))
        =
        -\int_{t_1}^{t_2}\int_{\R^d}(\partial_t u)^2\phi^2\,\ud x\ud t
        -2\int_{t_1}^{t_2}\int_{\R^d}
        (\nabla u\cdot\nabla\phi)\phi\partial_tu\,\ud x\ud t,
        \label{eqn:loc energy equality}
\end{align}
and consequently
\begin{align}
        \left|E_{\phi}(u(t_2))-E_\phi(u(t_1))\right| &\leq \int_{t_1}^{t_2}\int_{\R^d}(\partial_tu)^2\phi^2\,\ud x\ud t\\
        &\quad +
        2\left(
        \int_{t_1}^{t_2}\int_{\R^d}(\partial_tu)^2\phi^2\,\ud x\ud t
        \right)^{1/2}
        \left(
        \int_{t_1}^{t_2}\int_{\R^d}|\nabla u|^2|\nabla\phi|^2\,\ud x\ud t
        \right)^{1/2}.\qquad  
        \label{eqn:energy ineq I}
\end{align}
\end{lem}
\begin{proof}
Differentiating $E_\phi(u(t))$ and integrating by parts gives
\begin{align}\label{eq:localized-energy-proof}
        \frac{\ud}{\ud t}E_\phi(u(t))
        &=
        \int_{\R^d}
        \left(\nabla u\cdot\nabla\partial_tu-|u|^{p-1}u\partial_tu\right)\phi^2\,\ud x \notag\\
        &=
        -\int_{\R^d}
        \left(\Delta u+|u|^{p-1}u\right)\partial_tu\,\phi^2\,\ud x
        -2\int_{\R^d}(\nabla u\cdot\nabla\phi)\phi\partial_tu\,\ud x \notag\\
        &=
        -\int_{\R^d}(\partial_tu)^2\phi^2\,\ud x
        -2\int_{\R^d}(\nabla u\cdot\nabla\phi)\phi\partial_tu\,\ud x .
\end{align}
% Integrating \eqref{eq:localized-energy-proof} proves \eqref{eqn:loc energy equality}. Similarly,
% \begin{align}\label{eq:localized-modified-energy-proof}
%         \frac{\ud}{\ud t}\bar E_\phi(u(t))
%         &=
%         2\int_{\R^d}\nabla u\cdot\nabla\partial_tu\,\phi^2\,\ud x \notag\\
%         &=
%         -2\int_{\R^d}\Delta u\,\partial_tu\,\phi^2\,\ud x
%         -4\int_{\R^d}(\nabla u\cdot\nabla\phi)\phi\partial_tu\,\ud x \notag\\
%         &=
%         -2\int_{\R^d}(\partial_tu)^2\phi^2\,\ud x
%         +2\int_{\R^d}|u|^{p-1}u\partial_tu\,\phi^2\,\ud x \notag\\
%         &\quad
%         -4\int_{\R^d}(\nabla u\cdot\nabla\phi)\phi\partial_tu\,\ud x .
% \end{align}
The estimate \eqref{eqn:energy ineq I} follows from the Cauchy--Schwarz inequality.
\end{proof}
\begin{rem}\label{rem:localized-energy-noncompact-cutoff}
Note that the conclusions of Lemma~\ref{lem:localized-energy-density} remain valid for cutoff functions $\phi\in C^\infty(\R^d)\cap W^{1,\infty}(\R^d)$ satisfying $0\leq\phi\leq1$, provided the quantities in the right-hand side of \eqref{eqn:energy ineq I} are finite. 
\end{rem}
As a consequence of these localized energy inequalities, we can propagate localized nonlinear energy over short time intervals.
\begin{lem}[Short-time propagation of localized nonlinear energy]
\label{lem:short-time-nonlinear-energy-prop}
Let $u(t)$ be a solution to \eqref{eqn:NLH} with initial data $u(0)=u_0\in\dot H^1$. Let
$T_+=T_+(u_0)$ denote its maximal time of existence and assume that
\begin{align}\label{eq:short-time-nonlinear-energy-bound}
        \sup_{t\in[0,T_+)}\|u(t)\|_{\dot H^1}<\infty .
\end{align}
Let $0<\sigma_n<\ta_n<T_+$ satisfy $\sigma_n,\ta_n\to T_+$ and
$h_n:=\ta_n-\sigma_n\to0$ as $n\to\infty$. If $T_+=\infty$, this means
$\sigma_n,\ta_n\to\infty$. Let $\eta_n\in C_c^\infty(\R^d)$ satisfy
$0\leq \eta_n\leq1$, $\|\nabla\eta_n\|_{L^\infty}\lesssim \ell_n^{-1}$ for some $\ell_n>0$ and
$h_n/\ell_n^2\to0$ as $n\to\infty$. Then
\begin{align}\label{eq:short-time-nonlinear-energy-conclusion}
        \lim_{n\to \infty} |E_{\eta_n}(u(\ta_n))-E_{\eta_n}(u(\sigma_n))|=0.
\end{align}
Consequently, if $W\in\dot H^1(\R^d)$ is a stationary solution of \eqref{eq:yamabe}, possibly zero, and
\begin{align}\label{eq:short-time-nonlinear-energy-initial}
        \lim_{n\to \infty} |E_{\eta_n}(u(\sigma_n))-E_{\eta_n}(W)| = 0,
\end{align}
then
\begin{align}\label{eq:short-time-nonlinear-energy-final}
       \lim_{n\to \infty} |E_{\eta_n}(u(\ta_n))-E_{\eta_n}(W)| = 0.
\end{align}
\end{lem}
\begin{proof}
Since $\sigma_n>0$, Lemma~\ref{lem:localized-energy-density} applies on
$[\sigma_n,\ta_n]$. From \eqref{eqn:energy ineq I},
\begin{align}\label{eq:short-time-nonlinear-energy-est}
        \left|E_{\eta_n}(u(\ta_n))-E_{\eta_n}(u(\sigma_n))\right|
        &\leq
        \int_{\sigma_n}^{\ta_n}\int_{\R^d}|\partial_tu|^2\eta_n^2\,\ud x\ud t \notag\\
        &\quad+
        2\left(
        \int_{\sigma_n}^{\ta_n}\int_{\R^d}|\partial_tu|^2\eta_n^2\,\ud x\ud t
        \right)^{1/2}
        \left(
        \int_{\sigma_n}^{\ta_n}\int_{\R^d}|\nabla u|^2|\nabla\eta_n|^2\,\ud x\ud t
        \right)^{1/2}.
\end{align}
By \eqref{eq:tension-L2},
\begin{align}\label{eq:short-time-tension-tail}
        \int_{\sigma_n}^{\ta_n}\|\partial_tu(t)\|_{L^2}^2\,\ud t\to0
        \qquad\textrm{as }n\to\infty .
\end{align}
Also,
\begin{align}\label{eq:short-time-gradient-cutoff}
        \int_{\sigma_n}^{\ta_n}\int_{\R^d}|\nabla u|^2|\nabla\eta_n|^2\,\ud x\ud t
        \lesssim
        \frac{h_n}{\ell_n^2}
        \sup_{t\in[0,T_+)}\|u(t)\|_{\dot H^1}^2
        \to0
        \qquad\textrm{as }n\to\infty .
\end{align}
Combining \eqref{eq:short-time-nonlinear-energy-est},
\eqref{eq:short-time-tension-tail}, and
\eqref{eq:short-time-gradient-cutoff} proves
\eqref{eq:short-time-nonlinear-energy-conclusion}. The final implication follows from
\eqref{eq:short-time-nonlinear-energy-initial}, \eqref{eq:short-time-nonlinear-energy-conclusion} and the triangle inequality.
\end{proof}
We next recall some facts about profile decompositions (cf. Theorem 1.1 in
\cite{gerard1998description} and Proposition 3.4 in
\cite{ishiwata2018potential}).
\begin{prop}\label{prop:profile-decomp}
Let $\{u_n\}_{n\in\N}$ be bounded in $\dot H^1(\R^d)$. Then, after passing to a subsequence still denoted by $u_n$, there exist $\{\lambda_n^i\}_{n\in\N}\subset(0,\infty)$, $\{a_n^i\}_{n\in\N}\subset\R^d$, and $\psi^i\in\dot H^1(\R^d)$, $i\in\N$, satisfying the following:

$\mathrm{(a)}$ There holds
\begin{align}
\frac{\lambda_n^i}{\lambda_n^j}
+
\frac{\lambda_n^j}{\lambda_n^i}
+ \frac{|a_n^j-a_n^i|}{\lambda^i_n}
\to \infty
\quad \textrm{as } n\to\infty
\quad \textrm{for } i\neq j .
\end{align}

$\mathrm{(b)}$ For $J\in \mathbb N$, let
\begin{align}
r_n^J(x):=
u_n(x)
-
\sum_{i=1}^J
(\lambda_n^i)^{-\frac{d-2}{2}}
\psi^i\left(\frac{x-a_n^i}{\lambda_n^i}\right).
\label{eq:profile-decomposition-remainder}
\end{align}
Then, for every fixed $J\in\N$, we have
\begin{align}
\|\nabla u_n\|_{L^2}^2
&=
\sum_{i=1}^J \|\nabla \psi^i\|_{L^2}^2
+
\|\nabla r_n^J\|_{L^2}^2
+
o_n(1),\\
\|u_n\|_{L^{p+1}}^{p+1}
&=
\sum_{i=1}^J\|\psi^i\|_{L^{p+1}}^{p+1}
+
\|r_n^J\|_{L^{p+1}}^{p+1}
+
o_n(1),
\end{align}
as $n\to\infty$, and
\begin{align}
\lim_{J\to\infty}\limsup_{n\to\infty}\|r_n^J\|_{L^{p+1}}=0.
\end{align} 
Moreover, we may choose the decomposition so that, for every fixed $J$ and every $1\leq j\leq J$,
\begin{align}\label{eq:profile-remainder-frame-vanishing}
(\lambda_n^j)^{\frac{d-2}{2}}r_n^J(a_n^j+\lambda_n^j\cdot)
\rightharpoonup0
\qquad\textrm{weakly in }\dot H^1(\R^d).
\end{align}
\end{prop}
Since we work with cutoff functions, we recall the following lemma, which gives a Pythagorean identity for profile decompositions with cutoff functions and a lower bound for the nonlinear energy of such a decomposition when the profiles are stationary solutions of \eqref{eq:yamabe}. 
\begin{lem}\label{lem:weighted-profile-decoupling}
Let $\{v_n\}_{n\in \N}$ be bounded in $\dot H^1(\R^d)$ and admit the profile decomposition as in Proposition~\ref{prop:profile-decomp}. For every fixed $J$, write
\begin{align}\label{eq:weighted-profile-decomp}
        v_n=\sum_{j=1}^J \psi_n^j+r_n^J,
        \quad
        \psi_n^j(x):=(\lambda_n^j)^{-\frac{d-2}{2}}
        \psi^j\left(\frac{x-a_n^j}{\lambda_n^j}\right).
\end{align}
Let $\phi_n\in C^\infty(\R^d)$ satisfy $0\leq\phi_n\leq1$ and
$\sup_n\|\nabla\phi_n\|_{L^d}<\infty$. Then, for every fixed $J$,
\begin{align}\label{eq:weighted-BL-Lq}
        \int_{\R^d}|v_n|^{p+1}\phi_n^2\,\ud x
        =
        \sum_{j=1}^J\int_{\R^d}|\psi_n^j|^{p+1}\phi_n^2\,\ud x
        +
        \int_{\R^d}|r_n^J|^{p+1}\phi_n^2\,\ud x
        +
        o_n(1).
\end{align}
Moreover, assume that $\psi^1,\ldots,\psi^J$ are stationary solutions of \eqref{eq:yamabe} and that
\begin{align}\label{eqn:weighted-rem-small}
        \|r_n^J\|_{L^{p+1}(\supp\phi_n)}\to0,
        \textrm{ as }n\to\infty .
\end{align}
Then
\begin{align}\label{eqn:weighted-Lp-decoupling}
        \int_{\R^d}|\sum_{j=1}^J\psi_n^j|^{p+1}\phi_n^2\,\ud x
        =
        \sum_{j=1}^J\int_{\R^d}|\psi_n^j|^{p+1}\phi_n^2\,\ud x+o_n(1),
\end{align}
\begin{align}\label{eqn:weighted-grad-cross-decoupling}
        \int_{\R^d}\nabla \psi_n^j\cdot\nabla \psi_n^k\,\phi_n^2\,\ud x=o_n(1)
        \qquad\textrm{for }j\neq k,
\end{align}
and
\begin{align}\label{eqn:weighted-rem-cross-decoupling}
        \int_{\R^d}\nabla \psi_n^j\cdot\nabla r_n^J\,\phi_n^2\,\ud x=o_n(1)
        \qquad\textrm{for }1\leq j\leq J.
\end{align}
Consequently, if $S_n:=\sum_{j=1}^J\psi_n^j$, then
\begin{align}\label{eqn:weighted-energy-decoupling}
        E_{\phi_n}(S_n+r_n^J)
        \geq
        \sum_{j=1}^J E_{\phi_n}(\psi_n^j)-o_n(1).
\end{align}
\end{lem}

\begin{proof}
We first record the standard consequence of parameter orthogonality. For $j\neq k$,
\begin{align}\label{eqn:weighted-product-orth}
        \big\||\psi_n^k||\nabla \psi_n^j|\big\|_{L^{\frac{d}{d-1}}}
        +
        \int_{\R^d}|\psi_n^j|^p|\psi_n^k|\,\ud x
        \to0
        \qquad\textrm{as }n\to\infty .
\end{align}
Indeed, \eqref{eqn:weighted-product-orth} is first checked when $\psi^j,\psi^k\in C_c^\infty(\R^d)$, where it follows directly from the orthogonality condition in part $\mathrm{(a)}$ of Proposition~\ref{prop:profile-decomp}. The general case follows from an approximation argument. Next, we will use the following pointwise estimate
\begin{align}\label{eq:weighted-dec-pointwise-direct}
        |
        |\sum_{j=1}^J b_j+c|^{p+1}
        - \sum_{j=1}^J|b_j|^{p+1}
        -
        |c|^{p+1}|
        \lesssim_J \sum_{i\neq j}|b_i|^p|b_j| + \sum_{j=1}^J\left(|b_j|^p|c|+|c|^p|b_j|\right)
\end{align}
with $b_j=\psi_n^j(x)$ and $c=r_n^J(x)$. Since $0\leq \phi_n\leq1$, it is enough to show
\begin{align}\label{eq:weighted-profile-profile-interaction}
        \int_{\R^d}|\psi_n^i|^p|\psi_n^j|\,\ud x\to0
        \qquad\textrm{as }n\to\infty,
        \qquad\textrm{for }i\neq j,
\end{align}
and
\begin{align}\label{eq:weighted-profile-remainder-interaction}
        \int_{\R^d}|\psi_n^j|^p|r_n^J|\,\ud x
        +
        \int_{\R^d}|r_n^J|^p|\psi_n^j|\,\ud x
        \to0
        \qquad\textrm{as }n\to\infty,
        \qquad\textrm{for }1\leq j\leq J.
\end{align}

The limit \eqref{eq:weighted-profile-profile-interaction} follows from \eqref{eqn:weighted-product-orth}. To prove \eqref{eq:weighted-profile-remainder-interaction}, fix $j\leq J$ and set
\begin{align}\label{eq:weighted-rescaled-remainder-vanishing}
        \widetilde r_n^j(y):=
        (\lambda_n^j)^{\frac{d-2}{2}}r_n^J(a_n^j+\lambda_n^j y).
\end{align}
By \eqref{eq:profile-remainder-frame-vanishing}, $\widetilde r_n^j\rightharpoonup0$ weakly in $\dot H^1(\R^d)$, as $n\to \infty.$ In particular, $\widetilde r_n^j$ is bounded in $L^{p+1}(\R^d)$ and converges weakly to $0$ in $L^{p+1}(\R^d)$. We claim that, as $n\to \infty$,
\begin{align}\label{eq:weighted-rescaled-remainder-local-strong}
        \widetilde r_n^j\to0
        \qquad\textrm{strongly in }L^r_{\loc}(\R^d)
        \qquad\textrm{for every }1\leq r<p+1.
\end{align}
Indeed, on each compact set $K\subset\R^d$, the sequence $\widetilde r_n^j$ is bounded in $H^1(K)$, since it is bounded in $\dot H^1(\R^d)$ and in $L^{p+1}(\R^d)$. Hence Rellich's theorem gives strong convergence in $L^2(K)$ along subsequences. Since the weak $L^{p+1}$ limit is $0$, every such $L^2(K)$ limit is $0$, and therefore $\widetilde r_n^j\to0$ strongly in $L^2(K)$, as $n\to \infty.$ Interpolating with the uniform $L^{p+1}$ bound gives strong convergence in $L^r(K)$ for $2\leq r<p+1$, and the case $1\leq r<2$ follows from Hölder's inequality and the fact that $K$ is compact.

Changing variables $x=a_n^j+\lambda_n^j y$, we have
\begin{align}\label{eq:weighted-rescaled-interactions}
        \int_{\R^d}|\psi_n^j|^p|r_n^J|\,\ud x
        +
        \int_{\R^d}|r_n^J|^p|\psi_n^j|\,\ud x
        =
        \int_{\R^d}|\psi^j|^p|\widetilde r_n^j|\,\ud y
        +
        \int_{\R^d}|\widetilde r_n^j|^p|\psi^j|\,\ud y .
\end{align}
Let $\varepsilon>0$. Choose $\zeta\in C_c^\infty(\R^d)$ such that $\|\psi^j-\zeta\|_{L^{p+1}}<\varepsilon$. Using $|\psi^j|^p\lesssim|\zeta|^p+|\psi^j-\zeta|^p$, Hölder's inequality, the uniform $L^{p+1}$ bound on $\widetilde r_n^j$, and \eqref{eq:weighted-rescaled-remainder-local-strong}, we obtain
\begin{align}
        \limsup_{n\to\infty}
        \int_{\R^d}|\psi^j|^p|\widetilde r_n^j|\,\ud y
        \lesssim
        \limsup_{n\to\infty}
        \int_{\R^d}|\zeta|^p|\widetilde r_n^j|\,\ud y
        +
        \limsup_{n\to\infty}
        \|\psi^j-\zeta\|_{L^{p+1}}^p
        \|\widetilde r_n^j\|_{L^{p+1}}
        \lesssim
        \varepsilon^p .
\end{align}
Similarly, using $|\psi^j|\leq|\zeta|+|\psi^j-\zeta|$, Hölder's inequality, the uniform $L^{p+1}$ bound on $\widetilde r_n^j$, and \eqref{eq:weighted-rescaled-remainder-local-strong} with $r=p$, we obtain
\begin{align}
        \limsup_{n\to\infty}
        \int_{\R^d}|\widetilde r_n^j|^p|\psi^j|\,\ud y
        \leq
        \limsup_{n\to\infty}
        \|\zeta\|_{L^\infty}
        \|\widetilde r_n^j\|_{L^p(\supp\zeta)}^p
        +
        \limsup_{n\to\infty}
        \|\widetilde r_n^j\|_{L^{p+1}}^p
        \|\psi^j-\zeta\|_{L^{p+1}}
        \lesssim
        \varepsilon .
\end{align}
Letting $\varepsilon\to0$ in the last two estimates and using \eqref{eq:weighted-rescaled-interactions}, we obtain \eqref{eq:weighted-profile-remainder-interaction}. Therefore \eqref{eq:weighted-BL-Lq} holds. 

To see \eqref{eqn:weighted-Lp-decoupling} we use the following pointwise estimate
\begin{align}\label{eqn:weighted-BL-pointwise}
        |
        |\sum_{j=1}^J b_j|^{p+1}
        -
        \sum_{j=1}^J|b_j|^{p+1}
        |
        \lesssim
        \sum_{j\neq k}|b_j|^p|b_k|,
\end{align}
with $b_j=\psi^j_n$. Since \eqref{eqn:weighted-product-orth} holds and $0\leq\phi_n\leq1$, the error terms on the right-hand side of \eqref{eqn:weighted-BL-pointwise} are $o_n(1)$, which proves \eqref{eqn:weighted-Lp-decoupling}. We next prove \eqref{eqn:weighted-grad-cross-decoupling}. Since $\psi_n^j$ solves \eqref{eq:yamabe}, we have
\begin{align}\label{eqn:weighted-grad-cross-stationary}
        \int_{\R^d}\nabla \psi_n^j\cdot\nabla \psi_n^k\,\phi_n^2\,\ud x
        &=
        \int_{\R^d}|\psi_n^j|^{p-1}\psi_n^j\psi_n^k\phi_n^2\,\ud x
        -2\int_{\R^d}\psi_n^k\phi_n\,\nabla \psi_n^j\cdot\nabla\phi_n\,\ud x .
\end{align}
The first term on the right-hand side of \eqref{eqn:weighted-grad-cross-stationary} tends to $0$ as $n\to\infty$ by \eqref{eqn:weighted-product-orth}. For the second term on the right-hand side of \eqref{eqn:weighted-grad-cross-stationary}, using Hölder's inequality we obtain
\begin{align}
        \left|
        \int_{\R^d}\psi_n^k\phi_n\,\nabla \psi_n^j\cdot\nabla\phi_n\,\ud x
        \right|
        \leq
        \big\||\psi_n^k||\nabla \psi_n^j|\big\|_{L^{\frac{d}{d-1}}}
        \|\nabla\phi_n\|_{L^d}
        =
        o_n(1).
\end{align}
Thus \eqref{eqn:weighted-grad-cross-decoupling} holds. Similarly, using \eqref{eq:yamabe} we obtain
\begin{align}\label{eqn:weighted-rem-cross-stationary}
        \int_{\R^d}\nabla \psi_n^j\cdot\nabla r_n^J\,\phi_n^2\,\ud x
        &=
        \int_{\R^d}|\psi_n^j|^{p-1}\psi_n^jr_n^J\phi_n^2\,\ud x
        -2\int_{\R^d}r_n^J\phi_n\,\nabla \psi_n^j\cdot\nabla\phi_n\,\ud x .
\end{align}
Using \eqref{eqn:weighted-rem-small}, Hölder's inequality, and $\sup_n\|\nabla\phi_n\|_{L^d}<\infty$, we obtain
\begin{align}
        \left|
        \int_{\R^d}\nabla \psi_n^j\cdot\nabla r_n^J\,\phi_n^2\,\ud x
        \right|
        &\lesssim
        \|r_n^J\|_{L^{p+1}(\supp\phi_n)}
        \left(
        \|\psi^j\|_{L^{p+1}}^p
        +
        \|\nabla \psi^j\|_{L^2}\|\nabla\phi_n\|_{L^d}
        \right)
        =
        o_n(1).
\end{align}
This proves \eqref{eqn:weighted-rem-cross-decoupling}. Finally, by \eqref{eqn:weighted-rem-cross-decoupling},
\begin{align}
        \int_{\R^d}|\nabla S_n+\nabla r_n^J|^2\phi_n^2\,\ud x
        \geq
        \int_{\R^d}|\nabla S_n|^2\phi_n^2\,\ud x-o_n(1).
\end{align}
Moreover, since $S_n$ is bounded in $L^{p+1}$,
\begin{align}
        \int_{\R^d}\left||S_n+r_n^J|^{p+1}-|S_n|^{p+1}\right|\phi_n^2\,\ud x
        \lesssim
        \|r_n^J\|_{L^{p+1}(\supp\phi_n)}
        +
        \|r_n^J\|_{L^{p+1}(\supp\phi_n)}^{p+1}
        =
        o_n(1).
\end{align}
Expanding the energy, using \eqref{eqn:weighted-Lp-decoupling} and \eqref{eqn:weighted-grad-cross-decoupling}, we obtain \eqref{eqn:weighted-energy-decoupling}.
\end{proof}

The next lemma is the main place where the profile decomposition is used. It allows us to upgrade short-time control of localized nonlinear energy into propagation of localized $L^{p+1}$-norm.
\begin{lem}[Short-time propagation of the $L^{p+1}$ norm]
\label{lem:short-time-Lq-prop}
Let $u(t)$ be a solution to \eqref{eqn:NLH} with initial data $u(0)=u_0\in\dot H^1$. Let
$T_+=T_+(u_0)$ denote its maximal time of existence and assume
\begin{align}
        \sup_{t\in[0,T_+)}\|u(t)\|_{\dot H^1}<\infty .
\end{align}
Let $0<\sigma_n<\ta_n<T_+$ satisfy $\sigma_n,\ta_n\to T_+$ and
$h_n:=\ta_n-\sigma_n\to0$ as $n\to\infty$. Let $W\in\dot H^1(\R^d)$ be a stationary solution of \eqref{eq:yamabe}, possibly trivial. Let
$\Omega_n^1\subset\Omega_n^2$ be open sets. Let
$\phi_n\in C_c^\infty(\Omega_n^2)$ satisfy
\begin{align}\label{eq:Lq-prop-cutoff}
        0\leq\phi_n\leq1,\qquad
        \phi_n\equiv1 \textrm{ on } \Omega_n^1,\qquad
        \|\nabla\phi_n\|_{L^\infty}\lesssim \ell_n^{-1},\qquad
        \|\nabla\phi_n\|_{L^d}\lesssim1,
\end{align}
for some $\ell_n>0$ with $h_n/\ell_n^2\to0$ as $n\to\infty$. Assume
\begin{align}\label{eq:Lq-prop-initial-small}
        \lim_{n\to \infty}\|u(\sigma_n)-W\|_{L^{p+1}(\Omega_n^2)}=0
\end{align}
and
\begin{align}\label{eq:Lq-prop-energy-input}
        \lim_{n\to \infty}|E_{\phi_n}(u(\sigma_n))-E_{\phi_n}(W)|=0.
\end{align}
Then
\begin{align}\label{eq:Lq-prop-conclusion}
        \lim_{n\to \infty}\|u(\ta_n)-W\|_{L^{p+1}(\Omega_n^1)}=0.
\end{align}
\end{lem}
\begin{proof}
Set $v_n:=u(\ta_n)-W$. Suppose that \eqref{eq:Lq-prop-conclusion} fails. Passing to a subsequence, there exists
$\delta_0>0$ such that
\begin{align}\label{eq:Lq-weighted-lower-vn}
        \int_{\R^d}|v_n|^{p+1}\phi_n^2\,\ud x\geq\delta_0 .
\end{align}
Then applying Proposition \ref{prop:profile-decomp} to the sequence $\{v_n\}_{n\in \N}$ yields for every fixed $J$,
\begin{align}\label{eq:profile-decomp-vn}
        v_n=\sum_{j=1}^J \psi_n^j+r_n^J,
        \quad
        \psi_n^j(x):=(\lambda_n^j)^{-\frac{d-2}{2}}
        \psi^j\left(\frac{x-a_n^j}{\lambda_n^j}\right),
\end{align}
where the parameters are pairwise orthogonal as in part (a) in Proposition \ref{prop:profile-decomp} and
\begin{align}\label{eq:profile-rem-Lq}
        \lim_{J\to\infty}\limsup_{n\to\infty}\|r_n^J\|_{L^{p+1}}=0 .
\end{align}
After a diagonal extraction, the limits
\begin{align}\label{eq:profile-weighted-mass}
        m_j:=
        \lim_{n\to\infty}
        \int_{\R^d}|\psi_n^j|^{p+1}\phi_n^2\,\ud x
\end{align}
exist. By Lemma~\ref{lem:weighted-profile-decoupling},
\eqref{eq:Lq-weighted-lower-vn} and \eqref{eq:profile-rem-Lq} imply
\begin{align}\label{eq:selected-mass-sum}
        \sum_{j=1}^\infty m_j\geq\delta_0 .
\end{align}
Indeed, for every fixed $J$, \eqref{eq:weighted-BL-Lq} gives
\begin{align}
        \delta_0
        \leq
        \sum_{j=1}^J m_j
        +
        \limsup_{n\to\infty}\int_{\R^d}|r_n^J|^{p+1}\phi_n^2\,\ud x .
\end{align}
Since $0\leq\phi_n\leq1$ and \eqref{eq:profile-rem-Lq} holds, letting $J\to\infty$ yields \eqref{eq:selected-mass-sum}.
Let $\calJ:=\{j\in\N:m_j>0\}$. We now show that for each $j\in \calJ$,  $\psi^j$ is a stationary solution, i.e. it solves \eqref{eq:yamabe}. First observe that for each $j\in \N$ if we define
\begin{align}\label{eq:Gtau-profile-frame}
        G_{n}^{\ta}(y):=
        (\lambda_n^j)^{\frac{d-2}{2}}v_n(a_n^j+\lambda_n^j y),
\end{align}
then
\begin{align}\label{eq:Gtau-strong-profile}
        G_n^\ta\to \psi^j
        \qquad\textrm{strongly in }L^2_{\loc}(\R^d).
\end{align}
To see this, fix $j\leq J$. Then the profile $(\lam^j_n)^{\frac{d-2}{2}}\psi_n^j(a^j_n+\lam^j_n \cdot)\to \psi^j$ in $L^2_{\operatorname{loc}}$ as $n\to \infty$, while any other fixed profile for $k\neq j$, $(\lam^j_n)^{\frac{d-2}{2}}\psi_n^k(a^j_n+\lam^j_n \cdot)\to 0$ strongly in $L^2_{\loc}(\R^d)$ as $n\to \infty$ by asymptotic orthogonality of parameters. Moreover, by \eqref{eq:profile-remainder-frame-vanishing}, the rescaled remainder term converges weakly to zero in $\dot H^1$ and hence in
$L^{p+1}$. On every compact set $K\Subset\R^d$, it is bounded in $H^1(K)$, therefore, Rellich's theorem gives strong convergence to zero in $L^2(K)$.

Now fix $j\in\calJ$. We first rule out the case $h_n/(\lambda_n^j)^2\to0$ as $n\to \infty$. Suppose, toward a contradiction, that after passing to a subsequence $h_n/(\lambda_n^j)^2\to0$ as $n\to \infty$. Define
\begin{align}\label{eq:Gsigma-def}
        G_n^\sigma(y):=
        (\lambda_n^j)^{\frac{d-2}{2}}
        (u(\sigma_n,a_n^j+\lambda_n^j y)-W(a_n^j+\lambda_n^j y)),
        \quad
        \Phi_{j,n}(y):=\phi_n(a_n^j+\lambda_n^j y).
\end{align}
For every fixed $R<\infty$, Cauchy--Schwarz gives
\begin{align}\label{eq:rescaled-two-slices-close}
        \|G_n^\ta-G_n^\sigma\|_{L^2(B(0,R))}
        &\leq
        (\lambda_n^j)^{-1}
        \int_{\sigma_n}^{\ta_n}
        \|\partial_tu(t)\|_{L^2(B(a_n^j,R\lambda_n^j))}\,\ud t \notag\\
        &\leq
        \left(\frac{h_n}{(\lambda_n^j)^2}\right)^{1/2}
        \left(
        \int_{\sigma_n}^{\ta_n}\|\partial_tu(t)\|_{L^2}^2\,\ud t
        \right)^{1/2}
        \to0
        \qquad\textrm{as }n\to\infty .
\end{align}
Since $\supp\phi_n\subset\Omega_n^2$, scale invariance of the $L^{p+1}$-norm and
\eqref{eq:Lq-prop-initial-small} implies
\begin{align}\label{eq:Phi-Gsigma-small}
        \|\Phi_{j,n}G_n^\sigma\|_{L^{p+1}(\R^d)}
        \leq
        \|u(\sigma_n)-W\|_{L^{p+1}(\Omega_n^2)}
        \to0
        \qquad\textrm{as }n\to\infty .
\end{align}
Thus $\Phi_{j,n}G_n^\sigma\to0$ in $L^2_{\loc}$ as $n\to\infty$. Combining
\eqref{eq:Gtau-strong-profile}, \eqref{eq:rescaled-two-slices-close}, and
\eqref{eq:Phi-Gsigma-small} gives
\begin{align}\label{eq:Phi-Qj-L2-zero}
        \Phi_{j,n}\psi^j\to0
        \quad\textrm{in }L^2_{\loc}(\R^d)
        \quad\textrm{as }n\to\infty .
\end{align}
This implies that $\int_{\R^d}|\psi^j(y)|^{p+1}\Phi_{j,n}(y)^2\,\ud y\to0$ as $n\to \infty.$ To see this first let $\varepsilon>0$. Since $\psi^j\in L^{p+1}(\R^d)$, we may choose $R\geq 1$ large enough such that
\begin{align}\label{eq:Phi-Qj-tail-small}
        \int_{\R^d\setminus B(0,R)}|\psi^j(y)|^{p+1}\,\ud y<\varepsilon .
\end{align}
After fixing this $R$, we may choose $M\geq 1$ large enough such that
\begin{align}\label{eq:Phi-Qj-large-values-small}
        \int_{B(0,R)\cap\{|\psi^j|>M\}}|\psi^j(y)|^{p+1}\,\ud y<\varepsilon .
\end{align}
Using $0\leq\Phi_{j,n}\leq1$, we split
\begin{align}\label{eq:Phi-Qj-truncation-split}
        \int_{\R^d}|\psi^j|^{p+1}\Phi_{j,n}^2\,\ud y
        &\leq
        \int_{\R^d\setminus B(0,R)}|\psi^j|^{p+1}\,\ud y
        +
        \int_{B(0,R)\cap\{|\psi^j|>M\}}|\psi^j|^{p+1}\,\ud y \notag\\
        &\quad+
        \int_{B(0,R)\cap\{|\psi^j|\leq M\}}|\psi^j|^{p+1}\Phi_{j,n}^2\,\ud y \\
        &\leq 2\veps + M^{p-1}\int_{B(0,R)}|\psi^j|^2\Phi_{j,n}^2\,\ud y.
\end{align}
Taking $\limsup_{n\to\infty}$ in \eqref{eq:Phi-Qj-truncation-split} and using \eqref{eq:Phi-Qj-L2-zero}, we obtain
\begin{align}
        \limsup_{n\to\infty}
        \int_{\R^d}|\psi^j(y)|^{p+1}\Phi_{j,n}(y)^2\,\ud y
        \leq 2\varepsilon .
\end{align}
Letting $\varepsilon\to0$ gives
\begin{align}\label{eq:weighted-mass-zero-contradiction}
        \int_{\R^d}|\psi^j(y)|^{p+1}\Phi_{j,n}(y)^2\,\ud y\to0
        \qquad\textrm{as }n\to\infty .
\end{align}
This contradicts $m_j>0$. Therefore
\begin{align}\label{eq:h-over-lambda-liminf}
        \liminf_{n\to\infty}\frac{h_n}{(\lambda_n^j)^2}>0 .
\end{align}
Since $h_n\to0$ and $h_n/\ell_n^2\to0$ as $n\to\infty$, \eqref{eq:h-over-lambda-liminf} implies
\begin{align}\label{eq:selected-scale-small}
        \lambda_n^j\to0,
        \qquad
        \frac{\lambda_n^j}{\ell_n}\to0
        \qquad\textrm{as }n\to\infty .
\end{align}

We prove that $\psi^j$ is stationary. By \eqref{eq:h-over-lambda-liminf}, there exists
$\delta_j>0$ such that $\delta_j(\lambda_n^j)^2\leq h_n$ for all large $n\in \N$. For $s\in[-\delta_j,0]$, define
\begin{align}
        U_n(y,s)&:=
        (\lambda_n^j)^{\frac{d-2}{2}}
        u(\ta_n+(\lambda_n^j)^2s,a_n^j+\lambda_n^j y),
        \quad
        W_n(y):=
        (\lambda_n^j)^{\frac{d-2}{2}}
        W(a_n^j+\lambda_n^j y),
        \notag\\
        V_n(y,s)&:=U_n(y,s)-W_n.
        \label{eq:backward-rescaling}
\end{align}
Then $V_n$ solves
\begin{align}\label{eq:Vn-equation}
        \partial_sV_n=\Delta V_n+\left(|V_n+W_n|^{p-1}(V_n+W_n)-|W_n|^{p-1}W_n\right).
\end{align}
Moreover,
\begin{align}\label{eq:rescaled-tension-small}
        \int_{-\delta_j}^0\|\partial_sV_n(s)\|_{L^2_y}^2\,\ud s
        =
        \int_{\ta_n-\delta_j(\lambda_n^j)^2}^{\ta_n}
        \|\partial_tu(t)\|_{L^2_x}^2\,\ud t
        \to0
        \qquad\textrm{as }n\to\infty .
\end{align}
For every compact $K\Subset\R^d$,
\begin{align}\label{eq:Vn-time-constant}
        \sup_{s\in[-\delta_j,0]}
        \|V_n(s)-V_n(0)\|_{L^2(K)}
        \leq
        \delta_j^{1/2}
        \left(
        \int_{-\delta_j}^0\|\partial_sV_n(\rho)\|_{L^2}^2\,\ud\rho
        \right)^{1/2}
        \to0
        \qquad\textrm{as }n\to\infty .
\end{align}
Since $V_n(0)=G_n^\ta$, \eqref{eq:Gtau-strong-profile} gives
\begin{align}\label{eq:Vn-to-Qj-L2}
        V_n\to \psi^j
        \qquad\textrm{strongly in }L^\infty([-\delta_j,0];L^2_{\loc}(\R^d))
        \qquad\textrm{as }n\to\infty .
\end{align}
Interpolating with the uniform $L^\infty_sL^{p+1}_y$ bound gives
\begin{align}\label{eq:Vn-to-Qj-Lr}
        V_n\to \psi^j
        \qquad\textrm{strongly in }L^r_{\loc}((-\delta_j,0)\times\R^d).
        \quad\textrm{for every }2\leq r<p+1,
        \qquad\textrm{as }n\to\infty .
\end{align}
Also, since $\lambda_n^j\to0$ as $n\to\infty$ and $W\in L^{p+1}(\R^d)$,
\begin{align}\label{eq:Wn-to-zero}
        W_n\to0
        \qquad\textrm{strongly in }L^{p+1}_{\loc}(\R^d)
        \qquad\textrm{as }n\to\infty .
\end{align}
We next justify the convergence of the nonlinear term. Let $K\Subset\R^d$ and choose
$r$ such that $\max\{2,(p+1)/2\}<r<p+1$. By \eqref{eq:Vn-to-Qj-Lr} and
\eqref{eq:Wn-to-zero}, $V_n+W_n\to \psi^j$ strongly in
$L^r((-\delta_j,0)\times K)$, hence also strongly in
$L^{(p+1)/2}((-\delta_j,0)\times K)$. Since $V_n+W_n$ and $\psi^j$ are bounded in
$L^{p+1}((-\delta_j,0)\times K)$, the pointwise inequality
\begin{align}\label{eq:nonlinearity-L1-ineq}
        \left||a|^{p-1}a-|b|^{p-1}b\right|
        \lesssim
        \left(|a|^{p-1}+|b|^{p-1}\right)|a-b|
\end{align}
and Hölder's inequality imply
\begin{align}\label{eq:nonlinearity-main-L1}
        |V_n+W_n|^{p-1}(V_n+W_n)
        \to
        |\psi^j|^{p-1}\psi^j
        \qquad\textrm{in }L^1((-\delta_j,0)\times K)
        \qquad\textrm{as }n\to\infty .
\end{align}
Moreover, \eqref{eq:Wn-to-zero} gives
$|W_n|^{p-1}W_n\to0$ in $L^{\frac{p+1}{p}}_{\loc}$, and hence in $L^1_{\loc}$, as $n\to\infty$.
Therefore,
\begin{align}\label{eq:nonlinearity-limit}
        |V_n+W_n|^{p-1}(V_n+W_n)-|W_n|^{p-1}W_n
        \to
        |\psi^j|^{p-1}\psi^j
        \qquad\textrm{in }L^1_{\loc}
        \qquad\textrm{as }n\to\infty .
\end{align}
The estimate \eqref{eq:nonlinearity-L1-ineq} is valid for every $p>1$, so this covers both the cases $p\geq2$ and $1<p<2$.

We also record the weak convergence of the spatial gradients. Since $W\in\dot H^1(\R^d)$ and $u$ is uniformly bounded in $\dot H^1(\R^d)$, the sequence $V_n$ is bounded in
$L^\infty([-\delta_j,0];\dot H^1(\R^d))$. Hence, after passing to a subsequence,
\begin{align}\label{eq:Vn-weak-gradient-short-time-prop}
        \nabla V_n\rightharpoonup \nabla \psi^j
        \quad\textrm{weakly in }L^2_{\loc}((-\delta_j,0)\times\R^d).
\end{align}
Testing \eqref{eq:Vn-equation} against $\zeta(y)\rho(s)$, where
$\zeta\in C_c^\infty(\R^d)$ and $\rho\in C_c^\infty((-\delta_j,0))$, gives
\begin{align}\label{eq:Vn-weak-form}
        \int_{-\delta_j}^0\int_{\R^d}\partial_sV_n\,\zeta\rho\,\ud y\ud s
        &=
        -\int_{-\delta_j}^0\int_{\R^d}\nabla V_n\cdot\nabla\zeta\,\rho\,\ud y\ud s \\
        &\quad+
        \int_{-\delta_j}^0\int_{\R^d}
        \left(
        |V_n+W_n|^{p-1}(V_n+W_n)-|W_n|^{p-1}W_n
        \right)\zeta\rho\,\ud y\ud s .
\end{align}
The first term in \eqref{eq:Vn-weak-form} tends to $0$ as $n\to\infty$ by
\eqref{eq:rescaled-tension-small}. Passing to the limit in the remaining two terms by
\eqref{eq:Vn-weak-gradient-short-time-prop} and \eqref{eq:nonlinearity-limit}, we obtain
\begin{align}
        0
        =
        -\left(\int_{-\delta_j}^0\rho(s)\,\ud s\right)
        \int_{\R^d}\nabla \psi^j\cdot\nabla\zeta\,\ud y
        +
        \left(\int_{-\delta_j}^0\rho(s)\,\ud s\right)
        \int_{\R^d}|\psi^j|^{p-1}\psi^j\zeta\,\ud y .
\end{align}
Choosing $\rho$ such that $\int_{-\delta_j}^0 \rho(s)\,\ud s\neq 0$ gives
\begin{align}\label{eq:Qj-stationary}
        \int_{\R^d}\nabla \psi^j\cdot\nabla\zeta\,\ud y
        =
        \int_{\R^d}|\psi^j|^{p-1}\psi^j\zeta\,\ud y .
\end{align}
Thus $\psi^j$ is a stationary solution of \eqref{eq:yamabe}. Since $m_j>0$, each $\psi^j$ with $j\in\calJ$ is nontrivial. Hence $ \|\nabla \psi^j\|_{L^2}^2\geq\bar E_*.$ The Pythagorean expansion in the profile decomposition gives
\begin{align}\label{eq:profile-energy-sum-bound}
\sum_{j=1}^{\infty}\|\nabla \psi^j\|_{L^2}^2 \leq \limsup_{n\to\infty}\|v_n\|_{\dot H^1}^2<\infty .
\end{align}
Therefore $\calJ$ is finite. Relabel it as $\calJ=\{1,\ldots,J_*\}$. Set
$S_n:=\sum_{j=1}^{J_*}\psi_n^j$ and $z_n:=v_n-S_n$. We next show that the localized $L^{p+1}$-norm of $z_n$ tends to $0$ as $n\to \infty$. To see this, first note that for every fixed $J\geq J_*$, we have
\begin{align}
        z_n
        =
        \sum_{j=J_*+1}^J\psi_n^j+r_n^J .
\end{align}
Using \eqref{eq:weighted-BL-Lq} for this finite decomposition gives
\begin{align}
        \limsup_{n\to\infty}
        \int_{\R^d}|z_n|^{p+1}\phi_n^2\,\ud x
        \leq
        \sum_{j=J_*+1}^J m_j
        +
        \limsup_{n\to\infty}
        \int_{\R^d}|r_n^J|^{p+1}\phi_n^2\,\ud x .
\end{align}
Since $m_j=0$ for $j>J_*$ and $0\leq\phi_n\leq1$, letting $J\to\infty$ and using \eqref{eq:profile-rem-Lq} gives
\begin{align}\label{eq:zn-weighted-small-final}
        \int_{\R^d}|z_n|^{p+1}\phi_n^2\,\ud x\to0
        \qquad\textrm{as }n\to\infty,
\end{align}
and
\begin{align}\label{eq:selected-mass-positive}
        \sum_{j=1}^{J_*}m_j\geq\delta_0 .
\end{align}

Set $A_n:=W+S_n$. We include $W$ as the fixed profile corresponding to $\lambda_n^0=1$ and $a_n^0=0$. By \eqref{eq:selected-scale-small}, this fixed profile is asymptotically orthogonal to each selected profile $\psi_n^j$, $1\leq j\leq J_*$. Applying Lemma~\ref{lem:weighted-profile-decoupling} to the finite stationary family $W,\psi_n^1,\ldots,\psi_n^{J_*}$, with no remainder term, gives
\begin{align}\label{eq:weighted-energy-decoupling-An}
        E_{\phi_n}(A_n)
        \geq
        E_{\phi_n}(W)+\sum_{j=1}^{J_*}E_{\phi_n}(\psi_n^j)-o_n(1).
\end{align}
We next add $z_n$. Since $0\leq\phi_n\leq1$ and $p+1\geq2$, \eqref{eq:zn-weighted-small-final} implies
\begin{align}\label{eq:zn-weighted-small-use}
        \|z_n\phi_n\|_{L^{p+1}(\R^d)}
        +
        \|z_n\phi_n^2\|_{L^{p+1}(\R^d)}
        +
        \int_{\R^d}|z_n|^{p+1}\phi_n^2\,\ud x
        \to0
        \qquad\textrm{as }n\to\infty .
\end{align}
For $Q_n=W$ or $Q_n=\psi_n^j$, testing \eqref{eq:yamabe} with the function $z_n\phi_n^2$ gives
\begin{align}\label{eq:Qn-zn-gradient-identity}
        \int_{\R^d}\nabla Q_n\cdot\nabla z_n\,\phi_n^2\,\ud x
        =
        \int_{\R^d}|Q_n|^{p-1}Q_nz_n\phi_n^2\,\ud x
        -2\int_{\R^d}z_n\phi_n\nabla Q_n\cdot\nabla\phi_n\,\ud x .
\end{align}
Therefore, by Hölder's inequality, \eqref{eq:Lq-prop-cutoff}, and \eqref{eq:zn-weighted-small-use},
\begin{align}\label{eq:Qn-zn-gradient-small}
        \left|
        \int_{\R^d}\nabla Q_n\cdot\nabla z_n\,\phi_n^2\,\ud x
        \right|
        &\lesssim
        \|Q_n\|_{L^{p+1}}^p\|z_n\phi_n^2\|_{L^{p+1}}
        +
        \|\nabla Q_n\|_{L^2}\|\nabla\phi_n\|_{L^d}\|z_n\phi_n\|_{L^{p+1}}
        =
        o_n(1).\qquad
\end{align}
Summing \eqref{eq:Qn-zn-gradient-small} over $Q_n=W,\psi_n^1,\ldots,\psi_n^{J_*}$ gives
\begin{align}\label{eq:An-zn-gradient-small}
        \int_{\R^d}\nabla A_n\cdot\nabla z_n\,\phi_n^2\,\ud x=o_n(1).
\end{align}
Also, since $\{A_n\}_{n\in \N}$ is bounded in $L^{p+1}(\R^d)$, the pointwise inequality
\begin{align}\label{eq:potential-difference-simple}
        \left||a+b|^{p+1}-|a|^{p+1}\right|
        \lesssim
        |a|^p|b|+|b|^{p+1}
\end{align}
and Hölder's inequality imply
\begin{align}\label{eq:An-zn-potential-small}
        \int_{\R^d}
        \left||A_n+z_n|^{p+1}-|A_n|^{p+1}\right|\phi_n^2\,\ud x
        &\lesssim
        \|A_n\|_{L^{p+1}}^p\|z_n\phi_n^2\|_{L^{p+1}}
        +
        \int_{\R^d}|z_n|^{p+1}\phi_n^2\,\ud x
        =
        o_n(1).\qquad
\end{align}
Expanding the nonlinear energy and using \eqref{eq:An-zn-gradient-small} and \eqref{eq:An-zn-potential-small}, we obtain
\begin{align}\label{eq:add-zn-energy-lower}
        E_{\phi_n}(A_n+z_n)-E_{\phi_n}(A_n)
        &=
        \frac12\int_{\R^d}|\nabla z_n|^2\phi_n^2\,\ud x
        +
        \int_{\R^d}\nabla A_n\cdot\nabla z_n\,\phi_n^2\,\ud x \notag\\
        &\quad
        -\frac1{p+1}\int_{\R^d}
        \left(|A_n+z_n|^{p+1}-|A_n|^{p+1}\right)\phi_n^2\,\ud x
        \geq -o_n(1).
\end{align}
Combining \eqref{eq:weighted-energy-decoupling-An} and \eqref{eq:add-zn-energy-lower}, we obtain
\begin{align}\label{eq:weighted-energy-decoupling-An-zn}
        E_{\phi_n}(A_n+z_n)-E_{\phi_n}(W)
        \geq
        \sum_{j=1}^{J_*}E_{\phi_n}(\psi_n^j)-o_n(1).
\end{align}
Now, testing \eqref{eq:yamabe} for $\psi_n^j$ with the test function $\psi_n^j\phi_n^2$ yields
\begin{align}\label{eq:selected-profile-local-energy-identity}
        E_{\phi_n}(\psi_n^j)
        =
        \frac1d\int_{\R^d}|\psi_n^j|^{p+1}\phi_n^2\,\ud x
        -
        \int_{\R^d}\psi_n^j\phi_n\nabla\psi_n^j\cdot\nabla\phi_n\,\ud x .
\end{align}
Set $\Phi_{j,n}(y):=\phi_n(a_n^j+\lambda_n^j y)$. By \eqref{eq:Lq-prop-cutoff} and \eqref{eq:selected-scale-small}, $\|\nabla\Phi_{j,n}\|_{L^\infty}\leq\lambda_n^j\|\nabla\phi_n\|_{L^\infty}\to0$ as $n\to\infty$ and $\|\nabla\Phi_{j,n}\|_{L^d}\lesssim1$. Changing variables gives
\begin{align}\label{eq:selected-profile-cutoff-gradient-change}
        \int_{\R^d}\psi_n^j\phi_n\nabla\psi_n^j\cdot\nabla\phi_n\,\ud x
        =
        \int_{\R^d}\psi^j\Phi_{j,n}\nabla\psi^j\cdot\nabla\Phi_{j,n}\,\ud y .
\end{align}
For each fixed $R<\infty$,
\begin{align}\label{eq:selected-profile-cutoff-gradient-inside}
        \left|
        \int_{B(0,R)}\psi^j\Phi_{j,n}\nabla\psi^j\cdot\nabla\Phi_{j,n}\,\ud y
        \right|
        \leq
        \|\nabla\Phi_{j,n}\|_{L^\infty}
        \|\psi^j\|_{L^2(B(0,R))}
        \|\nabla\psi^j\|_{L^2(B(0,R))}
        \to0,
\end{align}
as $n\to \infty.$ On the complement region, Hölder's inequality gives
\begin{align}\label{eq:selected-profile-cutoff-gradient-outside}
        \left|
        \int_{\R^d\setminus B(0,R)}
        \psi^j\Phi_{j,n}\nabla\psi^j\cdot\nabla\Phi_{j,n}\,\ud y
        \right|
        \leq
        \|\psi^j\|_{L^{p+1}(\R^d\setminus B(0,R))}
        \|\nabla\psi^j\|_{L^2}
        \|\nabla\Phi_{j,n}\|_{L^d}.
\end{align}
Letting first $n\to\infty$ in \eqref{eq:selected-profile-cutoff-gradient-inside} and then $R\to\infty$ in \eqref{eq:selected-profile-cutoff-gradient-outside}, we obtain
\begin{align}\label{eq:selected-profile-cutoff-gradient-small}
        \int_{\R^d}\psi_n^j\phi_n\nabla\psi_n^j\cdot\nabla\phi_n\,\ud x=o_n(1).
\end{align}
Hence \eqref{eq:selected-profile-local-energy-identity} gives
\begin{align}\label{eq:selected-profile-energy-mass}
        E_{\phi_n}(\psi_n^j)
        =
        \frac1d\int_{\R^d}|\psi_n^j|^{p+1}\phi_n^2\,\ud x+o_n(1)
        \qquad\textrm{for }1\leq j\leq J_* .
\end{align}
Since $A_n+z_n=W+v_n=u(\ta_n)$, \eqref{eq:weighted-energy-decoupling-An-zn} and \eqref{eq:selected-profile-energy-mass} imply
\begin{align}\label{eq:positive-energy-defect-proof}
        E_{\phi_n}(u(\ta_n))-E_{\phi_n}(W)
        =
        E_{\phi_n}(W+v_n)-E_{\phi_n}(W) \geq
        \frac1d\sum_{j=1}^{J_*}
        \int_{\R^d}|\psi_n^j|^{p+1}\phi_n^2\,\ud x
        -o_n(1).
\end{align}
Then using \eqref{eq:selected-mass-positive}, we obtain
\begin{align}\label{eq:positive-energy-defect}
        \liminf_{n\to\infty}
        \left(E_{\phi_n}(u(\ta_n))-E_{\phi_n}(W)\right)
        \geq
        \frac{\delta_0}{d}>0 .
\end{align}
On the other hand, applying Lemma~\ref{lem:short-time-nonlinear-energy-prop} with $\eta_n=\phi_n$ and using \eqref{eq:Lq-prop-energy-input}, we obtain
\begin{align}
        \lim_{n\to\infty}
        \left|E_{\phi_n}(u(\ta_n))-E_{\phi_n}(W)\right|=0,
\end{align}
which contradicts \eqref{eq:positive-energy-defect}.
\end{proof}
Combining the preceding two short-time propagation lemmas, we now propagate localized energy closeness on balls and punctured balls.
\begin{lem}[Short-time propagation of the localized energy]
\label{lem:short-time-prop}
Let $u(t)$ be a solution to \eqref{eqn:NLH} with initial data
$u(0)=u_0\in\dot H^1$. Let $T_+=T_+(u_0)$ denote its maximal time of existence and assume that
\begin{align}\label{eq:short-time-prop-energy-bound}
        \sup_{t\in[0,T_+)}\|u(t)\|_{\dot H^1}<\infty .
\end{align}
Let $0<\sigma_n<\ta_n<T_+$ satisfy $\sigma_n,\ta_n\to T_+$ and
$h_n:=\ta_n-\sigma_n\to0$ as $n\to\infty$. If $T_+=\infty$, this means
$\sigma_n,\ta_n\to\infty$. Let $W\in\dot H^1(\R^d)$ be a stationary solution of
\eqref{eq:yamabe}, possibly zero.

Let $r_n>0$ satisfy
\begin{align}\label{eq:short-time-prop-rn}
        \frac{\ta_n-\sigma_n}{r_n^2}\to0
        \qquad\textrm{as }n\to\infty .
\end{align}
Assume
\begin{align}\label{eq:short-time-prop-ball-input}
        \lim_{n\to \infty} \left[\bar E(u(\sigma_n)-W;B(0,2r_n))
        +
        \|u(\sigma_n)-W\|_{L^{p+1}(B(0,2r_n))}^{p+1}\right] = 0.
\end{align}
Then
\begin{align}\label{eq:en-discR}
        \lim_{n\to \infty}
        \left[\bar E(u(\ta_n)-W;B(0,r_n))
        +
        \|u(\ta_n)-W\|_{L^{p+1}(B(0,r_n))}^{p+1}\right]
        =0.
\end{align}

Next, let $\veps_n>0$ satisfy $\veps_n<r_n$ and
\begin{align}\label{eq:short-time-prop-epsn}
        \frac{\ta_n-\sigma_n}{\veps_n^2}\to0
        \qquad\textrm{as }n\to\infty .
\end{align}
Let $L\in\N$, $L\geq1$, and let $x_{\ell,n}\in\R^d$, $1\leq\ell\leq L$, satisfy
\begin{align}\label{eq:short-time-prop-puncture-sep}
        B(x_{\ell,n},\veps_n)\subset B(0,r_n),
        \qquad
        |x_{\ell,n}-x_{m,n}|\geq5\veps_n
        \quad\textrm{for }\ell\neq m .
\end{align}
Let $I_n=B(0,2r_n)\setminus \bigcup_{\ell=1}^L{B(x_{\ell,n},\veps_n/2)}$ and assume that
\begin{align}\label{eq:short-time-prop-annulus-input}
        \lim_{n\to\infty}
        \Bigl[
        &\bar E(u(\sigma_n)-W; I_n)
        +
        \|u(\sigma_n)-W\|_{L^{p+1}
        (I_n)}^{p+1}
        \Bigr]
        =0.
\end{align}
Then on $J_n:=B(0,r_n)\setminus\bigcup_{\ell=1}^L {B(x_{\ell,n},\veps_n)}$, we obtain
\begin{align}\label{eq:en-annulus}
        \lim_{n\to \infty}
        \left[\bar E(u(\ta_n)-W;J_n)
        +
        \|u(\ta_n)-W\|_{L^{p+1}(J_n)}^{p+1}\right]
        =0.
\end{align}
\end{lem}
\begin{proof}
We prove the ball and punctured ball case simultaneously. In the ball case, set
$\Omega_n^0:=B(0,r_n)$, $\Omega_n^1:=B(0,3r_n/2)$, $\Omega_n^2:=B(0,2r_n)$, and $\ell_n:=r_n$. In the punctured ball case, set $\ell_n:=\veps_n$,
\begin{align}\label{eq:short-time-punctured-regions}
        \Omega_n^0
        &:=
        B(0,r_n)\setminus
        \bigcup_{\ell=1}^L \overline{B(x_{\ell,n},\veps_n)},\quad \Omega_n^1
        :=
        B(0,3r_n/2)\setminus
        \bigcup_{\ell=1}^L \overline{B(x_{\ell,n},3\veps_n/4)}, \text{ and }\notag\\
        \Omega_n^2
        &:=
        B(0,2r_n)\setminus
        \bigcup_{\ell=1}^L \overline{B(x_{\ell,n},\veps_n/2)}.
\end{align}
Then $\Omega_n^0\subset\Omega_n^1\subset\Omega_n^2$ and
$(\ta_n-\sigma_n)/\ell_n^2\to0$ as $n\to\infty$. Using \eqref{eq:short-time-prop-puncture-sep} in the punctured ball case, choose
$\phi_n,\psi_n\in C_c^\infty(\Omega_n^2)$ such that
\begin{align}\label{eq:short-time-main-cutoffs}
        0\leq\phi_n,\psi_n\leq1,
        \qquad
        \phi_n\equiv1\textrm{ on }\Omega_n^1,
        \qquad
        \psi_n\equiv1\textrm{ on }\Omega_n^0,
        \qquad
        \supp\psi_n\subset\Omega_n^1,
\end{align}
and
\begin{align}\label{eq:short-time-main-cutoff-bounds}
        \|\nabla\phi_n\|_{L^\infty}
        +
        \|\nabla\psi_n\|_{L^\infty}
        \lesssim\ell_n^{-1},
        \qquad
        \|\nabla\phi_n\|_{L^d}
        +
        \|\nabla\psi_n\|_{L^d}
        \lesssim1 .
\end{align}
Now let $\eta_n=\phi_n$ or $\eta_n=\psi_n$, and set $e_n:=u(\sigma_n)-W$. Since $\supp\eta_n\subset\Omega_n^2$ and $0\leq\eta_n\leq1$,
\begin{align}\label{eq:short-time-main-initial-E}
        \left|E_{\eta_n}(u(\sigma_n))-E_{\eta_n}(W)\right|
        &\leq
        \frac12\int_{\Omega_n^2}|\nabla e_n|^2\,\ud x
        +
        \int_{\Omega_n^2}|\nabla W||\nabla e_n|\,\ud x \notag\\
        &\quad+
        \frac1{p+1}\int_{\Omega_n^2}
        \left||W+e_n|^{p+1}-|W|^{p+1}\right|\,\ud x \notag\\
        &\lesssim
        \|\nabla e_n\|_{L^2(\Omega_n^2)}^2
        +
        \|\nabla W\|_{L^2}\|\nabla e_n\|_{L^2(\Omega_n^2)} \notag\\
        &\quad+
        \|W\|_{L^{p+1}}^p\|e_n\|_{L^{p+1}(\Omega_n^2)}
        +
        \|e_n\|_{L^{p+1}(\Omega_n^2)}^{p+1}.
\end{align}
Each term on the right-hand side of the above display tends to $0$ as $n\to \infty$. Thus, in either case, Lemma~\ref{lem:short-time-nonlinear-energy-prop} applies with $\eta_n=\phi_n$ and with $\eta_n=\psi_n$, and gives
\begin{align}\label{eq:short-time-main-final-E}
        E_{\eta_n}(u(\ta_n))-E_{\eta_n}(W)\to0
        \quad\textrm{as }n\to\infty.
\end{align}
Applying Lemma~\ref{lem:short-time-Lq-prop} with $\Omega_n^1$, $\Omega_n^2$, $\ell_n$, and $\phi_n$ implies
\begin{align}\label{eq:short-time-main-Lq}
        \|u(\ta_n)-W\|_{L^{p+1}(\Omega_n^1)}\to0
        \quad\textrm{as }n\to\infty .
\end{align}
Set $w_n:=u(\ta_n)-W$. By \eqref{eq:short-time-main-final-E} with $\eta_n=\psi_n$,
\begin{align}\label{eq:short-time-main-Epsi}
        E_{\psi_n}(W+w_n)-E_{\psi_n}(W)\to0
        \quad\textrm{as }n\to\infty .
\end{align}
Expanding,
\begin{align}\label{eq:short-time-main-expand}
        E_{\psi_n}(W+w_n)-E_{\psi_n}(W)
        &=
        \frac12\int_{\R^d}|\nabla w_n|^2\psi_n^2\,\ud x
        +
        \int_{\R^d}\nabla W\cdot\nabla w_n\,\psi_n^2\,\ud x \notag\\
        &\quad
        -\frac1{p+1}\int_{\R^d}
        \left(|W+w_n|^{p+1}-|W|^{p+1}\right)\psi_n^2\,\ud x .
\end{align}
Using \eqref{eq:yamabe} for $W$ implies
\begin{align}\label{eq:short-time-main-cross-W}
        \int_{\R^d}\nabla W\cdot\nabla w_n\,\psi_n^2\,\ud x
        =
        \int_{\R^d}|W|^{p-1}Ww_n\psi_n^2\,\ud x
        -
        2\int_{\R^d}w_n\psi_n\nabla W\cdot\nabla\psi_n\,\ud x .
\end{align}
Because $\supp\psi_n\subset\Omega_n^1$ and \eqref{eq:short-time-main-Lq} holds,
\begin{align}\label{eq:short-time-main-Taylor-zero}
        \int_{\R^d}
        |
        |W+w_n|^{p+1}-|W|^{p+1}-(p+1)|W|^{p-1}Ww_n
        |\psi_n^2\,\ud x
        \to0
        \qquad\textrm{as }n\to\infty .
\end{align}
Also,
\begin{align}\label{eq:short-time-main-cutoff-W-zero}
        \left|
        \int_{\R^d}w_n\psi_n\nabla W\cdot\nabla\psi_n\,\ud x
        \right|
        \leq
        \|w_n\|_{L^{p+1}(\Omega_n^1)}
        \|\nabla W\|_{L^2}
        \|\nabla\psi_n\|_{L^d}
        \to0
        \qquad\textrm{as }n\to\infty .
\end{align}
Combining \eqref{eq:short-time-main-Epsi},
\eqref{eq:short-time-main-expand}, \eqref{eq:short-time-main-cross-W},
\eqref{eq:short-time-main-Taylor-zero}, and
\eqref{eq:short-time-main-cutoff-W-zero}, we obtain
\begin{align}\label{eq:short-time-main-grad-small}
        \int_{\R^d}|\nabla w_n|^2\psi_n^2\,\ud x\to0
        \qquad\textrm{as }n\to\infty .
\end{align}
Since $\psi_n\equiv1$ on $\Omega_n^0$,
\begin{align}\label{eq:short-time-main-Omega0}
        \bar E(u(\ta_n)-W;\Omega_n^0)\to0
        \qquad\textrm{as }n\to\infty.
\end{align}
Since $\Omega_n^0\subset\Omega_n^1$, \eqref{eq:short-time-main-Lq} also gives
\begin{align}\label{eq:short-time-main-Lq-Omega0}
        \|u(\ta_n)-W\|_{L^{p+1}(\Omega_n^0)}\to0
        \qquad\textrm{as }n\to\infty.
\end{align}
Combining \eqref{eq:short-time-main-Lq-Omega0} with
\eqref{eq:short-time-main-Omega0} gives the desired conclusions.
\end{proof}
By the same localized nonlinear energy estimate, the nonlinear energy
can also be propagated backwards in time. More precisely, if
$\eta_n$ satisfies the hypotheses of Lemma~\ref{lem:short-time-nonlinear-energy-prop}
and
\begin{align}\label{eq:short-time-nonlinear-energy-final-backward-input}
        \lim_{n\to\infty}
        |E_{\eta_n}(u(\ta_n))-E_{\eta_n}(W)|=0,
\end{align}
then
\begin{align}\label{eq:short-time-nonlinear-energy-initial-backward-output}
        \lim_{n\to\infty}
        |E_{\eta_n}(u(\sigma_n))-E_{\eta_n}(W)|=0.
\end{align}
Indeed, this follows immediately from
\eqref{eq:short-time-nonlinear-energy-conclusion} and the triangle
inequality.

We next establish the corresponding backwards-in-time analogue of
Lemma~\ref{lem:short-time-Lq-prop}.
\begin{lem}[Short-time propagation of the $L^{p+1}$ norm backwards in time]
\label{lem:short-time-Lp-prop-backward}
Let $u(t)$ be a solution to \eqref{eqn:NLH} with initial data
$u(0)=u_0\in\dot H^1$. Let $T_+=T_+(u_0)$ denote its maximal time of
existence and assume
\begin{align}\label{eq:backward-Lp-prop-energy-bound}
        \sup_{t\in[0,T_+)}\|u(t)\|_{\dot H^1}<\infty .
\end{align}
Let $0<\sigma_n<\ta_n<T_+$ satisfy $\sigma_n,\ta_n\to T_+$ and
$h_n:=\ta_n-\sigma_n\to0$ as $n\to\infty$. If $T_+=\infty$, this means
$\sigma_n,\ta_n\to\infty$. Let $W\in\dot H^1(\R^d)$ be a stationary solution
of \eqref{eq:yamabe}, possibly zero. Let $\Omega_n^1\subset\Omega_n^2$ be
open sets. Let $\phi_n\in C_c^\infty(\Omega_n^2)$ satisfy
\begin{align}\label{eq:backward-Lp-prop-cutoff}
        0\leq\phi_n\leq1,\quad
        \phi_n\equiv1 \textrm{ on } \Omega_n^1,\quad
        \|\nabla\phi_n\|_{L^\infty}\lesssim \ell_n^{-1},\quad
        \|\nabla\phi_n\|_{L^d}\lesssim1,
\end{align}
for some $\ell_n>0$ with $h_n/\ell_n^2\to0$ as $n\to\infty$. Assume
\begin{align}\label{eq:backward-Lp-prop-final-small}
        \lim_{n\to \infty}\|u(\ta_n)-W\|_{L^{p+1}(\Omega_n^2)}=0
\end{align}
and
\begin{align}\label{eq:backward-Lp-prop-energy-input}
        \lim_{n\to \infty}|E_{\phi_n}(u(\ta_n))-E_{\phi_n}(W)|=0.
\end{align}
Then
\begin{align}\label{eq:backward-Lp-prop-conclusion}
        \lim_{n\to \infty}\|u(\sigma_n)-W\|_{L^{p+1}(\Omega_n^1)}=0.
\end{align}
\end{lem}

\begin{proof}
Set $v_n:=u(\sigma_n)-W$. Suppose that
\eqref{eq:backward-Lp-prop-conclusion} fails. Passing to a subsequence,
there exists $\delta_0>0$ such that
\begin{align}\label{eq:backward-Lp-weighted-lower-vn}
        \int_{\R^d}|v_n|^{p+1}\phi_n^2\,\ud x\geq\delta_0 .
\end{align}
Applying Proposition~\ref{prop:profile-decomp} to $\{v_n\}_{n\in\N}$, we obtain
for every fixed $J$,
\begin{align}\label{eq:backward-profile-decomp-vn}
        v_n=\sum_{j=1}^J\psi_n^j+r_n^J,
        \qquad
        \psi_n^j(x):=(\lambda_n^j)^{-\frac{d-2}{2}}
        \psi^j\left(\frac{x-a_n^j}{\lambda_n^j}\right),
\end{align}
with parameters satisfying part $\mathrm{(a)}$ of
Proposition~\ref{prop:profile-decomp} and
\begin{align}\label{eq:backward-profile-rem-Lp}
        \lim_{J\to\infty}\limsup_{n\to\infty}
        \|r_n^J\|_{L^{p+1}}=0 .
\end{align}
After a diagonal extraction, the limits
\begin{align}\label{eq:backward-profile-weighted-mass}
        m_j:=
        \lim_{n\to\infty}
        \int_{\R^d}|\psi_n^j|^{p+1}\phi_n^2\,\ud x
\end{align}
exist. By Lemma~\ref{lem:weighted-profile-decoupling},
\eqref{eq:backward-Lp-weighted-lower-vn}, and
\eqref{eq:backward-profile-rem-Lp}, we obtain
\begin{align}\label{eq:backward-selected-mass-sum}
        \sum_{j=1}^\infty m_j\geq\delta_0 .
\end{align}
Let $\calJ:=\{j\in\N:m_j>0\}$. Fix $j\in\calJ$. We first rule out the case $h_n/(\lambda_n^j)^2\to0$ as $n\to\infty$. Define
\begin{align}\label{eq:Gsigma-profile-frame-backward}
        G_n^\sigma(y):=
        (\lambda_n^j)^{\frac{d-2}{2}}
        v_n(a_n^j+\lambda_n^j y).
\end{align}
Then, as in the proof of Lemma~\ref{lem:short-time-Lq-prop},
\begin{align}\label{eq:Gsigma-strong-profile-backward}
        G_n^\sigma\to\psi^j \quad\textrm{strongly in }L^2_{\loc}(\R^d).
\end{align}
Suppose, toward a contradiction, that after passing to a subsequence
$h_n/(\lambda_n^j)^2\to0$ as $n\to \infty$. Define
\begin{align}\label{eq:Gtau-def-backward}
        G_n^\ta(y)
        :=
        (\lambda_n^j)^{\frac{d-2}{2}}
        (u(\ta_n,a_n^j+\lambda_n^j y)-W(a_n^j+\lambda_n^j y)), \quad  \Phi_{j,n}(y)&:=\phi_n(a_n^j+\lambda_n^j y).
\end{align}
For every fixed $0<R<\infty$, Cauchy--Schwarz gives, 
\begin{align}\label{eq:backward-two-slices-close}
        \|G_n^\ta-G_n^\sigma\|_{L^2(B(0,R))}
        &\leq
        (\lambda_n^j)^{-1}
        \int_{\sigma_n}^{\ta_n}
        \|\partial_tu(t)\|_{L^2(B(a_n^j,R\lambda_n^j))}\,\ud t \notag\\
        &\leq
        \left(\frac{h_n}{(\lambda_n^j)^2}\right)^{1/2}
        \left(
        \int_{\sigma_n}^{\ta_n}\|\partial_tu(t)\|_{L^2}^2\,\ud t
        \right)^{1/2}
        \to0,
\end{align}
as $n\to \infty.$ Since $\supp\phi_n\subset\Omega_n^2$, scale invariance of the
$L^{p+1}$-norm and \eqref{eq:backward-Lp-prop-final-small} imply
\begin{align}\label{eq:Phi-Gtau-small-backward}
        \|\Phi_{j,n}G_n^\ta\|_{L^{p+1}(\R^d)}
        \leq
        \|u(\ta_n)-W\|_{L^{p+1}(\Omega_n^2)}
        \to0,
\end{align}
as $n\to \infty.$ Combining \eqref{eq:Gsigma-strong-profile-backward},
\eqref{eq:backward-two-slices-close}, and
\eqref{eq:Phi-Gtau-small-backward}, we obtain
$\Phi_{j,n}\psi^j\to0$ in $L^2_{\loc}(\R^d)$. By the same truncation argument
used in the proof of Lemma~\ref{lem:short-time-Lq-prop}, this gives
\begin{align}\label{eq:backward-weighted-mass-zero-contradiction}
        \int_{\R^d}|\psi^j(y)|^{p+1}\Phi_{j,n}(y)^2\,\ud y\to0,
\end{align}
as $n\to \infty$ which contradicts $m_j>0$. Therefore
\begin{align}\label{eq:backward-h-over-lambda-liminf}
        \liminf_{n\to\infty}\frac{h_n}{(\lambda_n^j)^2}>0 .
\end{align}
Since $h_n\to0$ and $h_n/\ell_n^2\to0$ as $n\to \infty$, it follows that, as $n\to \infty$,
\begin{align}\label{eq:backward-selected-scale-small}
        \lambda_n^j\to0,
        \qquad
        \frac{\lambda_n^j}{\ell_n}\to0 .
\end{align}

We next prove that $\psi^j$ is stationary. By
\eqref{eq:backward-h-over-lambda-liminf}, there exists $\delta_j>0$ such that
$\delta_j(\lambda_n^j)^2\leq h_n$ for all large $n$. For
$s\in[0,\delta_j]$, define
\begin{align}\label{eq:forward-rescaling-backward-lemma}
        U_n(y,s)&:=
        (\lambda_n^j)^{\frac{d-2}{2}}
        u(\sigma_n+(\lambda_n^j)^2s,a_n^j+\lambda_n^j y),\quad 
        W_n(y):=
        (\lambda_n^j)^{\frac{d-2}{2}}
        W(a_n^j+\lambda_n^j y),
        \notag\\
        V_n(y,s)&:=U_n(y,s)-W_n .
\end{align}
Then $V_n$ solves
\begin{align}\label{eq:Vn-equation-backward-time}
        \partial_sV_n
        =
        \Delta V_n+
        \left(
        |V_n+W_n|^{p-1}(V_n+W_n)-|W_n|^{p-1}W_n
        \right).
\end{align}
Moreover, as $n\to \infty$,
\begin{align}\label{eq:rescaled-tension-small-backward}
        \int_0^{\delta_j}\|\partial_sV_n(s)\|_{L^2_y}^2\,\ud s
        =
        \int_{\sigma_n}^{\sigma_n+\delta_j(\lambda_n^j)^2}
        \|\partial_tu(t)\|_{L^2_x}^2\,\ud t
        \to0 .
\end{align}
Thus, for every compact $K\Subset\R^d$, as $n\to \infty$,
\begin{align}\label{eq:Vn-time-constant-backward}
        \sup_{s\in[0,\delta_j]}\|V_n(s)-V_n(0)\|_{L^2(K)}
        \leq
        \delta_j^{1/2}
        \left(
        \int_0^{\delta_j}\|\partial_sV_n(\rho)\|_{L^2}^2\,\ud\rho
        \right)^{1/2}
        \to0 .
\end{align}
Since $V_n(0)=G_n^\sigma$, \eqref{eq:Gsigma-strong-profile-backward}
implies that as $n \to \infty$,
\begin{align}\label{eq:Vn-to-Qj-L2-backward}
        V_n\to\psi^j
        \qquad\textrm{strongly in }L^\infty([0,\delta_j];L^2_{\loc}(\R^d)).
\end{align}
Interpolating with the uniform $L^\infty_sL^{p+1}_y$ bound gives strong
convergence in $L^r_{\loc}((0,\delta_j)\times\R^d)$ for every $2\leq r<p+1$. Since
$\lambda_n^j\to0$ and $W\in L^{p+1}(\R^d)$, we also have
$W_n\to0$ strongly in $L^{p+1}_{\loc}(\R^d)$. Thus the nonlinear term in
\eqref{eq:Vn-equation-backward-time} converges to
$|\psi^j|^{p-1}\psi^j$ in $L^1_{\loc}$. Moreover, after passing to a subsequence,
\begin{align}
\nabla V_n\rightharpoonup\nabla\psi^j
\quad\textrm{weakly in }
L^2_{\loc}((0,\delta_j)\times\R^d),
\end{align}
as $n\to \infty.$

Testing \eqref{eq:Vn-equation-backward-time} against
$\zeta(y)\rho(s)$, where $\zeta\in C_c^\infty(\R^d)$ and
$\rho\in C_c^\infty((0,\delta_j))$ with
$\int_0^{\delta_j}\rho(s)\,\ud s=1$, and then passing to the limit using
\eqref{eq:rescaled-tension-small-backward}, the weak convergence of
$\nabla V_n$, and the convergence of the nonlinear term, gives
\begin{align}\label{eq:Qj-stationary-backward}
        \int_{\R^d}\nabla\psi^j\cdot\nabla\zeta\,\ud y
        =
        \int_{\R^d}|\psi^j|^{p-1}\psi^j\zeta\,\ud y .
\end{align}
Thus $\psi^j$ is a weak stationary solution of \eqref{eq:yamabe}.

Since $m_j>0$, each $\psi^j$ with $j\in\calJ$ is nontrivial. Hence the
energy lower bound for nontrivial stationary solutions and the
Pythagorean expansion imply that $\calJ$ is finite. Relabel
$\calJ=\{1,\ldots,J_*\}$. Set $S_n:=\sum_{j=1}^{J_*}\psi_n^j$ and
$z_n:=v_n-S_n$. As in the proof of Lemma~\ref{lem:short-time-Lq-prop},
Lemma~\ref{lem:weighted-profile-decoupling} implies that as $n\to \infty$
\begin{align}\label{eq:backward-zn-weighted-small-final}
        \int_{\R^d}|z_n|^{p+1}\phi_n^2\,\ud x\to0
\end{align}
and
\begin{align}\label{eq:backward-selected-mass-positive}
        \sum_{j=1}^{J_*}m_j\geq\delta_0 .
\end{align}
Set $A_n:=W+S_n$. Since the fixed profile $W$ is asymptotically orthogonal
to the selected profiles by \eqref{eq:backward-selected-scale-small}, the same
weighted energy expansion as in the proof of Lemma~\ref{lem:short-time-Lq-prop}
gives
\begin{align}\label{eq:backward-positive-energy-defect-proof}
        E_{\phi_n}(u(\sigma_n))-E_{\phi_n}(W)
        &=
        E_{\phi_n}(W+v_n)-E_{\phi_n}(W) \notag\\
        &\geq
        \frac1d\sum_{j=1}^{J_*}
        \int_{\R^d}|\psi_n^j|^{p+1}\phi_n^2\,\ud x
        -o_n(1).
\end{align}
Using \eqref{eq:backward-selected-mass-positive}, we obtain
\begin{align}\label{eq:backward-positive-energy-defect}
        \liminf_{n\to\infty}
        \left(E_{\phi_n}(u(\sigma_n))-E_{\phi_n}(W)\right)
        \geq
        \frac{\delta_0}{d}>0 .
\end{align}
On the other hand, Lemma~\ref{lem:short-time-nonlinear-energy-prop} with
$\eta_n=\phi_n$ gives
$\lim_{n\to \infty}|E_{\phi_n}(u(\sigma_n))-E_{\phi_n}(u(\ta_n))|=0$. Together with
\eqref{eq:backward-Lp-prop-energy-input}, this implies
\begin{align}\label{eq:backward-energy-closeness-sigma}
        \lim_{n\to\infty}
        |E_{\phi_n}(u(\sigma_n))-E_{\phi_n}(W)|=0,
\end{align}
which contradicts \eqref{eq:backward-positive-energy-defect}. Therefore
\eqref{eq:backward-Lp-prop-conclusion} holds.
\end{proof}

Combining the backwards $L^{p+1}$ estimate with the localized nonlinear
energy estimate gives backwards-in-time propagation of the localized energy.
\begin{lem}[Short-time propagation of the localized energy backwards in time]
\label{lem:short-time-prop-backward}
Let $u(t)$ be a solution to \eqref{eqn:NLH} with initial data
$u(0)=u_0\in\dot H^1$. Let $T_+=T_+(u_0)$ denote its maximal time of
existence and assume that
\begin{align}\label{eq:short-time-prop-backward-energy-bound}
        \sup_{t\in[0,T_+)}\|u(t)\|_{\dot H^1}<\infty .
\end{align}
Let $0<\sigma_n<\ta_n<T_+$ satisfy $\sigma_n,\ta_n\to T_+$ and
$h_n:=\ta_n-\sigma_n\to0$ as $n\to\infty$. If $T_+=\infty$, this means
$\sigma_n,\ta_n\to\infty$. Let $W\in\dot H^1(\R^d)$ be a stationary solution
of \eqref{eq:yamabe}, possibly zero.

Let $r_n>0$ satisfy
\begin{align}\label{eq:short-time-prop-backward-rn}
        \frac{\ta_n-\sigma_n}{r_n^2}\to0
        \qquad\textrm{as }n\to\infty .
\end{align}
Assume
\begin{align}\label{eq:short-time-prop-backward-ball-input}
        \lim_{n\to \infty} \left[\bar E(u(\ta_n)-W;B(0,2r_n))
        +
        \|u(\ta_n)-W\|_{L^{p+1}(B(0,2r_n))}^{p+1}\right] = 0.
\end{align}
Then
\begin{align}\label{eq:en-discR-backward}
        \lim_{n\to \infty}\left[\bar E(u(\sigma_n)-W;B(0,r_n))+
        \|u(\sigma_n)-W\|_{L^{p+1}(B(0,r_n))}^{p+1}\right]
        =0.
\end{align}

Next, let $\veps_n>0$ satisfy $\veps_n<r_n$ and
\begin{align}\label{eq:short-time-prop-backward-epsn}
        \frac{\ta_n-\sigma_n}{\veps_n^2}\to0
        \qquad\textrm{as }n\to\infty .
\end{align}
Let $L\in\N$, $L\geq1$, and let $x_{\ell,n}\in\R^d$, $1\leq\ell\leq L$,
satisfy
\begin{align}\label{eq:short-time-prop-backward-puncture-sep}
        B(x_{\ell,n},\veps_n)\subset B(0,r_n),
        \qquad
        |x_{\ell,n}-x_{m,n}|\geq5\veps_n
        \quad\textrm{for }\ell\neq m .
\end{align}
Denote $J_n = B(0,2r_n)\setminus \cup_{\ell=1}^L B(x_{\ell,n},\veps_n/2)$ and assume that
\begin{align}\label{eq:short-time-prop-backward-annulus-input}
        \lim_{n\to\infty}
        \Bigl[
        &\bar E(u(\ta_n)-W; J_n)+
        \|u(\ta_n)-W\|_{L^{p+1}
        (J_n)}^{p+1}
        \Bigr]
        =0.
\end{align}
Then on $I_n = B(0,r_n)\setminus
        \cup_{\ell=1}^L B(x_{\ell,n},\veps_n)$ we have
\begin{align}\label{eq:en-annulus-backward}
        \lim_{n\to\infty}
        \Bigl[
        &\bar E(u(\sigma_n)-W; I_n)
        +
        \|u(\sigma_n)-W\|_{L^{p+1}
        (I_n)}^{p+1}
        \Bigr]
        =0.
\end{align}
\end{lem}

\begin{proof}
We prove the ball and punctured ball case simultaneously. In the ball case,
set $\Omega_n^0:=B(0,r_n)$, $\Omega_n^1:=B(0,3r_n/2)$,
$\Omega_n^2:=B(0,2r_n)$, and $\ell_n:=r_n$. In the punctured ball case,
set $\ell_n:=\veps_n$,
\begin{align}\label{eq:short-time-backward-punctured-regions}
        \Omega_n^0
        &:=
        B(0,r_n)\setminus
        \bigcup_{\ell=1}^L \overline{B(x_{\ell,n},\veps_n)},
        \quad
        \Omega_n^1
        :=
        B(0,3r_n/2)\setminus
        \bigcup_{\ell=1}^L \overline{B(x_{\ell,n},3\veps_n/4)},
        \notag\\
        \Omega_n^2
        &:=
        B(0,2r_n)\setminus
        \bigcup_{\ell=1}^L \overline{B(x_{\ell,n},\veps_n/2)}.
\end{align}
Then $\Omega_n^0\subset\Omega_n^1\subset\Omega_n^2$ and
$(\ta_n-\sigma_n)/\ell_n^2\to0$, as $n\to \infty.$ Using
\eqref{eq:short-time-prop-backward-puncture-sep} in the punctured ball case,
choose $\phi_n,\psi_n\in C_c^\infty(\Omega_n^2)$ such that
\begin{align}\label{eq:short-time-backward-main-cutoffs}
        0\leq\phi_n,\psi_n\leq1,
        \qquad
        \phi_n\equiv1\textrm{ on }\Omega_n^1,
        \qquad
        \psi_n\equiv1\textrm{ on }\Omega_n^0,
        \qquad
        \supp\psi_n\subset\Omega_n^1,
\end{align}
and
\begin{align}\label{eq:short-time-backward-main-cutoff-bounds}
        \|\nabla\phi_n\|_{L^\infty}
        +
        \|\nabla\psi_n\|_{L^\infty}
        \lesssim\ell_n^{-1},
        \qquad
        \|\nabla\phi_n\|_{L^d}
        +
        \|\nabla\psi_n\|_{L^d}
        \lesssim1 .
\end{align}
Let $\eta_n=\phi_n$ or $\eta_n=\psi_n$, and set
$e_n:=u(\ta_n)-W$. Since $\supp\eta_n\subset\Omega_n^2$ and
$0\leq\eta_n\leq1$,
\begin{align}\label{eq:short-time-backward-final-E}
        \left|E_{\eta_n}(u(\ta_n))-E_{\eta_n}(W)\right|
        &\leq
        \frac12\int_{\Omega_n^2}|\nabla e_n|^2\,\ud x
        +
        \int_{\Omega_n^2}|\nabla W||\nabla e_n|\,\ud x \notag\\
        &\quad+
        \frac1{p+1}\int_{\Omega_n^2}
        \left||W+e_n|^{p+1}-|W|^{p+1}\right|\,\ud x \notag\\
        &\lesssim
        \|\nabla e_n\|_{L^2(\Omega_n^2)}^2
        +
        \|\nabla W\|_{L^2}\|\nabla e_n\|_{L^2(\Omega_n^2)} \notag\\
        &\quad+
        \|W\|_{L^{p+1}}^p\|e_n\|_{L^{p+1}(\Omega_n^2)}
        +
        \|e_n\|_{L^{p+1}(\Omega_n^2)}^{p+1}.
\end{align}
The right-hand side tends to $0$ by
\eqref{eq:short-time-prop-backward-ball-input} in the ball case and by
\eqref{eq:short-time-prop-backward-annulus-input} in the punctured ball case.
Therefore, using Lemma~\ref{lem:short-time-nonlinear-energy-prop} with
$\eta_n=\phi_n$ and $\eta_n=\psi_n$, we obtain
\begin{align}\label{eq:short-time-backward-sigma-E}
        E_{\eta_n}(u(\sigma_n))-E_{\eta_n}(W)\to0
        \qquad\textrm{as }n\to\infty .
\end{align}
Applying Lemma~\ref{lem:short-time-Lp-prop-backward} with
$\Omega_n^1$, $\Omega_n^2$, $\ell_n$, and $\phi_n$, we obtain
\begin{align}\label{eq:short-time-backward-main-Lp}
        \|u(\sigma_n)-W\|_{L^{p+1}(\Omega_n^1)}\to0
        \qquad\textrm{as }n\to\infty .
\end{align}
Set $w_n:=u(\sigma_n)-W$. By
\eqref{eq:short-time-backward-sigma-E} with $\eta_n=\psi_n$,
\begin{align}\label{eq:short-time-backward-main-Epsi}
        E_{\psi_n}(W+w_n)-E_{\psi_n}(W)\to0 .
\end{align}
Expanding,
\begin{align}\label{eq:short-time-backward-main-expand}
        E_{\psi_n}(W+w_n)-E_{\psi_n}(W)
        &=
        \frac12\int_{\R^d}|\nabla w_n|^2\psi_n^2\,\ud x
        +
        \int_{\R^d}\nabla W\cdot\nabla w_n\,\psi_n^2\,\ud x \notag\\
        &\quad
        -\frac1{p+1}\int_{\R^d}
        \left(|W+w_n|^{p+1}-|W|^{p+1}\right)\psi_n^2\,\ud x .
\end{align}
Using \eqref{eq:yamabe} for $W$ gives
\begin{align}\label{eq:short-time-backward-main-cross-W}
        \int_{\R^d}\nabla W\cdot\nabla w_n\,\psi_n^2\,\ud x
        =
        \int_{\R^d}|W|^{p-1}Ww_n\psi_n^2\,\ud x
        -
        2\int_{\R^d}w_n\psi_n\nabla W\cdot\nabla\psi_n\,\ud x .
\end{align}
Because $\supp\psi_n\subset\Omega_n^1$ and
\eqref{eq:short-time-backward-main-Lp} holds,
\begin{align}\label{eq:short-time-backward-main-Taylor-zero}
        \int_{\R^d}
        \left|
        |W+w_n|^{p+1}-|W|^{p+1}
        -(p+1)|W|^{p-1}Ww_n
        \right|\psi_n^2\,\ud x
        \to0 .
\end{align}
Also,
\begin{align}\label{eq:short-time-backward-main-cutoff-W-zero}
        \left|
        \int_{\R^d}w_n\psi_n\nabla W\cdot\nabla\psi_n\,\ud x
        \right|
        \leq
        \|w_n\|_{L^{p+1}(\Omega_n^1)}
        \|\nabla W\|_{L^2}
        \|\nabla\psi_n\|_{L^d}
        \to0 .
\end{align}
Combining \eqref{eq:short-time-backward-main-Epsi},
\eqref{eq:short-time-backward-main-expand},
\eqref{eq:short-time-backward-main-cross-W},
\eqref{eq:short-time-backward-main-Taylor-zero}, and
\eqref{eq:short-time-backward-main-cutoff-W-zero}, we obtain
\begin{align}\label{eq:short-time-backward-main-grad-small}
        \lim_{n\to \infty}\int_{\R^d}|\nabla w_n|^2\psi_n^2\,\ud x = 0,
\end{align}
which implies
\begin{align}\label{eq:short-time-backward-main-Omega0}
        \lim_{n\to \infty}\bar E(u(\sigma_n)-W;\Omega_n^0)=0.
\end{align}
Combining \eqref{eq:short-time-backward-main-Omega0} with
\eqref{eq:short-time-backward-main-Lp} gives the desired conclusions.
\end{proof}
We next use Lemma~\ref{lem:short-time-prop} to show that the singular set in Theorem~\ref{thm:lwp} is nonempty when $T_+<\infty$.
\begin{lem}[Non-empty singular set when $T_+<\infty$]
\label{lem:finite-time-no-escape}
Let $u$ solve \eqref{eqn:NLH} with maximal time $T_+<\infty$ and
\begin{align}\label{eq:no-escape-energy-bound}
        \sup_{0\leq t<T_+}\bar E(u(t))<\infty .
\end{align}
Then
\begin{align}\label{eq:no-escape-conclusion}
        \lim_{R\to\infty}\limsup_{t\to T_+}
        \bar E(u(t);\R^d\setminus B(0,R))=0 .
\end{align}
Consequently, in Theorem~\ref{thm:lwp}, the singular set is not empty, i.e. $L\geq1$.
\end{lem}
\begin{proof}
We first prove \eqref{eq:no-escape-conclusion}. Suppose not. Then there exists
$\delta_0>0$ such that, for every $R>0$ and every $0<\sigma<T_+$, there is
$\ta\in(\sigma,T_+)$ with
\begin{align}\label{eq:no-escape-failure}
        \bar E(u(\ta);\R^d\setminus B(0,R))\geq\delta_0 .
\end{align}
Choose $\sigma_n\to T_+$. Since $u(\sigma_n)\in\dot H^1(\R^d)\subset L^{p+1}(\R^d)$, we may choose $R_n\to\infty$ such that
\begin{align}\label{eq:no-escape-initial-tail}
        \bar E(u(\sigma_n);\R^d\setminus B(0,R_n/2))
        +
        \|u(\sigma_n)\|_{L^{p+1}(\R^d\setminus B(0,R_n/2))}^{p+1}
        \to0
        \qquad\textrm{as }n\to\infty.
\end{align}
We may further choose the sequence $R_n$ such that
\begin{align}\label{eq:no-escape-short-time-scale}
        \frac{T_+-\sigma_n}{R_n^2}\to0
        \qquad\textrm{as }n\to\infty .
\end{align}
By \eqref{eq:no-escape-failure}, after passing to a subsequence if necessary, choose
$\ta_n\in(\sigma_n,T_+)$ such that
\begin{align}\label{eq:no-escape-final-tail-lower}
        \bar E(u(\ta_n);\R^d\setminus B(0,R_n))\geq\delta_0 .
\end{align}
For each $n$, since $u(\ta_n)\in\dot H^1(\R^d)$, choose $S_n>R_n$ so large that
\begin{align}\label{eq:no-escape-final-far-tail}
        \bar E(u(\ta_n);\R^d\setminus B(0,S_n))\leq n^{-1}.
\end{align}
Increasing $S_n$ if necessary, we may assume $S_n>2R_n$.

We apply Lemma~\ref{lem:short-time-prop} with $W=0$, $L=1$, $r_n:=S_n$, $\veps_n:=R_n$, and $x_{1,n}:=0$. Moreover,
$B(0,2S_n)\setminus B(0,R_n/2)\subset\R^d\setminus B(0,R_n/2)$, so
\eqref{eq:no-escape-initial-tail} implies \eqref{eq:short-time-prop-annulus-input}. Therefore, Lemma~\ref{lem:short-time-prop} gives
\begin{align}\label{eq:no-escape-propagated-annulus}
        \bar E(u(\ta_n);B(0,S_n)\setminus B(0,R_n))\to0
        \qquad\textrm{as }n\to\infty .
\end{align}
Combining \eqref{eq:no-escape-final-far-tail} and
\eqref{eq:no-escape-propagated-annulus}, we obtain
\begin{align}\label{eq:no-escape-final-tail-upper}
        \bar E(u(\ta_n);\R^d\setminus B(0,R_n))
        &\leq
        \bar E(u(\ta_n);B(0,S_n)\setminus B(0,R_n))
        +
        \bar E(u(\ta_n);\R^d\setminus B(0,S_n))\to0,
\end{align}
as $n\to \infty.$ This contradicts \eqref{eq:no-escape-final-tail-lower}. Hence
\eqref{eq:no-escape-conclusion} holds.

It remains to prove $L\geq1$. Suppose, toward a contradiction, that $L=0$ in
Theorem~\ref{thm:lwp}. Then $u$ converges strongly to the weak limit $u^*$ on every compact
subset of $\R^d$. By \eqref{eq:no-escape-conclusion}, for every $\eta>0$ there is
$R>0$ such that
\begin{align}\label{eq:no-escape-tail-u}
        \limsup_{t\to T_+}\bar E(u(t);\R^d\setminus B(0,R))\leq\eta .
\end{align}
Since $u^*\in\dot H^1(\R^d)$, increasing $R$ if necessary gives
\begin{align}\label{eq:no-escape-tail-ustar}
        \bar E(u^*;\R^d\setminus B(0,R))\leq\eta .
\end{align}
Using the convergence in $C^1(\overline{B(0,R)})$, \eqref{eq:no-escape-tail-u}, and
\eqref{eq:no-escape-tail-ustar}, we obtain
\begin{align}\label{eq:no-escape-global-strong}
        \limsup_{t\to T_+}\|\nabla(u(t)-u^*)\|_{L^2}^2
        &\leq
        \limsup_{t\to T_+}\|\nabla(u(t)-u^*)\|_{L^2(B(0,R))}^2 \notag\\
        &\quad+
        2\limsup_{t\to T_+}\bar E(u(t);\R^d\setminus B(0,R))
        +
        2\bar E(u^*;\R^d\setminus B(0,R))\notag\\
        &\leq 4\eta .
\end{align}
Letting $\eta\to0$ yields
\begin{align}\label{eq:no-escape-strong-H1-limit}
        u(t)\to u^*
        \qquad\textrm{strongly in }\dot H^1(\R^d)
        \qquad\textrm{as }t\to T_+ .
\end{align}
The local well-posedness theorem applied with initial data $u^*$ at time
$T_+$ gives a solution on $[T_+,T_++\delta)$ for some $\delta>0$. By the continuation criterion in Theorem~\ref{thm:lwp}, $u$ extends past $T_+$, contradicting the maximality of $T_+$. Therefore $L\geq1$.
\end{proof}

\subsection{Concentration properties of the heat flow}
The goal of this subsection is to prove that energy cannot concentrate outside the self-similar scale, as expected in a Type-II blowup scenario. Similar results are known for other energy-critical equations, including the energy-critical nonlinear wave equation \cite{jia-dkm}, wave maps \cite{christodoulou1993regularity,shath1992regularity}, and the harmonic map heat flow \cite{lawrie-harmonic-map-nonradial}. Because \eqref{eqn:NLH} has no finite speed of propagation, the hyperbolic arguments do not apply directly. Moreover, the nonlinear energy for \eqref{eqn:NLH} is not coercive, so the standard energy arguments for the harmonic map heat flow are not directly available. We will instead use profile decompositions. We first treat the finite-time blowup case.

\begin{lem}[No self-similar nonlinear energy concentration in the finite-time blowup case]
\label{lem:ss-bu-nonlinear}
Let $u(t)$ be a solution of~\eqref{eqn:NLH} with initial data
$u_0\in\dot H^1$, $T_+=T_+(u_0)<\infty$, and
$\sup_{t\in[0,T_+)}\bar E(u(t))<\infty$. Let $x_0\in\calS$ be a
singular point as in \eqref{defn:singular-points}, and let $r>0$ be
sufficiently small so that
$B(x_0,2r)\cap(\calS\setminus\{x_0\})=\emptyset$. Set
$\rho(t):=\sqrt{T_+-t}$. For every $0<a<b<\infty$ and every
$\Theta\in C_c^\infty(B(0,b)\setminus B(0,a))$ satisfying
$0\leq\Theta\leq1$, define $ \Theta_t(x):= \Theta((x-x_0)/{\rho(t)}).$ Then
\begin{align}\label{eqn:self-similar-annulus-vanish}
        \lim_{t\to T_+}E_{\Theta_t}(u(t))=0.
\end{align}
\end{lem}

\begin{proof}
Fix $s<T_+$ close to $T_+$. By \eqref{eqn:energy ineq I}, for
$s<\tau<T_+$,
\begin{align}\label{eqn:self-similar-energy-cauchy}
\begin{split}
\left|
E_{\Theta_s}(u(\tau))
-
E_{\Theta_s}(u(s))
\right|
&\lesssim
\int_s^\tau\int_{\R^d}|\partial_tu|^2\ud x\ud t\\
&\quad+
\left(\int_s^\tau\int_{\R^d}|\partial_tu|^2\ud x\ud t\right)^{1/2}
\left(\int_s^\tau\int_{\R^d}
|\nabla u|^2|\nabla\Theta_s|^2\ud x\ud t\right)^{1/2}.
\end{split}
\end{align}
Since $|\nabla\Theta_s|\lesssim\rho(s)^{-1}$,
$\sup_{0\leq t<T_+}\bar E(u(t))<\infty$, and
\eqref{eq:tension-L2} holds, we obtain
\begin{align}\label{eqn:self-similar-cauchy}
        \sup_{\tau\in(s,T_+)}
        \left|
        E_{\Theta_s}(u(\tau))
        -
        E_{\Theta_s}(u(s))
        \right|
        =o_s(1).
\end{align}
For $s$ sufficiently close to $T_+$,
$\supp\Theta_s\Subset B(x_0,r)\setminus\calS$. Letting
$\tau\to T_+$ in \eqref{eqn:self-similar-cauchy} and using
Theorem~\ref{thm:lwp}, we obtain
\begin{align}
        E_{\Theta_s}(u(s))
        =
        E_{\Theta_s}(u^*)+o_s(1).
\end{align}
Since $E_{\Theta_s}(u^*)\to0$ as $s\to T_+$,
\eqref{eqn:self-similar-annulus-vanish} follows.
\end{proof}
We next show that no self-similar concentration of
$L^{p+1}$-norm occurs. This argument makes use of profile decompositions.
\begin{lem}[No self-similar $L^{p+1}$ concentration in the finite-time blowup case]
\label{lem:ss-bu-lp}
Under the hypotheses of Lemma~\ref{lem:ss-bu-nonlinear}, for every $\alp>0$,
\begin{align}\label{eqn:no-conc-finite-lp}
\lim_{t\to T_+}
\int_{B(x_0,r)\setminus B(x_0,\alp\sqrt{T_+-t})}|u(t)|^{p+1}\ud x
=
\int_{B(x_0,r)}|u^*|^{p+1}\ud x.
\end{align}
\end{lem}
\begin{proof}
After translating, assume $x_0=0$. It suffices to prove
\begin{align}\label{eqn:finite-shrinking-neck}
\lim_{\veps\to 0}\limsup_{t\to T_+}
\int_{B(0,\veps)\setminus B(0,\alp\rho(t))}|u(t)|^{p+1}\ud x=0,
\end{align}
where $\rho(t):=\sqrt{T_+-t}$.  Suppose the limit in \eqref{eqn:finite-shrinking-neck} fails. Then there exist $\delta>0$, $\veps_n\to0$, and $t_n\to T_+$ such that, with $\rho_n:=\sqrt{T_+-t_n}$ and $\eta_n:={\rho_n}/{\veps_n}$, we have $\eta_n\to0$ as $n\to\infty$ and
\begin{align}\label{eqn:finite-neck-lower}
\int_{B(0,\veps_n)\setminus B(0,\alp\rho_n)}|u(t_n)|^{p+1}\ud x\geq\delta.
\end{align}
The choice $\eta_n\to0$ follows from the definition of $\limsup$, since $t_n$ may be chosen sufficiently close to $T_+$. Since $t_n\to T_+$, with $T_+>0$ and $\veps_n\to0$, for $n\gg 1$ we have $t_n-\veps_n^2>0$. Define
\begin{align}\label{eqn:finite-rescaled-flow}
v_n(y,s):=\veps_n^{\frac{d-2}{2}}u(t_n+\veps_n^2s,\veps_n y),
\qquad
-1\leq s\leq0,
\end{align}
and set $w_n(y):=v_n(y,0)$. Then $\{w_n\}_{n\in \N}$ is bounded in $\dot H^1(\R^d)$ and
\begin{align}\label{eqn:finite-rescaled-tension}
\int_{-1}^0\|\partial_sv_n(s)\|_{L^2_y}^2\ud s
=
\int_{t_n-\veps_n^2}^{t_n}\|\partial_tu(t)\|_{L^2_x}^2\ud t
\to0
\qquad\textrm{as }n\to\infty .
\end{align}
Moreover,
\begin{align}\label{eqn:finite-neck-lower-rescaled}
\int_{B(0,1)\setminus B(0,\alp\eta_n)}|w_n|^{p+1}\ud y\geq\delta.
\end{align}
By Proposition \ref{prop:profile-decomp} we obtain the following profile decomposition of $w_n$,
\begin{align}\label{eqn:finite-profile-decomp}
w_n=\sum_{j=1}^{J}W_n^j+r_n^J,
\quad
W_n^j(y):=(\lambda_{j,n})^{-\frac{d-2}{2}}
W^j\left(\frac{y-a_{j,n}}{\lambda_{j,n}}\right),
\end{align}
with pairwise orthogonal parameters and
\begin{align}\label{eqn:finite-profile-rem}
\lim_{J\to\infty}\limsup_{n\to\infty}\|r_n^J\|_{L^{p+1}}=0.
\end{align}
We now claim that every nontrivial profile with $\limsup_n\lambda_{j,n}<\infty$ is stationary. To see this, fix such a profile and write $\lambda_n=\lambda_{j,n}$, $a_n=a_{j,n}$, and $W=W^j$. Choose $\sigma_0>0$ such that $\sigma_0\lambda_n^2\leq1$ for all large $n$, and define
\begin{align}\label{eqn:finite-profile-flow}
V_n(z,\sigma):=
\lambda_n^{\frac{d-2}{2}}
v_n(a_n+\lambda_n z,\lambda_n^2\sigma),
\qquad
-\sigma_0\leq\sigma\leq0.
\end{align}
Since $v_n$ solves \eqref{eqn:NLH} on $\R^d\times[-1,0]$ and $\sigma_0\lambda_n^2\leq1$, the parabolic scale invariance of \eqref{eqn:NLH} implies that $V_n$ solves $\partial_\sigma V_n=\Delta_zV_n+|V_n|^{p-1}V_n$ on $\R^d\times[-\sigma_0,0]$. Furthermore,
\begin{align}\label{eqn:finite-profile-tension}
\int_{-\sigma_0}^0\|\partial_\sigma V_n(\sigma)\|_{L^2_z}^2\ud\sigma\to0
\qquad\textrm{as }n\to\infty .
\end{align}
By \eqref{eq:profile-remainder-frame-vanishing} and parameter orthogonality,
\begin{align}
V_n(0)\rightharpoonup W
\qquad\textrm{weakly in }\dot H^1(\R^d).
\end{align}
Hence, after passing to a subsequence, $V_n(0)\to W$ strongly in $L^2_{\loc}(\R^d)$ by Rellich's theorem.
% Also $V_n(0)\to W$ strongly in $L^2_{\loc}(\R^d)$. This follows from the fact that $\lam_{j,n}^{\frac{d-2}{2}}W_n^j(a_{j,n}+\lam_{j,n}\cdot)\to W^j$ strongly in $L^2_{\loc}$ while  for each $k\neq j$, $\lam_{j,n}^{\frac{d-2}{2}}W_n^k(a_{j,n}+\lam_{j,n}\cdot) \to 0$ strongly in $L^2_{\loc}$ as $n\to \infty$ by orthogonality of parameters, and the remainder under the same rescaling converges weakly to zero in $\dot H^1$ and hence converges strongly to zero in $L^2_{\loc}$ by Rellich's theorem as $n\to \infty$. 
Thus, for every compact $K\Subset\R^d$,
\begin{align}
\sup_{-\sigma_0\leq\sigma\leq0}
\|V_n(\sigma)-W\|_{L^2(K)}
\leq
\|V_n(0)-W\|_{L^2(K)}
+
\sigma_0^{1/2}\|\partial_\sigma V_n\|_{L^2_\sigma L^2_z}
\to0
\qquad\textrm{as }n\to\infty .
\end{align}
Since $V_n\to W$ strongly in $L^2_{\loc}$ and $\{V_n\}_{n\in \N}$ is uniformly bounded in $L^\infty_\sigma\dot H^1_z$, Sobolev's inequality gives a uniform $L^\infty_\sigma L^{p+1}_z$ bound. Hence $V_n\to W$ strongly in $L^{q}_{\loc}$ for every $2\leq q<p+1$. On compact sets, the same conclusion for $1\leq q<2$ follows from Hölder's inequality. In particular, $V_n\to W$ strongly in $L^p_{\loc}$. Furthermore, the sequence $V_n$ is bounded in $L^\infty([-\sigma_0,0];\dot H^1(\R^d))$. Hence, after passing to a subsequence,
\begin{align}\label{eq:Vn-weak-gradient-finite-time}
        \nabla V_n\rightharpoonup \nabla W
        \qquad\textrm{weakly in }L^2_{\loc}((-\sigma_0,0)\times\R^d).
\end{align}
Taking the weak formulation of \eqref{eqn:NLH} satisfied by $V_n$ against $\theta(\sigma)\zeta(z)$, where $\theta\in C_c^\infty((-\sigma_0,0))$, $\int\theta\ud\sigma=1$, and $\zeta\in C_c^\infty(\R^d)$, and then using \eqref{eqn:finite-profile-tension}, we obtain
\begin{align}\label{eqn:finite-profile-stationary}
\int_{\R^d}\nabla W\cdot\nabla\zeta\ud z
=
\int_{\R^d}|W|^{p-1}W\zeta\ud z
\qquad
\textrm{for all }\zeta\in C_c^\infty(\R^d).
\end{align}
Thus $W$ solves \eqref{eq:yamabe}.  

After passing to a diagonal subsequence, we may assume that, for every $j\in\N$, the limit of $\lambda_{j,n}$ exists in $[0,\infty]$. Discarding the trivial profiles, every profile for which $\lim_n\lambda_{j,n}<\infty$ is a nontrivial stationary solution by the preceding argument. Since every nontrivial stationary solution has Dirichlet energy at least $\bar E_*$ and
\begin{align}
        \sum_{j=1}^{\infty}\|\nabla W^j\|_{L^2}^2
        \leq
        \limsup_{n\to\infty}\|w_n\|_{\dot H^1}^2<\infty,
\end{align}
there are only finitely many such profiles. After reordering, denote them by $W^1,\ldots,W^N$. Thus $\lambda_{j,n}\to\infty$ for every $j>N$. For every fixed $R<\infty$ and every $j>N$, note that
\begin{align}\label{eqn:finite-large-scale-profile-vanish}
        \|W_n^j\|_{L^{p+1}(B(0,R))}^{p+1}
        =
        \int_{B(-a_{j,n}/\lambda_{j,n},R/\lambda_{j,n})}
        |W^j(z)|^{p+1}\,\ud z
        \to0,
\end{align}
as $n\to \infty.$ Set $S_n:=\sum_{j=1}^N W_n^j$ and $z_n:=w_n-S_n$. For every fixed $R<\infty$ and every fixed $J\geq N+1$,
\begin{align}
        \|z_n\|_{L^{p+1}(B(0,R))}
        \leq
        \sum_{j=N+1}^J
        \|W_n^j\|_{L^{p+1}(B(0,R))}
        +
        \|r_n^J\|_{L^{p+1}}.
\end{align}
Sending $n\to \infty$, using \eqref{eqn:finite-large-scale-profile-vanish}, and then letting $J\to\infty$ in \eqref{eqn:finite-profile-rem}, we obtain for every fixed $R<\infty$
\begin{align}\label{eqn:finite-profile-remainder-local-small}
        \|z_n\|_{L^{p+1}(B(0,R))}\to0,
\end{align}
as $n\to \infty.$ In particular,
\begin{align}\label{eqn:finite-local-profile-approx}
        \|w_n-S_n\|_{L^{p+1}(B(0,R))}\to0
\end{align}
as $n\to \infty$. From \eqref{eqn:finite-neck-lower-rescaled}, we see that $N\geq1$. Choose $R>e$. For all sufficiently large $n$, let $\Phi_{n,R}$ be radial and satisfy $0\leq\Phi_{n,R}\leq1$, $\Phi_{n,R}=0$ for $|y|\leq\frac{\alp\eta_n}{R}$, $\Phi_{n,R}=1$ for $\alp\eta_n\leq |y|\leq1$, $\Phi_{n,R}=0$ for $|y|\geq R$, and $|\nabla\Phi_{n,R}(y)|\lesssim\frac1{|y|\log R}$. Then,
\begin{align}\label{eqn:finite-log-Ld}
        \|\nabla\Phi_{n,R}\|_{L^d}
        \lesssim
        (\log R)^{-\frac{d-1}{d}}.
\end{align}
By \eqref{eqn:finite-neck-lower-rescaled} and \eqref{eqn:finite-local-profile-approx}, for all sufficiently large $n$,
\begin{align}\label{eqn:finite-Sn-mass-lower}
        \int_{\R^d}|S_n|^{p+1}\Phi_{n,R}^2\,\ud y
        \geq
        c_\delta
\end{align}
for some $c_\delta>0$ depending only on $\delta$. Since
$|S_n|^{p+1}\leq N^p\sum_{j=1}^N|W_n^j|^{p+1}$, it follows that
\begin{align}\label{eqn:finite-profile-weighted-mass}
        \sum_{j=1}^{N}
        \int_{\R^d}|W_n^j|^{p+1}\Phi_{n,R}^2\,\ud y
        \geq
        m_0,
        \qquad
        m_0:=N^{-p}c_\delta>0.
\end{align}
Since $\supp\Phi_{n,R}\subset B(0,R)$, \eqref{eqn:finite-profile-remainder-local-small} allows us to apply Lemma~\ref{lem:weighted-profile-decoupling} to $w_n=S_n+z_n$, which gives
\begin{align}\label{eqn:finite-weighted-energy-expansion}
        E_{\Phi_{n,R}}(w_n)
        \geq
        \sum_{j=1}^{N}E_{\Phi_{n,R}}(W_n^j)-o_n(1).
\end{align}
Testing \eqref{eq:yamabe} for $W_n^j$ with $W_n^j\Phi_{n,R}^2$ gives
\begin{align}
        E_{\Phi_{n,R}}(W_n^j)
        =
        \frac1d\int_{\R^d}|W_n^j|^{p+1}\Phi_{n,R}^2\,\ud y
        -
        \int_{\R^d}W_n^j\Phi_{n,R}\nabla W_n^j\cdot\nabla\Phi_{n,R}\,\ud y.
\end{align}
Moreover,
\begin{align}
        \left|
        \int_{\R^d}W_n^j\Phi_{n,R}\nabla W_n^j\cdot\nabla\Phi_{n,R}\,\ud y
        \right|
        \leq
        \|W^j\|_{L^{p+1}}
        \|\nabla W^j\|_{L^2}
        \|\nabla\Phi_{n,R}\|_{L^d}.
\end{align}
Choose $R>e$ sufficiently large that
\begin{align}\label{eqn:finite-log-error-small}
        \sum_{j=1}^{N}
        \|W^j\|_{L^{p+1}}
        \|\nabla W^j\|_{L^2}
        \|\nabla\Phi_{n,R}\|_{L^d}
        \leq
        \frac{m_0}{2d}.
\end{align}
Then \eqref{eqn:finite-profile-weighted-mass}, \eqref{eqn:finite-weighted-energy-expansion}, and \eqref{eqn:finite-log-error-small} imply
\begin{align}\label{eqn:finite-positive-local-energy}
        \liminf_{n\to\infty}E_{\Phi_{n,R}}(w_n)
        \geq
        \frac{m_0}{2d}>0.
\end{align}
Set $\phi_{n,R}(x):=\Phi_{n,R}(x/\veps_n)$. By scale invariance,
\begin{align}\label{eqn:finite-energy-scaling}
E_{\phi_{n,R}}(u(t_n))=E_{\Phi_{n,R}}(w_n).
\end{align}
We claim that
\begin{align}\label{eqn:finite-energy-zero-contr}
E_{\phi_{n,R}}(u(t_n))\to0
\qquad\textrm{as }n\to\infty .
\end{align}
For $t_n<\tau<T_+$, \eqref{eqn:energy ineq I} gives
\begin{align}
\begin{split}
|E_{\phi_{n,R}}(u(\tau))-E_{\phi_{n,R}}(u(t_n))|
&\lesssim
\int_{t_n}^{\tau}\|\partial_tu(t)\|_{L^2}^2\ud t\\
&\quad+
\left(\int_{t_n}^{\tau}\|\partial_tu(t)\|_{L^2}^2\ud t\right)^{1/2}
\left(\int_{t_n}^{\tau}\int_{\R^d}|\nabla u|^2|\nabla\phi_{n,R}|^2\ud x\ud t\right)^{1/2}.
\end{split}
\end{align}
Since $(T_+-t_n)\|\nabla\phi_{n,R}\|_{L^\infty}^2\lesssim_{\alp,R}1$ and \eqref{eq:tension-L2}, we see that $|E_{\phi_{n,R}}(u(\tau))-E_{\phi_{n,R}}(u(t_n))| = o_n(1)$ uniformly in $\tau\in(t_n,T_+)$. For all sufficiently large $n$, we have $\supp\phi_{n,R}\subset\overline{B(0,R\veps_n)}\setminus\{0\}$, and hence $\supp\phi_{n,R}\cap\calS=\emptyset$. Letting $\tau\to T_+$ gives
\begin{align}
E_{\phi_{n,R}}(u(t_n))=E_{\phi_{n,R}}(u^*)+o_n(1).
\end{align}
Since $\supp\phi_{n,R}\subset B(0,R\veps_n)$, the absolute continuity of the $\dot H^1$ and $L^{p+1}$ norms of $u^*$ gives $E_{\phi_{n,R}}(u^*)\to0$ as $n\to \infty.$ Thus \eqref{eqn:finite-energy-zero-contr} follows, contradicting \eqref{eqn:finite-positive-local-energy}. Hence \eqref{eqn:finite-shrinking-neck} holds.
\end{proof}
The preceding two lemmas yield the following main result in the finite-time blowup case.
\begin{lem}[No self-similar energy concentration in the finite-time blowup case]
\label{lem:ss-bu}
Under the hypotheses of Lemma~\ref{lem:ss-bu-nonlinear}, for every $\alp>0$,
\begin{align}\label{eqn:no-conc-finite}
\lim_{t\to T_+}
\bar E\left(u(t);B(x_0,r)\setminus B(x_0,\alp\sqrt{T_+-t})\right)
=
\bar E(u^*;B(x_0,r)).
\end{align}
In particular, there exist $T_0<T_+$ and functions $\nu,\xi:[T_0,T_+)\to(0,\infty)$ such that
\begin{align}\label{eqn:nu-xi-curves}
&\lim_{t\to T_+}\left(
\frac{\xi(t)}{\sqrt{T_+-t}}
+
\frac{\sqrt{T_+-t}}{\nu(t)}
\right)=0,\\
&\lim_{t\to T_+} \left[\bar{E}\left(u(t);B(x_0,\nu(t))\setminus B(x_0,\xi(t))\right) + \|u(t)\|^{p+1}_{L^{p+1}(B(x_0,\nu(t))\setminus B(x_0,\xi(t)))}\right]=0.
\end{align}
\end{lem}
\begin{proof}
Set $\rho(t):=\sqrt{T_+-t}$. Fix $0<a<b<\infty$, and choose a radial
$\Theta\in C_c^\infty(B(0,2b)\setminus B(0,a/2))$ such that
$0\leq\Theta\leq1$ and
$\Theta\equiv1$ on $B(0,b)\setminus B(0,a)$. Set
$\Theta_t(x):=\Theta((x-x_0)/\rho(t))$. Subtracting \eqref{eqn:no-conc-finite-lp} with $\alp=2b$ from the same identity with $\alp=a/2$ gives
\begin{align}\label{eqn:self-similar-Lp-annulus-vanish}
        \lim_{t\to T_+}
        \int_{B(x_0,2b\rho(t))\setminus B(x_0,a\rho(t)/2)}
        |u(t)|^{p+1}\ud x
        =0.
\end{align}
Consequently, Lemma~\ref{lem:ss-bu-nonlinear} gives
\begin{align}
        \bar E\left(
        u(t);
        B(x_0,b\rho(t))\setminus B(x_0,a\rho(t))
        \right)
        &\leq
        \bar E_{\Theta_t}(u(t)) \notag\\
        &=
        2E_{\Theta_t}(u(t))
        +
        \frac{2}{p+1}
        \int_{\R^d}|u(t)|^{p+1}\Theta_t^2\ud x
        \to0
\end{align}
as $t\to T_+.$ Thus
\begin{align}\label{eqn:self-similar-positive-annulus-vanish}
        \lim_{t\to T_+}
        \left[
        \bar E\left(
        u(t);
        B(x_0,b\rho(t))\setminus B(x_0,a\rho(t))
        \right)
        +
        \int_{B(x_0,b\rho(t))\setminus B(x_0,a\rho(t))}
        |u(t)|^{p+1}\ud x
        \right]
        =0.\quad 
\end{align}
Fix $\alp>0$. For $0<\veps<r$, choose
$\chi_\veps\in C_c^\infty(B(x_0,r+\veps))$ such that
$0\leq\chi_\veps\leq1$ and
$\chi_\veps\equiv1$ on $B(x_0,r)$. Choose a radial
$\eta\in C^\infty(\R^d)$ such that $0\leq\eta\leq1$,
$\eta(z)=0$ for $|z|\leq\alp/2$, and
$\eta(z)=1$ for $|z|\geq\alp$, and set
\begin{align}
        \phi_{t,\veps}(x):=
        \chi_\veps(x)
        \eta\left(\frac{x-x_0}{\rho(t)}\right).
\end{align}
Since
\begin{align}
        (T_+-t)\|\nabla\phi_{t,\veps}\|_{L^\infty}^2
        \lesssim_{\alp,\veps}1,
\end{align}
the argument in the proof of Lemma~\ref{lem:ss-bu-nonlinear}, with
$\phi_{t,\veps}$ frozen at time $t$, gives
\begin{align}\label{eqn:finite-main-cutoff-energy}
        \lim_{t\to T_+}E_{\phi_{t,\veps}}(u(t))
        =
        E_{\chi_\veps}(u^*).
\end{align}
Moreover, Lemma~\ref{lem:ss-bu-lp},
\eqref{eqn:self-similar-Lp-annulus-vanish}, and the strong convergence
away from $\calS$ give
\begin{align}\label{eqn:finite-main-cutoff-Lp}
        \lim_{t\to T_+}
        \int_{\R^d}|u(t)|^{p+1}\phi_{t,\veps}^2\ud x
        =
        \int_{\R^d}|u^*|^{p+1}\chi_\veps^2\ud x.
\end{align}
Combining \eqref{eqn:finite-main-cutoff-energy} and
\eqref{eqn:finite-main-cutoff-Lp}, we obtain
\begin{align}\label{eqn:finite-main-cutoff-Dirichlet}
        \lim_{t\to T_+}
        \bar E_{\phi_{t,\veps}}(u(t))
        =
        \bar E_{\chi_\veps}(u^*).
\end{align}

Set
$A_t:=B(x_0,r)\setminus B(x_0,\alp\rho(t))$. For $t$ sufficiently
close to $T_+$,
\begin{align}
0
&\leq
\bar E_{\phi_{t,\veps}}(u(t))
-
\bar E(u(t);A_t) \notag\\
&\leq
\bar E\left(
u(t);
B(x_0,\alp\rho(t))\setminus B(x_0,\alp\rho(t)/2)
\right)
+
\bar E\left(
u(t);
B(x_0,r+\veps)\setminus B(x_0,r)
\right).
\end{align}
Using \eqref{eqn:self-similar-positive-annulus-vanish},
\eqref{eqn:finite-main-cutoff-Dirichlet}, and the strong convergence
away from $\calS$, we conclude that
\begin{align}
&\limsup_{t\to T_+}
\left|
\bar E(u(t);A_t)-\bar E(u^*;B(x_0,r))
\right|\leq
2\bar E\left(
u^*;
B(x_0,r+\veps)\setminus B(x_0,r)
\right).
\end{align}
Letting $\veps\to0$ proves \eqref{eqn:no-conc-finite}. Finally, a diagonal argument applied to \eqref{eqn:self-similar-positive-annulus-vanish} gives a function $m(t)\to\infty$ as $t\to T_+$ such that $m(t)\rho(t)\to0$ as $t\to T_+$ and the limit in \eqref{eqn:self-similar-positive-annulus-vanish} remains valid with $a=m(t)^{-1}$ and $b=m(t)$. Taking $\xi(t):=\rho(t)/m(t)$ and $\nu(t):=m(t)\rho(t)$ gives \eqref{eqn:nu-xi-curves}.
\end{proof}
% \begin{proof}
% First, note that for any measurable subset $\Omega\subset\R^d$,
% \begin{align}\label{eqn:Ebar-E-Lp-identity}
% \bar E(u;\Omega)
% =
% 2E(u;\Omega)+\frac{2}{p+1}\int_\Omega |u|^{p+1}\ud x.
% \end{align}
% Applying Lemma~\ref{lem:ss-bu-nonlinear} and Lemma~\ref{lem:ss-bu-lp} to $\Omega_{\alp,r}(t):=B(x_0,r)\setminus B(x_0,\alp\sqrt{T_+-t})$, and then using \eqref{eqn:Ebar-E-Lp-identity}, gives \eqref{eqn:no-conc-finite}. This, along with Lemma~\ref{lem:ss-bu-lp}, implies \eqref{eqn:nu-xi-curves}.
% \end{proof}
We use the same strategy as above to prove the absence of self-similar energy concentration in the global-in-time case.
\begin{lem}[No self-similar nonlinear energy concentration in the global case]
\label{lem:vanishing-nonlinear-energy}
Let $u(t)$ be the solution of~\eqref{eqn:NLH} with initial data $u_0\in\dot H^1$, $T_+=\infty$, and $\sup_{t\in[0,\infty)}\bar E(u(t))<\infty$. Let $\phi_T\in C^\infty(\R^d)\cap W^{1,\infty}(\R^d)$, $T\geq1$, be such that
\begin{align}\label{eqn:global-cutoff-family-assumptions}
0\leq\phi_T\leq1,
\qquad
\lim_{T\to \infty}\phi_T(x)=0 \textrm{ for a.e. }x\in\R^d,
\qquad
\sup_{T\geq1}T\|\nabla\phi_T\|_{L^\infty}^2<\infty.
\end{align}
Then
\begin{align}\label{eqn:global-nonlinear-energy-evacuation}
\lim_{T\to \infty}E_{\phi_T}(u(T)) = 0.
\end{align}
The same conclusion holds along any sequence $T_n\to\infty$ if $\phi_n\in C^\infty(\R^d)\cap W^{1,\infty}(\R^d)$ satisfies $0\leq\phi_n\leq1$, $\lim_{n\to \infty}\phi_n(x) = 0$ for a.e. $x\in\R^d$, and $\sup_nT_n\|\nabla\phi_n\|_{L^\infty}^2<\infty$.
\end{lem}

\begin{proof}
Given $\veps>0$, choose $T_0>0$ such that
\begin{align}\label{eqn:global-tail-small}
\left(\int_{T_0}^{\infty}\|\partial_tu(t)\|_{L^2}^2\ud t\right)^{1/2}\leq\veps.
\end{align}
For $T>T_0$, Remark~\ref{rem:localized-energy-noncompact-cutoff} and \eqref{eqn:energy ineq I} give
\begin{align}
\begin{split}
\left|E_{\phi_T}(u(T))-E_{\phi_T}(u(T_0))\right|
&\leq
\int_{T_0}^{T}\|\partial_tu(t)\|_{L^2}^2\ud t\\
&\quad+
2\left(\int_{T_0}^{T}\|\partial_tu(t)\|_{L^2}^2\ud t\right)^{1/2}
\left(\int_{T_0}^{T}\int_{\R^d}|\nabla u|^2|\nabla\phi_T|^2\ud x\ud t\right)^{1/2}.
\end{split}
\end{align}
Since $\sup_{t\geq0}\bar E(u(t))<\infty$, $|\nabla\phi_T|^2\lesssim T^{-1}$, and \eqref{eq:tension-L2} holds, we obtain
\begin{align}\label{eqn:global-energy-difference-small}
\limsup_{T\to\infty}
|E_{\phi_T}(u(T))-E_{\phi_T}(u(T_0))|
\lesssim
\veps^2+\veps,
\end{align}
where the implicit constant depends on $\sup_{t\geq0}\bar E(u(t))$ and $\sup_{T\geq1}T\|\nabla\phi_T\|_{L^\infty}^2$. For fixed $T_0$, dominated convergence gives $E_{\phi_T}(u(T_0))\to0$ as $T\to\infty$. Taking $T\to\infty$ in \eqref{eqn:global-energy-difference-small}, then letting $\veps\to0$, proves \eqref{eqn:global-nonlinear-energy-evacuation}. The sequential version is proved by the same argument with $T=T_n$ and $\phi_T=\phi_n$.
\end{proof}

\begin{lem}[No self-similar $L^{p+1}$ concentration in the global case]
\label{lem:ss-global-weaker}
Let $u(t)$ be the solution of~\eqref{eqn:NLH} with initial data $u_0\in\dot H^1$, $T_+=\infty$, and $\sup_{t\geq0}\bar E(u(t))<\infty$. Then, for every $y\in\R^d$ and every $\alp>0$,
\begin{align}\label{eqn:two-star-limit}
\lim_{t\to\infty}\int_{|x-y|\geq\alp\sqrt t}|u(t,x)|^{p+1}\ud x=0.
\end{align}
\end{lem}

\begin{proof}
By translation invariance, assume $y=0$. Suppose \eqref{eqn:two-star-limit} fails. Then there exist $\bar\alp>0$, $\delta>0$, and $t_n\to\infty$ such that
\begin{align}\label{eqn:global-Lp-lower}
\int_{|x|\geq\bar\alp\sqrt{t_n}}|u(t_n,x)|^{p+1}\ud x\geq\delta.
\end{align}
Set $u_n:=u(t_n)$. Applying Proposition~\ref{prop:profile-decomp} to $(u_n)$ in $\dot H^1(\R^d)$ gives, for every $J\in\N$, 
\begin{align}\label{eqn:global-profile-decomp}
u_n
=
\sum_{j=1}^{J}W_n^j+r_n^J,
\qquad
W_n^j(x):=(\lambda_{j,n})^{-\frac{d-2}{2}}
W^j\left(\frac{x-a_{j,n}}{\lambda_{j,n}}\right),
\end{align}
with pairwise orthogonal parameters and
\begin{align}\label{eqn:global-profile-rem}
\lim_{J\to\infty}\limsup_{n\to\infty}\|r_n^J\|_{L^{p+1}}=0.
\end{align}

We claim that every nontrivial profile $W^j$ is stationary. Fix $j$ and define, for $0\leq s\leq1$,
\begin{align}\label{eqn:global-profile-flow}
U_n^j(z,s):=
(\lambda_{j,n})^{\frac{d-2}{2}}
u(t_n+\lambda_{j,n}^2s,a_{j,n}+\lambda_{j,n}z).
\end{align}
By the parabolic scaling invariance of \eqref{eqn:NLH}, $U_n^j$ solves
$\partial_sU_n^j=\Delta_zU_n^j+|U_n^j|^{p-1}U_n^j$ on $\R^d\times[0,1]$. Furthermore,
\begin{align}\label{eqn:global-rescaled-tension-small}
\int_0^1\|\partial_sU_n^j(s)\|_{L^2_z}^2\ud s
=
\int_{t_n}^{t_n+\lambda_{j,n}^2}\|\partial_tu(t)\|_{L^2_x}^2\ud t
\leq
\int_{t_n}^{\infty}\|\partial_tu(t)\|_{L^2_x}^2\ud t
\to0
\qquad\textrm{as }n\to\infty.\quad 
\end{align}
By \eqref{eq:profile-remainder-frame-vanishing} and parameter orthogonality, $U_n^j(0)\rightharpoonup W^j$ weakly in $\dot H^1(\R^d)$. Hence, after passing to a subsequence, $U_n^j(0)\to W^j$ strongly in $L^2_{\loc}(\R^d)$ by Rellich's theorem. Thus, for every compact $A\Subset\R^d$,
\begin{align}
\sup_{0\leq s\leq1}\|U_n^j(s)-W^j\|_{L^2(A)}
\leq
\|U_n^j(0)-W^j\|_{L^2(A)}
+
\|\partial_sU_n^j\|_{L^2_sL^2_z}
\to0
\qquad\textrm{as }n\to\infty .
\end{align}
Since $U^j_n\to W^j$ strongly in $L^2_{\loc}$ and $\{U^j_n\}_{n\in \N}$ is uniformly bounded in $L^\infty_s\dot H^1_z$, Sobolev gives a uniform $L^\infty_s L^{p+1}_z$ bound. Hence $U^j_n\to W^j$ strongly in $L^{q}_{\loc}$ for every $2\leq q<p+1$. On compact sets, the same conclusion for $1\leq q<2$ follows from Hölder's inequality. In particular, $U^j_n\to W^j$ strongly in $L^p_{\loc}$. Moreover, after passing to a subsequence,
\begin{align}\label{eq:global-profile-weak-gradient}
        \nabla U_n^j\rightharpoonup \nabla W^j
        \qquad\textrm{weakly in }L^2_{\loc}((0,1)\times\R^d).
\end{align}
The local strong convergence and the uniform $L^\infty_sL^{p+1}_z$ bound imply
\begin{align}
        |U_n^j|^{p-1}U_n^j
        \to
        |W^j|^{p-1}W^j
        \qquad\textrm{strongly in }L^1_{\loc}((0,1)\times\R^d).
\end{align}
Taking the weak formulation of the equation for $U_n^j$ against $\theta(s)\zeta(z)$, where $\theta\in C_c^\infty((0,1))$, $\int_0^1\theta\,\ud s=1$, and $\zeta\in C_c^\infty(\R^d)$, and then using \eqref{eqn:global-rescaled-tension-small}, we obtain
\begin{align}\label{eqn:global-profile-stationary}
        \int_{\R^d}\nabla W^j\cdot\nabla\zeta\,\ud z
        =
        \int_{\R^d}|W^j|^{p-1}W^j\zeta\,\ud z
        \qquad
        \textrm{for all }\zeta\in C_c^\infty(\R^d).
\end{align}
Hence $W^j$ solves \eqref{eq:yamabe}. Since every nontrivial stationary solution has $\bar E(W^j)\geq\bar E_*$, and since the Pythagorean expansion in the profile decomposition gives
\begin{align}
        \sum_{j=1}^{\infty}\|\nabla W^j\|_{L^2}^2
        \leq
        \limsup_{n\to\infty}\|u_n\|_{\dot H^1}^2<\infty,
\end{align}
there are only finitely many nontrivial profiles. After reordering, relabel them $W^1,\ldots,W^K$. Then $W^j\equiv0$ for every $j>K$. Set $S_n:=\sum_{j=1}^{K}W_n^j$ and $r_n:=u_n-S_n$. For every fixed $J\geq K$, we have $r_n=r_n^J$ because the profiles $W^j$ with $j>K$ are identically zero. Therefore, by \eqref{eqn:global-profile-rem},
\begin{align}\label{eqn:global-remainder-Lp}
        \limsup_{n\to\infty}\|r_n\|_{L^{p+1}}
        =
        \lim_{J\to \infty}\limsup_{n\to\infty}\|r_n^J\|_{L^{p+1}}
         =0.
\end{align}
By \eqref{eqn:global-Lp-lower}, we must have $K\geq1$. Fix $R>e$. Choose a radial logarithmic exterior cutoff $\phi_{n,R}\in C^\infty(\R^d)$ satisfying $0\leq\phi_{n,R}\leq1$, $\phi_{n,R}=0$ for $|x|\leq\frac{\bar\alp\sqrt{t_n}}{R}$, $\phi_{n,R}=1$ for $|x|\geq\bar\alp\sqrt{t_n}$, and
\begin{align}\label{eqn:global-log-gradient}
|\nabla\phi_{n,R}(x)|
\lesssim
\frac1{|x|\log R}
\qquad\textrm{on }
\frac{\bar\alp\sqrt{t_n}}{R}\leq |x|\leq\bar\alp\sqrt{t_n}.
\end{align}
Then
\begin{align}\label{eqn:global-log-Ld}
\|\nabla\phi_{n,R}\|_{L^d}\lesssim(\log R)^{-\frac{d-1}{d}},
\qquad
t_n\|\nabla\phi_{n,R}\|_{L^\infty}^2\lesssim_{\bar\alp,R}1.
\end{align}
Since $\phi_{n,R}=1$ on $\{|x|\geq\bar\alp\sqrt{t_n}\}$, \eqref{eqn:global-Lp-lower} and \eqref{eqn:global-remainder-Lp} give, for all large $n$,
\begin{align}\label{eqn:global-Sn-mass-lower}
\int_{\R^d}|S_n|^{p+1}\phi_{n,R}^2\ud x\geq c_\delta
\end{align}
for some $c_\delta>0$. Hence
\begin{align}\label{eqn:global-profile-Lp-lower}
\sum_{j=1}^{K}\int_{\R^d}|W_n^j|^{p+1}\phi_{n,R}^2\ud x
\geq
m_0,
\qquad
m_0:=K^{-p}c_\delta>0.
\end{align}
By Lemma \ref{lem:weighted-profile-decoupling}, we have
\begin{align}\label{eqn:global-weighted-energy-lower1}
E_{\phi_{n,R}}(u_n)
\geq
\sum_{j=1}^{K}E_{\phi_{n,R}}(W_n^j)-o_n(1).
\end{align}
Testing \eqref{eq:yamabe} for $W_n^j$ with $W_n^j\phi_{n,R}^2$ gives
\begin{align}
E_{\phi_{n,R}}(W_n^j)
=
\frac1d\int_{\R^d}|W_n^j|^{p+1}\phi_{n,R}^2\ud x
-
\int_{\R^d}W_n^j\phi_{n,R}\nabla W_n^j\cdot\nabla\phi_{n,R}\ud x.
\end{align}
The error term satisfies
\begin{align}
\left|
\int_{\R^d}W_n^j\phi_{n,R}\nabla W_n^j\cdot\nabla\phi_{n,R}\ud x
\right|
\leq
\|W^j\|_{L^{p+1}}\|\nabla W^j\|_{L^2}\|\nabla\phi_{n,R}\|_{L^d}.
\end{align}
Using \eqref{eqn:global-log-Ld}, choose $R>e$ sufficiently large so that
\begin{align}\label{eqn:global-log-error-small}
\sum_{j=1}^{K}\|W^j\|_{L^{p+1}}\|\nabla W^j\|_{L^2}\|\nabla\phi_{n,R}\|_{L^d}
\leq
\frac{m_0}{2d}.
\end{align}
Then \eqref{eqn:global-profile-Lp-lower}--\eqref{eqn:global-log-error-small} imply
\begin{align}\label{eqn:global-positive-localized-energy}
\liminf_{n\to\infty}E_{\phi_{n,R}}(u(t_n))\geq\frac{m_0}{2d}>0.
\end{align}

On the other hand, for this fixed $R$, the family $\phi_{n,R}$ satisfies the sequential version of Lemma~\ref{lem:vanishing-nonlinear-energy}: for every fixed $x\in\R^d$, $\phi_{n,R}(x)\to0$ as $n\to\infty$, and \eqref{eqn:global-log-Ld} gives $t_n\|\nabla\phi_{n,R}\|_{L^\infty}^2\lesssim_{\bar\alp,R}1$. Therefore
\begin{align}\label{eqn:global-nonlinear-contr}
E_{\phi_{n,R}}(u(t_n))\to0
\qquad\textrm{as }n\to\infty,
\end{align}
contradicting \eqref{eqn:global-positive-localized-energy}. This proves \eqref{eqn:two-star-limit}.
\end{proof}

\begin{lem}[No self-similar Dirichlet energy concentration in the global case]
\label{lem:ss-global}
Let $u(t)$ be the solution of~\eqref{eqn:NLH} with initial data $u_0\in\dot H^1$, $T_+=\infty$, and $\sup_{t\geq0}\bar E(u(t))<\infty$. Then, for every $y\in\R^d$ and every $\alp>0$,
\begin{align}\label{eqn:ss-global-final}
\lim_{t\to\infty}
\int_{|x-y|\geq\alp\sqrt t}|\nabla u(t,x)|^2\ud x=0.
\end{align}
\end{lem}

\begin{proof}
Let $\chi\in C_c^\infty(\R^d)$ be the standard cutoff function, so that $\chi\equiv1$ on $B(0,1)$ and $\chi\equiv0$ outside $B(0,2)$. Define
\begin{align}\label{eqn:global-final-cutoff}
\phi_t(x):=
1-\chi\left(\frac{2(x-y)}{\alp\sqrt t}\right).
\end{align}
Then $\phi_t\equiv0$ on $B\left(y,\frac{\alp}{2}\sqrt t\right)$, $\phi_t\equiv1$ on $\R^d\setminus B(y,\alp\sqrt t)$, $t\|\nabla\phi_t\|_{L^\infty}^2\lesssim_\alp1$, and $\phi_t(x)\to0$ as $t\to\infty$ for every fixed $x\in\R^d$. Thus Lemma~\ref{lem:vanishing-nonlinear-energy} gives
\begin{align}\label{eqn:global-final-E-zero}
E_{\phi_t}(u(t))\to0,\textrm{ as }t\to\infty .
\end{align}
By Lemma~\ref{lem:ss-global-weaker},
\begin{align}\label{eqn:global-final-Lp-zero}
\int_{\R^d}|u(t)|^{p+1}\phi_t^2\ud x
\leq
\int_{|x-y|\geq\frac{\alp}{2}\sqrt t}|u(t)|^{p+1}\ud x
\to0\textrm{, as }t\to\infty .
\end{align}
Using
\begin{align}\label{eqn:global-final-Ebar-identity}
\bar E_{\phi_t}(u(t))
=
2E_{\phi_t}(u(t))
+
\frac{2}{p+1}\int_{\R^d}|u(t)|^{p+1}\phi_t^2\ud x,
\end{align}
we obtain $\bar E_{\phi_t}(u(t))\to0$ as $t\to\infty$. Since $\phi_t\equiv1$ on $\{|x-y|\geq\alp\sqrt t\}$, \eqref{eqn:ss-global-final} follows.
\end{proof}

\subsection{Sequential Compactness}
We recall the standard sequential compactness theorem for Palais--Smale sequences of the nonlinear energy (cf. \cite{struwe-global} and \cite{duShi}).
\begin{lem}
\label{lem:whole-space-PS}
Let $d\geq3$ and $p=\frac{d+2}{d-2}$. Let
$\{f_n\}_{n\in\N}\subset\dot H^1(\R^d)$ satisfy
\begin{align}\label{eq:whole-space-PS-assumption}
        \sup_{n\in\N}
        \|\nabla f_n\|_{L^2}
        <\infty \text{ and } \lim_{n\to \infty} \left\|
        \Delta f_n+|f_n|^{p-1}f_n
        \right\|_{\dot H^{-1}}
        =0.
\end{align}
Then, after passing to a subsequence, there exist an integer $J\geq0$ and nonzero stationary solutions
$W^1,\ldots,W^J\in\dot H^1(\R^d)$, and sequences $a_n^j\in\R^d,$ $\lambda_n^j\in(0,\infty),$ $1\leq j\leq J,$ with the following properties. For every $j\neq k$, 
\begin{align}\label{eq:whole-space-PS-orthogonality}
         \lim_{n\to \infty}  \left(\frac{\lambda_n^j}{\lambda_n^k}
        +
        \frac{\lambda_n^k}{\lambda_n^j}
        +
        \frac{|a_n^j-a_n^k|^2}
        {\lambda_n^j\lambda_n^k}\right)=\infty.
\end{align}
Moreover,
\begin{align}\label{eq:whole-space-PS-conclusion}
        \lim_{n\to \infty}\|
        f_n-
        \sum_{j=1}^J
        W^j[a_n^j,\lambda_n^j]
        \|_{\dot H^1}
        =0.
\end{align}
In particular,
\begin{align}\label{eq:whole-space-PS-dirichlet-decoupling}
        \|\nabla f_n\|_{L^2}^2
        =
        \sum_{j=1}^J
        \|\nabla W^j\|_{L^2}^2
        +
        o_n(1)
\end{align}
and
\begin{align}\label{eq:whole-space-PS-Lp-decoupling}
        \|f_n\|_{L^{p+1}}^{p+1}
        =
        \sum_{j=1}^J
        \|W^j\|_{L^{p+1}}^{p+1}
        +
        o_n(1).
\end{align}
\end{lem}

\begin{proof}
Apply Proposition~\ref{prop:profile-decomp} to $\{f_n\}_{n\in\N}$.
For each nonzero profile $\psi^j$, let
$(a_n^j,\lambda_n^j)$ denote its translation and scaling parameters.
The orthogonality condition in
\eqref{eq:whole-space-PS-orthogonality} is equivalent to the condition in
part $\mathrm{(a)}$ of Proposition~\ref{prop:profile-decomp}. We first show that every nonzero profile is stationary. Fix
$\varphi\in C_c^\infty(\R^d)$ and define
$ \varphi_n^j(x)
        :=
        (\lambda_n^j)^{-\frac{d-2}{2}}
        \varphi(
        {x-a_n^j}/{\lambda_n^j}
        ).$
Then $ \|\nabla\varphi_n^j\|_{L^2}
        =
        \|\nabla\varphi\|_{L^2},$ and hence \eqref{eq:whole-space-PS-assumption} gives
\begin{align}\label{eq:whole-space-PS-rescaled-test}
        -\int_{\R^d}
        \nabla f_n\cdot\nabla\varphi_n^j\,\ud x
        +
        \int_{\R^d}
        |f_n|^{p-1}f_n\varphi_n^j\,\ud x
        \to0,
\end{align}
as $n\to \infty.$ Set
\begin{align}
g_n^j(y):=
(\lambda_n^j)^{\frac{d-2}{2}}
f_n(a_n^j+\lambda_n^j y).
\end{align}
Then $g_n^j\rightharpoonup\psi^j$ weakly in
$\dot H^1(\R^d)$ and, by local compactness,
$g_n^j\to\psi^j$ strongly in $L^r_{\loc}(\R^d)$ for every
$1\leq r<p+1$. Therefore, changing variables in
\eqref{eq:whole-space-PS-rescaled-test} and passing to the limit gives
\begin{align}
        -\int_{\R^d}
        \nabla\psi^j\cdot\nabla\varphi\,\ud x
        +
        \int_{\R^d}
        |\psi^j|^{p-1}\psi^j\varphi\,\ud x
        =
        0.
\end{align}
Thus, $\Delta\psi^j+|\psi^j|^{p-1}\psi^j=0$ in the sense of distributions. Every nonzero stationary solution $W$ satisfies
\begin{align}
        \|\nabla W\|_{L^2}^2
        =
        \|W\|_{L^{p+1}}^{p+1}
        \geq
        \bar E_*.
\end{align}
The Pythagorean expansion in
Proposition~\ref{prop:profile-decomp} therefore shows that there are only
finitely many nonzero profiles. Denote them by
$W^1,\ldots,W^J$, and set $r_n
        :=
        f_n-
        \sum_{j=1}^J
        W^j[a_n^j,\lambda_n^j].$
Since all nonzero profiles have been retained, every profile with index
larger than $J$ is zero. Since $r_n^L=r_n$ for every $L\geq J$, the last conclusion in Proposition~\ref{prop:profile-decomp} gives
$\|r_n\|_{L^{p+1}}\to0$, as $n\to \infty.$ Moreover, the parameter orthogonality and part $\mathrm{(b)}$ of Proposition~\ref{prop:profile-decomp} give
\begin{align}\label{eq:whole-space-PS-gradient-remainder}
        \int_{\R^d}
        \nabla f_n\cdot\nabla r_n\,\ud x
        =
        \|\nabla r_n\|_{L^2}^2
        +
        o_n(1).
\end{align}
Using $r_n$ as a test function against the expression $\Delta f_n+|f_n|^{p-1}f_n$, we obtain
\begin{align}
        \|\nabla r_n\|_{L^2}^2
        \leq
        \left\|
        \Delta f_n+|f_n|^{p-1}f_n
        \right\|_{\dot H^{-1}}
        \|\nabla r_n\|_{L^2}+
        \|f_n\|_{L^{p+1}}^p
        \|r_n\|_{L^{p+1}}
        +
        o_n(1).
\end{align}
The sequence $\{f_n\}$ is bounded in $L^{p+1}(\R^d)$ by Sobolev's
inequality, while $\{r_n\}$ is bounded in $\dot H^1(\R^d)$. Hence
\eqref{eq:whole-space-PS-assumption} and $\|r_n\|_{L^{p+1}}\to0$ as $n\to \infty$ imply $\|\nabla r_n\|_{L^2}\to0$ as $n\to\infty$. This proves \eqref{eq:whole-space-PS-conclusion}. By part
$\mathrm{(b)}$ of Proposition~\ref{prop:profile-decomp},
\begin{align}
        \|\nabla f_n\|_{L^2}^2
        =
        \sum_{j=1}^J
        \|\nabla W^j\|_{L^2}^2
        +
        \|\nabla r_n\|_{L^2}^2
        +
        o_n(1)
\end{align}
and
\begin{align}
        \|f_n\|_{L^{p+1}}^{p+1}
        =
        \sum_{j=1}^J
        \|W^j\|_{L^{p+1}}^{p+1}
        +
        \|r_n\|_{L^{p+1}}^{p+1}
        +
        o_n(1).
\end{align}
Therefore, \eqref{eq:whole-space-PS-conclusion} and
$\|r_n\|_{L^{p+1}}\to 0$ as $n\to \infty$ give
\eqref{eq:whole-space-PS-dirichlet-decoupling} and
\eqref{eq:whole-space-PS-Lp-decoupling}.
\end{proof}

\section{Analysis of the Remainder Term in the Profile Decomposition}
\label{sec:terminal-compactness}
In this section, we apply the profile decomposition to a bounded energy solution of the nonlinear heat equation along an arbitrary sequence of times. We prove that every nonzero profile is a stationary solution of \eqref{eq:yamabe}; for global solutions, this argument was carried out by Ishiwata in \cite{ishiwata2018potential}. We then use Lemmas~\ref{lem:ss-bu} and \ref{lem:ss-global} to localize the profiles to a suitable parabolic region and prove that the remainder, which is initially known to converge to zero in $L^{p+1}(\R^d)$, also converges to zero in $\dot H^1(\R^d)$.

\begin{lem}
\label{lem:arbitrary-stationary-profiles}
Let $u(t)$ be the solution of \eqref{eqn:NLH}, and assume that $\sup_{0\leq t<T_+}\|\nabla u(t)\|_{L^2}<\infty$. Let $t_n\to T_+$. If $T_+<\infty$, set $v_n:=u(t_n)-u^*$, where $u^*$
and $\calS$ are as in Theorem~\ref{thm:lwp}; if $T_+=\infty$, set $v_n:=u(t_n)$.

Then, after passing to a subsequence, there exist an integer $J\geq0$ and nonzero solutions $W^1,\ldots,W^J\in\dot H^1(\R^d)$ of \eqref{eq:yamabe}. There also exist sequences $a_n^j\in\R^d$ and $\lambda_n^j\in(0,\infty)$, $1\leq j\leq J$, and a sequence $r_n\in\dot H^1(\R^d)$ with the following properties. For every $j,k\in\{1,\ldots,J\}$ with $j\neq k$,
\begin{align}\label{eq:arbitrary-profile-orthogonality}
        \lim_{n\to \infty}\left(\frac{\lambda_n^j}{\lambda_n^k}
        +
        \frac{\lambda_n^k}{\lambda_n^j}
        +
        \frac{|a_n^j-a_n^k|}{\lambda_n^j}\right) =\infty.
\end{align}
Moreover, $v_n=\sum_{j=1}^J W^j[a_n^j,\lambda_n^j]+r_n$, and
$\|r_n\|_{L^{p+1}}\to0$ as $n\to\infty$. Furthermore,
\begin{align}\label{eq:arbitrary-stationary-pythagorean}
        \|\nabla v_n\|_{L^2}^2
        &=
        \sum_{j=1}^J\|\nabla W^j\|_{L^2}^2
        +
        \|\nabla r_n\|_{L^2}^2
        +
        o_n(1),\\
        \|v_n\|_{L^{p+1}}^{p+1}
        &=
        \sum_{j=1}^J\|W^j\|_{L^{p+1}}^{p+1}
        +
        \|r_n\|_{L^{p+1}}^{p+1}
        +
        o_n(1).
\end{align}
For every $1\leq j\leq J$,
$(\lambda_n^j)^{\frac{d-2}{2}}
r_n(a_n^j+\lambda_n^j\cdot)\rightharpoonup0$ weakly in
$\dot H^1(\R^d)$ as $n\to\infty$.

If $T_+<\infty$, then for every $1\leq j\leq J$ there exists
$x^{i_j}\in\calS$ such that
\begin{align}\label{eq:arbitrary-profile-finite-localization}
        a_n^j\to x^{i_j},
        \qquad
        \frac{|a_n^j-x^{i_j}|+\lambda_n^j}{\sqrt{T_+-t_n}}
        \to0
        \qquad
        \textrm{as }n\to\infty.
\end{align}
If $T_+=\infty$, then for every $1\leq j\leq J$,
\begin{align}\label{eq:arbitrary-profile-global-localization}
        \frac{|a_n^j|+\lambda_n^j}{\sqrt{t_n}}
        \to0
        \qquad
        \textrm{as }n\to\infty.
\end{align}
When $J=0$, the sums above are understood to be equal to zero.
\end{lem}
\begin{proof}
Apply Proposition~\ref{prop:profile-decomp} to the finite-energy sequence $\{v_n\}_{n\in\N}$. We denote the resulting profiles by $\psi^j$, their
parameters by $(a_n^j,\lambda_n^j)$, and the remainders by $r_n$. We first prove \eqref{eq:arbitrary-profile-finite-localization} and \eqref{eq:arbitrary-profile-global-localization}.

Assume that $T_+<\infty$. By
Lemma~\ref{lem:finite-time-no-escape}, the singular set is nonempty. Write
$\calS=\{x^1,\ldots,x^K\}$ and set
$\rho_n:=\sqrt{T_+-t_n}$. We claim that, for every $\alpha>0$,
\begin{align}\label{eq:arbitrary-profile-finite-concentration}
        \lim_{n\to \infty}\bar E(
        v_n;
        \R^d\setminus
        \cup_{i=1}^K B(x^i,\alpha\rho_n)
        )
        =0.
\end{align}
Choose $r>0$ sufficiently small such that the balls
$B(x^i,2r)$, $1\leq i\leq K$, are pairwise disjoint. Fix $R>0$ sufficiently large that $ \bigcup_{i=1}^K B(x^i,2r)\subset B(0,R).$ We decompose
\begin{align}
        \R^d\setminus\cup_{i=1}^K B(x^i,r)=(B(0,R)\setminus\cup_{i=1}^K B(x^i,r)) \cup(\R^d\setminus B(0,R)).
\end{align}
The strong convergence in Theorem~\ref{thm:lwp} controls the first
region. Lemma~\ref{lem:finite-time-no-escape} and
$\bar E(u^*;\R^d\setminus B(0,R))\to0$ control the second region.
Letting first $n\to\infty$ and then $R\to\infty$ gives
\begin{align}\label{eqn:exterior-estimate}
        \lim_{n\to \infty}\bar E(
        v_n;
        \R^d\setminus
        \cup_{i=1}^K B(x^i,r)
        )
        =0.
\end{align}
Fix $1\leq i\leq K$, and denote the annular region
$A_{i,n}:=B(x^i,r)\setminus B(x^i,\alpha\rho_n)$. We claim that
\begin{align}
        \lim_{n\to\infty}\bar E(u(t_n);A_{i,n})
        &=
        \bar E(u^*;B(x^i,r)),
        \quad
        \lim_{n\to\infty}\bar E(u^*;A_{i,n})
        =
        \bar E(u^*;B(x^i,r)),\\
        \lim_{n\to\infty}\int_{A_{i,n}}
        \nabla u(t_n)\cdot\nabla u^*\,\ud x
        &=
        \bar E(u^*;B(x^i,r)),
\end{align}
as $n\to \infty$. The first two limits follow from Lemma~\ref{lem:ss-bu} and the absolute continuity of the $\dot H^1$ and $L^{p+1}$ norms of $u^*$. For the third, fix $0<\delta<r$. Strong convergence away from $\calS$ gives the limit on $B(x^i,r)\setminus B(x^i,\delta)$, while Cauchy--Schwarz bounds the integral over $A_{i,n}\cap B(x^i,\delta)$ by
$\sup_n\|\nabla u(t_n)\|_{L^2}\|\nabla u^*\|_{L^2(B(x^i,\delta))}$.
Letting $n\to\infty$ and then $\delta\to0$, and using the absolute continuity of the $\dot{H}^1$-norm of $u^*$, proves the third limit.

Thus $\lim_{n\to\infty}\bar E(v_n;A_{i,n})=0$. Combining this with \eqref{eqn:exterior-estimate} proves \eqref{eq:arbitrary-profile-finite-concentration}. Fix now a nonzero profile $\psi^j$, and set
$V_n^j(y):=(\lambda_n^j)^{\frac{d-2}{2}}
v_n(a_n^j+\lambda_n^j y)$. By
Proposition~\ref{prop:profile-decomp},
$V_n^j\rightharpoonup\psi^j$ weakly in $\dot H^1(\R^d)$. Set
$R_n^j:=\rho_n/\lambda_n^j$ and
$z_{i,n}^j:=(x^i-a_n^j)/\lambda_n^j$ for $1\leq i\leq K$. By scale
invariance, \eqref{eq:arbitrary-profile-finite-concentration} becomes
\begin{align}\label{eq:arbitrary-profile-frame-concentration}
        \lim_{n\to \infty}\bar E(
        V_n^j;
        \R^d\setminus
        \cup_{i=1}^K B(z_{i,n}^j,\alpha R_n^j)
        )
        =0
\end{align}
for every $\alpha>0$. We first show that $R_n^j\to\infty$ as $n\to \infty$. Suppose not. Then, after passing to
a subsequence, $R_n^j$ is bounded and, for every $1\leq i\leq K$, either
$z_{i,n}^j$ converges in $\R^d$ or $|z_{i,n}^j|\to\infty$. Since
$\psi^j$ is nonzero, we can choose $\chi\in C_c^\infty(\R^d)$ such that
the support of $\chi$ is disjoint from all the finite limits of
$z_{i,n}^j$ and
$\int_{\R^d}|\nabla\psi^j|^2\chi^2\,\ud y>0$. For $\alpha>0$
sufficiently small, $\supp\chi$ is contained in the complement of
$\cup_{i=1}^K B(z_{i,n}^j,\alpha R_n^j)$ for all sufficiently large
$n$. Hence \eqref{eq:arbitrary-profile-frame-concentration} and lower semicontinuity give
\begin{align}
        0
        <
        \int_{\R^d}|\nabla\psi^j|^2\chi^2\,\ud y
        \leq
        \liminf_{n\to\infty}
        \int_{\R^d}|\nabla V_n^j|^2\chi^2\,\ud y
        =
        0,
\end{align}
which is a contradiction. Thus $R_n^j\to\infty$.

We next claim that $\min_{1\leq i\leq K}|z_{i,n}^j|/R_n^j\to0$. If not then, after passing to
a subsequence, there exists $c>0$ such that $|z_{i,n}^j|\geq cR_n^j$ for every $1\leq i\leq K$. Choose
$\chi\in C_c^\infty(\R^d)$ such that $\int_{\R^d}|\nabla\psi^j|^2\chi^2\,\ud y>0$, and fix
$0<\alpha<c/2$. Since $R_n^j\to\infty$, the support of $\chi$ is
contained in the complement of $\cup_{i=1}^K B(z_{i,n}^j,\alpha R_n^j)$ for all sufficiently large
$n$. The same argument as in the previous paragraph gives a contradiction.

After passing to a subsequence, there is therefore an index
$i_j\in\{1,\ldots,K\}$ such that
$|z_{i_j,n}^j|/R_n^j\to0$ as $n\to \infty.$ Consequently,
\begin{align}
       \lim_{n\to\infty} \frac{\lambda_n^j}{\rho_n}
        =
        \lim_{n\to\infty}\frac{1}{R_n^j}
        =0
        ,\qquad 
        \lim_{n\to\infty}\frac{|a_n^j-x^{i_j}|}{\rho_n}
        =
        \lim_{n\to\infty}\frac{|z_{i_j,n}^j|}{R_n^j}
        =0.
\end{align}
Since the set of profiles is countable, a diagonal argument allows us
to assume that \eqref{eq:arbitrary-profile-finite-localization} holds
simultaneously for every nonzero profile. 

Assume now that $T_+=\infty$, and set $\rho_n:=\sqrt{t_n}$. By
Lemma~\ref{lem:ss-global}, applied with $y=0$, for every $\alpha>0$,
\begin{align}
        \lim_{n\to \infty}\bar E(
        v_n;
        \R^d\setminus B(0,\alpha\rho_n)
        )
        =0.
\end{align}
For a nonzero profile $\psi^j$, set
$V_n^j(y):=(\lambda_n^j)^{\frac{d-2}{2}}
v_n(a_n^j+\lambda_n^j y)$,
$R_n^j:=\rho_n/\lambda_n^j$, and
$z_n^j:=-a_n^j/\lambda_n^j$. By scale invariance,
\begin{align}
        \lim_{n\to\infty}\bar E(
        V_n^j;
        \R^d\setminus B(z_n^j,\alpha R_n^j)
        )
        =0
\end{align}
for every $\alpha>0$. Arguing as in the case when $T_+<\infty$, this time with only one center, gives $R_n^j\to\infty$ and $|z_n^j|/R_n^j\to0$ as $n\to\infty$, or equivalently,
\begin{align}
\lim_{n\to\infty}\frac{\lambda_n^j}{\sqrt{t_n}}=0
        ,\qquad \lim_{n\to\infty}\frac{|a_n^j|}{\sqrt{t_n}}=0.
\end{align}
Thus, \eqref{eq:arbitrary-profile-global-localization} holds.

We next prove that every nonzero profile is stationary. Fix such a profile
$\psi^j$, and write $a_n:=a_n^j$ and
$\lambda_n:=\lambda_n^j$. If $T_+<\infty$, set $I=[-1,0]$, while if
$T_+=\infty$, set $I=[0,1]$. Define
\begin{align}\label{eq:arbitrary-profile-rescaled-flow}
        U_n(y,s)
        :=
        \lambda_n^{\frac{d-2}{2}}
        u(t_n+\lambda_n^2s,a_n+\lambda_n y),
        \qquad
        s\in I.
\end{align}
In the finite-time case,
$\lambda_n^2/(T_+-t_n)\to0$ by
\eqref{eq:arbitrary-profile-finite-localization}, so
\eqref{eq:arbitrary-profile-rescaled-flow} is defined on
$\R^d\times I$ for all sufficiently large $n$. In either case,
\eqref{eq:tension-L2} gives
\begin{align}\label{eq:arbitrary-profile-rescaled-tension}
        \int_I\|\partial_sU_n(s)\|_{L^2}^2\,\ud s
        =
        \begin{cases}
        \displaystyle
        \int_{t_n-\lambda_n^2}^{t_n}
        \|\partial_tu(t)\|_{L^2}^2\,\ud t,
        & T_+<\infty,\\[3mm]
        \displaystyle
        \int_{t_n}^{t_n+\lambda_n^2}
        \|\partial_tu(t)\|_{L^2}^2\,\ud t,
        & T_+=\infty
        \end{cases}
\end{align}
and so 
\begin{align}
\lim_{n\to \infty} \int_I\|\partial_sU_n(s)\|_{L^2}^2\,\ud s  = 0.
\end{align}
By the asymptotic orthogonality of parameters in
Proposition~\ref{prop:profile-decomp} and
\eqref{eq:profile-remainder-frame-vanishing}, we obtain
\begin{align}
        \lambda_n^{\frac{d-2}{2}}
        v_n(a_n+\lambda_n\cdot)
        \to
        \psi^j
        \qquad
        \textrm{strongly in }L^2_{\loc}(\R^d).
\end{align}
If $T_+<\infty$, then for every bounded domain $K\Subset\R^d$, as $n\to \infty$ we have
\begin{align}
        \|
        \lambda_n^{\frac{d-2}{2}}
        u^*(a_n+\lambda_n\cdot)
        \|_{L^{p+1}(K)}^{p+1}
        =
        \int_{a_n+\lambda_nK}|u^*(x)|^{p+1}\,\ud x
        \to0,
\end{align}
and
\begin{align}
        \|
        \nabla[
        \lambda_n^{\frac{d-2}{2}}
        u^*(a_n+\lambda_n\cdot)
        ]
        \|_{L^2(K)}^2
        =
        \int_{a_n+\lambda_nK}|\nabla u^*(x)|^2\,\ud x
        \to0,
\end{align}
as $n\to \infty.$ Indeed, $a_n\to x^{i_j}$ and $\lambda_n\to0$, and both limits follow from the absolute continuity of the $L^{p+1}$ and $\dot H^1$ norms of $u^*$. Since $p+1>2$, the first limit also gives convergence to zero in
$L^2(K)$. Consequently,
$U_n(0)\to\psi^j$ strongly in $L^2_{\loc}(\R^d)$.

By \eqref{eq:arbitrary-profile-rescaled-tension},
\begin{align}
        \sup_{s\in I}\|U_n(s)-U_n(0)\|_{L^2}
        \leq
        |I|^{1/2}\|\partial_sU_n\|_{L^2(I\times\R^d)}
        \to0,
\end{align}
as $n\to \infty.$ Hence $U_n\to\psi^j$ strongly in
$L^2_{\loc}(\R^d\times I)$ as $n\to \infty.$ If $p>2$, interpolation with the uniform $L^{p+1}$ bound gives strong convergence in $L^p_{\loc}$; if $p\leq2$,
the same follows from H\"older's inequality. Moreover,
$\nabla U_n\rightharpoonup\nabla\psi^j$ weakly in
$L^2_{\loc}(\R^d\times I)$ as $n\to \infty.$ Passing to the limit in the weak
formulation of \eqref{eqn:NLH}, integrated against $\theta(s)\zeta(y)$ with $\theta\in C_c^\infty(\R)$, $\supp\theta\subset I$, and $\int_I\theta\,\ud s=1$, gives
\begin{align}
        \int_{\R^d}
        \nabla\psi^j\cdot\nabla\zeta\,\ud y
        =
        \int_{\R^d}
        |\psi^j|^{p-1}\psi^j\zeta\,\ud y
        \qquad
        \textrm{for every }\zeta\in C_c^\infty(\R^d).
\end{align}
Thus $\psi^j$ is a nonzero stationary solution of \eqref{eq:yamabe}. In particular,
$\|\nabla\psi^j\|_{L^2}^2 =\|\psi^j\|_{L^{p+1}}^{p+1}\geq\bar E_*$. By the Pythagorean expansion in
Proposition~\ref{prop:profile-decomp},
\begin{align}
        \sum_{j=1}^\infty
        \|\nabla\psi^j\|_{L^2}^2
        \leq
        \limsup_{n\to\infty}
        \|\nabla v_n\|_{L^2}^2
        <\infty.
\end{align}
Hence only finitely many profiles are nonzero. After reordering, denote
them by $W^1,\ldots,W^J$. If $J\geq1$, set $r_n:=r_n^J$, while if
$J=0$, set $r_n:=v_n$.

All profiles with index larger than $J$ are zero. Thus
$r_n^L=r_n$ for every sufficiently large $L$, and the last conclusion of
Proposition~\ref{prop:profile-decomp} gives
$\|r_n\|_{L^{p+1}}\to0$. The decomposition, the two Pythagorean
expansions, the parameter separation, and the weak convergence
of the rescaled remainder term now follow from parts $\mathrm{(a)}$ and
$\mathrm{(b)}$ of Proposition~\ref{prop:profile-decomp}. This completes
the proof.
\end{proof}
The preceding lemma gives the following decomposition:
\begin{lem}
\label{lem:no-ghost-finite-splitting}
Assume that $T_+<\infty$. Let $\tau_n\to T_+$ and set
$v_n:=u(\tau_n)-u^*$. After passing to a subsequence, we have
\begin{align}\label{eq:no-ghost-finite-Dirichlet-splitting}
        \|\nabla u(\tau_n)\|_{L^2}^2
        =
        \|\nabla u^*\|_{L^2}^2
        +
        \|\nabla v_n\|_{L^2}^2
        +
        o_n(1)
\end{align}
and
\begin{align}\label{eq:no-ghost-finite-BL-splitting}
        \|u(\tau_n)\|_{L^{p+1}}^{p+1}
        =
        \|u^*\|_{L^{p+1}}^{p+1}
        +
        \|v_n\|_{L^{p+1}}^{p+1}
        +
        o_n(1).
\end{align}
Consequently,
\begin{align}\label{eq:no-ghost-finite-energy-splitting}
        E(u(\tau_n))
        =
        E(u^*)
        +
        E(v_n)
        +
        o_n(1).
\end{align}
\end{lem}
\begin{proof}
Since $u(\tau_n)\rightharpoonup u^*$ weakly in $\dot H^1(\R^d)$,
we have $v_n\rightharpoonup0$ weakly in $\dot H^1(\R^d)$. Therefore,
\begin{align}
        \int_{\R^d}
        \nabla u^*\cdot\nabla v_n\,\ud x
        \to0,
\end{align}
which proves
\eqref{eq:no-ghost-finite-Dirichlet-splitting}. We next prove
\eqref{eq:no-ghost-finite-BL-splitting}. By the strong convergence away
from $\calS$ in Theorem~\ref{thm:lwp} and the finiteness of $\calS$,
after passing to a subsequence,
\begin{align}
        u(\tau_n,x)\to u^*(x)
        \qquad
        \textrm{for almost every }x\in\R^d.
\end{align}
Thus $v_n(x)\to0$ for almost every $x\in\R^d$. Set $q:=p+1$ and $D_n:=|u(\tau_n)|^q-|u^*|^q-|v_n|^q$.
For every $\veps>0$, there exists $C_\veps>0$ such that
\begin{align}
        \big|
        |a+b|^q-|a|^q-|b|^q
        \big|
        \leq
        \veps|a|^q+C_\veps|b|^q
\end{align}
for every $a,b\in\R$. Applying this inequality with
$a=v_n(x)$ and $b=u^*(x)$ gives
\begin{align}
        |D_n|
        \leq
        \veps|v_n|^q+C_\veps|u^*|^q.
\end{align}
Set $H_n
        :=
        \big(
        |D_n|-\veps|v_n|^q
        \big)^+.$ Then $H_n\to0$ almost everywhere and
$0\leq H_n\leq C_\veps|u^*|^q$. Hence, by the dominated convergence
theorem,
\begin{align}
        \lim_{n\to\infty}
        \int_{\R^d}H_n\,\ud x
        =
        0.
\end{align}
Since $\{v_n\}_{n\in\N}$ is bounded in $L^q(\R^d)$,
\begin{align}
        \limsup_{n\to\infty}
        \int_{\R^d}|D_n|\,\ud x
        \lesssim
        \veps.
\end{align}
Letting $\veps\to0$ proves
\eqref{eq:no-ghost-finite-BL-splitting}. The identity
\eqref{eq:no-ghost-finite-energy-splitting} follows from
\eqref{eqn:energy}.
\end{proof}
We now prove that the remainder in Lemma~\ref{lem:arbitrary-stationary-profiles} converges to zero in $\dot H^1(\R^d)$. 
\begin{lem}
\label{lem:arbitrary-time-no-ghost}
Let $t_n\to T_+$ and let $r_n$ be the remainder in
Lemma~\ref{lem:arbitrary-stationary-profiles}. Then
\begin{align}\label{eq:arbitrary-no-ghost}
        \|\nabla r_n\|_{L^2}\to0.
\end{align}
\end{lem}
The proof differs in the cases $d\geq5$ and $d\in\{3,4\}$. For $d\in\{3,4\}$, we use the {\L}ojasiewicz--Simon inequality for the nonlinear energy functional, which we recall below.
\begin{lem}[{\L}ojasiewicz--Simon inequality]
\label{lem:no-ghost-LS}
Assume that $d\in\{3,4\}$. Let $W\in\dot H^1(\R^d)$ be a stationary solution
of \eqref{eq:yamabe}, possibly zero. Then there exist constants $C_W>0$,
$\sigma_W>0$, and $\theta_W\in[1/2,1)$ such that
\begin{align}\label{eq:no-ghost-LS}
        |E(v)-E(W)|^{\theta_W}
        \leq
        C_W
        \|\Delta v+|v|^{p-1}v\|_{\dot H^{-1}}
\end{align}
whenever $v\in\dot H^1(\R^d)$ satisfies
$\|v-W\|_{\dot H^1}<\sigma_W$.
\end{lem}
\begin{proof}
We apply \cite[Theorem~1]{feehan2020lojasiewicz} to the Hilbert
space $\dot H^1(\R^d)$ and verify its hypotheses. For $f\in\dot H^1(\R^d)$, define $\calL(f)=-\Delta f\in
\dot H^{-1}(\R^d)$ by
\begin{align}
        \calL(f)[g]
        :=
        \int_{\R^d}\nabla f\cdot\nabla g\,\ud x,
        \text{ for any } g\in\dot H^1(\R^d).
\end{align}
Then
\begin{align}
        \|\calL(f)\|_{\dot H^{-1}}
        =
        \|\nabla f\|_{L^2},
        \qquad
        \calL(f)[f]
        =
        \|\nabla f\|_{L^2}^2>0
\end{align}
for every nonzero $f\in\dot H^1(\R^d)$. Thus $\calL$ is a continuous
definite embedding in the sense of \cite{feehan2020lojasiewicz}.

Moreover, $p=5$ when $d=3$ and $p=3$ when $d=4$. Thus the nonlinear
energy $E:\dot H^1(\R^d)\to\R$ is a continuous polynomial and is therefore
analytic, with derivative
\begin{align}
        E'(v)=-\Delta v-|v|^{p-1}v\in\dot H^{-1}(\R^d).
\end{align}

Since $W$ solves \eqref{eq:yamabe}, it is a critical point of $E$. It remains to verify the Fredholm property of the Hessian of $E$. Set $V_W:=p|W|^{p-1}$. Since $(p-1)d/2=p+1$ and
$W\in L^{p+1}(\R^d)$, we have $V_W\in L^{d/2}(\R^d)$. For every
$V\in L^{d/2}(\R^d)$ and $f\in\dot H^1(\R^d)$, H\"older's inequality and
Sobolev's inequality give
\begin{align}\label{eq:no-ghost-multiplication-bound}
        \|Vf\|_{\dot H^{-1}}
        \lesssim
        \|V\|_{L^{d/2}}
        \|\nabla f\|_{L^2}.
\end{align}
We next show that multiplication by $V_W$ is compact from
$\dot H^1(\R^d)$ to $\dot H^{-1}(\R^d)$. Choose
$V_m\in C_c^\infty(\R^d)$ such that
\begin{align}
        \lim_{m\to\infty}
        \|V_W-V_m\|_{L^{d/2}}
        =
        0.
\end{align}
Let $\{f_n\}_{n\in\N}$ be bounded in $\dot H^1(\R^d)$ and suppose that
$f_n\rightharpoonup0$ weakly in $\dot H^1(\R^d)$. By
\eqref{eq:no-ghost-multiplication-bound},
\begin{align}
        \limsup_{n\to\infty}
        \|(V_W-V_m)f_n\|_{\dot H^{-1}}
        \lesssim
        \|V_W-V_m\|_{L^{d/2}}
        \sup_{n\in\N}\|\nabla f_n\|_{L^2}.
\end{align}
For fixed $m$, Rellich's theorem and the weak convergence give
$\|f_n\|_{L^2(\supp V_m)}\to0$, as $n\to \infty$ and hence
\begin{align}
        \|V_mf_n\|_{\dot H^{-1}}
        \lesssim_m
        \|f_n\|_{L^2(\supp V_m)}
        \to0,
\end{align}
as $n\to \infty.$ Letting first $n\to\infty$ and then $m\to\infty$ proves the desired compactness. Since
$-\Delta:\dot H^1(\R^d)\to\dot H^{-1}(\R^d)$ is an isomorphism, the Hessian
\begin{align}
        E''(W)=-\Delta-p|W|^{p-1}
        :\dot H^1(\R^d)\to\dot H^{-1}(\R^d)
\end{align}
is Fredholm of index zero. Since $E'(v)=-(\Delta v+|v|^{p-1}v)$, the conclusion follows from \cite[Theorem~1]{feehan2020lojasiewicz}.
\end{proof}
Using the {\L}ojasiewicz--Simon inequality we obtain the following:
\begin{lem}[Vanishing of the remainder when $d\in\{3,4\}$]
\label{lem:arbitrary-time-no-ghost-low-dimensional}
Assume that $d\in\{3,4\}$. Let $t_n\to T_+$ and let $r_n$ be the
remainder in Lemma~\ref{lem:arbitrary-stationary-profiles}. Then
$\|\nabla r_n\|_{L^2}\to0$ as $n\to \infty.$
\end{lem}
\begin{proof}
Let
\begin{align}
        \mathscr A
        :=
        \{
        \|\nabla W\|_{L^2}^2:
        W\in\dot H^1(\R^d)
        \textrm{ is a nonzero solution of }\eqref{eq:yamabe}
        \}.
\end{align}
Let $\calC$ denote the set of stationary solutions of
\eqref{eq:yamabe}, including zero. If $W,V\in\calC$ and
$\|V-W\|_{\dot H^1}<\sigma_W$, then \eqref{eq:no-ghost-LS} gives
$E(V)=E(W)$. Thus $E$ is locally constant on $\calC$. Since $\dot H^1(\R^d)$ is a separable Hilbert space, $\calC$ is second countable, and the following open cover of $\calC$
\begin{align}
        \{
        V\in\dot H^1(\R^d):
        \|V-W\|_{\dot H^1}<\sigma_W
        \},
        \qquad
        W\in\calC,
\end{align}
has a countable subcover. Since the nonlinear energy $E$ is constant on the intersection of
each of these neighborhoods with $\calC$, the set
$\{E(W):W\in\calC\}$ is countable. Moreover,
\begin{align}
        E(W)
        =
        \frac1d\|\nabla W\|_{L^2}^2
\end{align}
for every stationary solution $W$. Therefore, $\mathscr A$ is countable. For a constant $C>0$, set
\begin{align}\label{eq:no-ghost-mass-set}
        \mathscr M_C
        :=
        \{
        A_1+\cdots+A_L:
        L\in\N\cup\{0\},
        \ A_j\in\mathscr A,\quad1\leq j\leq L,
        \ A_1+\cdots+A_L\leq C
        \}.
\end{align}
Here the sum is zero when $L=0$. Since $\mathscr A$ is countable, the set of its $L$-fold sums is countable for every $L\geq0$, and hence $\mathscr M_C$ is countable. We will apply this with $C$ determined by the energy bound of $u(t)$ to prove that $b(t):=\|u(t)\|_{L^{p+1}}^{p+1}$ has a limit as $t\to T_+$.

We next construct a sequence of times along which the remainder in the profile decomposition converges strongly to zero
in $\dot H^1$. Assume first that $T_+<\infty$. Write $\calS=\{x^1,\ldots,x^K\}$. Let $\delta_n>0$ satisfy
$\lim_{n\to\infty}\delta_n=0$. By \eqref{eq:tension-L2}, there exists
$s_n\in[T_+-2\delta_n,T_+-\delta_n]$ such that, with
$\rho_n:=\sqrt{T_+-s_n}$,
\begin{align}\label{eq:no-ghost-finite-good-time}
\rho_n^2\|\partial_tu(s_n)\|_{L^2}^2
\leq
2\int_{T_+-2\delta_n}^{T_+-\delta_n}
\|\partial_tu(t)\|_{L^2}^2\,\ud t
\to0,
\end{align}
as $n\to \infty.$ Choose $\zeta_n\in C_c^\infty(\R^d)$ such that $\zeta_n=1$ on
$\cup_{i=1}^K B(x^i,\rho_n)$, $\supp\zeta_n \subset \cup_{i=1}^K B(x^i,2\rho_n),$ $0\leq \zeta_n\leq 1$, $|\nabla\zeta_n|\lesssim\rho_n^{-1}$, and $|\Delta\zeta_n|\lesssim\rho_n^{-2}$. For all sufficiently large $n$,
the balls in the support of $\zeta_n$ are pairwise disjoint.
We claim that
\begin{align}\label{eq:no-ghost-finite-good-decomp}
        \lim_{n\to\infty}
        \left[
        \|u(s_n)-u^*-\zeta_nu(s_n)\|_{\dot H^1}
        +
        \|u(s_n)-u^*-\zeta_nu(s_n)\|_{L^{p+1}}
        \right]
        =
        0.
\end{align}
Set $w_n:=\zeta_nu(s_n)$ and $A_n:=\bigcup_{i=1}^K (B(x^i,2\rho_n)\setminus B(x^i,\rho_n))$. Applying Lemmas~\ref{lem:ss-bu} and \ref{lem:ss-bu-lp} with
$\alp=1$ and $\alp=2$, and subtracting the corresponding limits, gives
\begin{align}\label{eq:no-ghost-finite-transition-small}
        \lim_{n\to\infty}
        [
        \bar E(u(s_n);A_n)
        +
        \|u(s_n)\|_{L^{p+1}(A_n)}^{p+1}
        ]
        =
        0.
\end{align}
This is because, for every $1\leq i\leq K$ and every fixed sufficiently small
$r>0$, Lemma~\ref{lem:ss-bu} gives
\begin{align}
        \lim_{n\to\infty}
        \bar E(
        u(s_n);
        B(x^i,r)\setminus B(x^i,\rho_n)
        )
        &=
        \bar E(u^*;B(x^i,r)),
        \notag\\
        \lim_{n\to\infty}
        \bar E(
        u(s_n);
        B(x^i,r)\setminus B(x^i,2\rho_n)
        )
        &=
        \bar E(u^*;B(x^i,r)).
\end{align}
Subtracting these two limits gives
\begin{align}
        \lim_{n\to\infty}
        \bar E(
        u(s_n);
        B(x^i,2\rho_n)\setminus B(x^i,\rho_n)
        )
        =
        0.
\end{align}
The same argument using Lemma~\ref{lem:ss-bu-lp} gives
\begin{align}
        \lim_{n\to\infty}
        \|u(s_n)\|_{L^{p+1}
        (B(x^i,2\rho_n)\setminus B(x^i,\rho_n))}
        =
        0.
\end{align}
Since $\calS$ is finite, this proves
\eqref{eq:no-ghost-finite-transition-small}. Lemmas~\ref{lem:ss-bu} and \ref{lem:ss-bu-lp}, the strong convergence away from $\calS$, and the absolute continuity of the $\dot{H}^1$ and $L^{p+1}$ norms of $u^*$ give
\begin{align}
        \|(1-\zeta_n)\nabla u(s_n)-\nabla u^*\|_{L^2}
        +
        \|(1-\zeta_n)u(s_n)-u^*\|_{L^{p+1}}
        \to0,
\end{align}
as $n\to \infty.$ Moreover, H\"older's inequality and \eqref{eq:no-ghost-finite-transition-small} give
\begin{align}
        \|u(s_n)\nabla\zeta_n\|_{L^2}
        \lesssim
        \|u(s_n)\|_{L^{p+1}(A_n)}
        \|\nabla\zeta_n\|_{L^d}
        \to0,
\end{align}
as $n\to \infty.$ Since
\begin{align}
        \nabla(u(s_n)-u^*-\zeta_nu(s_n))
        =
        (1-\zeta_n)\nabla u(s_n)-\nabla u^*
        -u(s_n)\nabla\zeta_n,
\end{align}
this proves \eqref{eq:no-ghost-finite-good-decomp}. In particular,
$\{w_n\}_{n\in\N}$ is bounded in $\dot H^1(\R^d)$. Using \eqref{eqn:NLH}, we compute
\begin{align}\label{eq:no-ghost-cutoff-residual}
        \Delta w_n+|w_n|^{p-1}w_n
        &=
        \zeta_n\partial_tu(s_n)
        +
        2\nabla\zeta_n\cdot\nabla u(s_n)
        +
        u(s_n)\Delta\zeta_n
        \notag\\
        &\quad+
        |\zeta_nu(s_n)|^{p-1}\zeta_nu(s_n)
        -
        \zeta_n|u(s_n)|^{p-1}u(s_n).
\end{align}
Let $\varphi\in C_c^\infty(\R^d)$ satisfy
$\|\nabla\varphi\|_{L^2}\leq1$. Since
$\supp\zeta_n\subset\cup_{i=1}^K B(x^i,2\rho_n)$ and the derivatives of
$\zeta_n$ are supported in $A_n$, we have
\begin{align}
        \|\zeta_n\varphi\|_{L^2}
        \lesssim
        \rho_n\|\nabla\varphi\|_{L^2},
        \quad 
        \|\nabla\zeta_n\|_{L^d}
        +
        \|\Delta\zeta_n\|_{L^{d/2}}
        \lesssim
        1.
\end{align}
Moreover, since $0\leq\zeta_n\leq1$,
\begin{align}
        \big|
        |\zeta_nu(s_n)|^{p-1}\zeta_nu(s_n)
        -
        \zeta_n|u(s_n)|^{p-1}u(s_n)
        \big|
        \lesssim
        |u(s_n)|^p\mathbf{1}_{A_n}.
\end{align}
Thus H\"older's inequality and Sobolev's inequality give
\begin{align}\label{eq:no-ghost-finite-residual-estimate}
        \|\Delta w_n+|w_n|^{p-1}w_n\|_{\dot H^{-1}}
        &\lesssim
        \rho_n\|\partial_tu(s_n)\|_{L^2}
        +
        \|\nabla u(s_n)\|_{L^2(A_n)}
        \notag\\
        &\quad+
        \|u(s_n)\|_{L^{p+1}(A_n)}
        +
        \|u(s_n)\|_{L^{p+1}(A_n)}^p.
\end{align}
By \eqref{eq:no-ghost-finite-good-time} and
\eqref{eq:no-ghost-finite-transition-small},
\begin{align}\label{eq:no-ghost-finite-PS}
        \lim_{n\to\infty}
        \|\Delta w_n+|w_n|^{p-1}w_n\|_{\dot H^{-1}}
        =
        0.
\end{align}
After passing to a subsequence, Lemma~\ref{lem:whole-space-PS} and
\eqref{eq:no-ghost-finite-good-decomp} give a finite sum of stationary
solutions whose difference from $u(s_n)-u^*$ converges to zero in
$\dot H^1(\R^d)$.

Assume now that $T_+=\infty$. Let $\tau_n\to\infty$. By
\eqref{eq:tension-L2}, we may choose $s_n\in[\tau_n,2\tau_n]$ such that
\begin{align}\label{eq:no-ghost-global-good-time}
        s_n\|\partial_tu(s_n)\|_{L^2}^2
        \leq
        2\int_{\tau_n}^{2\tau_n}
        \|\partial_tu(t)\|_{L^2}^2\,\ud t
        \to0,
\end{align}
as $n\to \infty.$ For every fixed $\alpha>0$, Lemmas~\ref{lem:ss-global} and
\ref{lem:ss-global-weaker} give
\begin{align}
        \bar E(u(s_n);\R^d\setminus B(0,\alpha\sqrt{s_n}))
        +
        \|u(s_n)\|_{L^{p+1}
        (\R^d\setminus B(0,\alpha\sqrt{s_n}))}^{p+1}
        \to0.
\end{align}
Choose integers $N_1<N_2<\cdots$ with $N_k\to\infty$ so that, for
$n\geq N_k$, the preceding sum with $\alpha=1/k$ is at most $1/k$ and
$\sqrt{s_n}/k\geq k$. Set $\alpha_n:=1/k$ for
$N_k\leq n<N_{k+1}$ and set the remaining finitely many $\alpha_n$
equal to $1$. Then $\alpha_n\to0$, $\alpha_n\sqrt{s_n}\to\infty$, and,
with $\rho_n:=\alpha_n\sqrt{s_n}$,
\begin{align}\label{eq:no-ghost-global-tail}
        \lim_{n\to \infty} \left[\bar E(u(s_n);\R^d\setminus B(0,\rho_n))
        +
        \|u(s_n)\|_{L^{p+1}(\R^d\setminus B(0,\rho_n))}^{p+1}\right]
        =0.
\end{align}
Choose $\zeta_n\in C_c^\infty(\R^d)$ such that $0\leq\zeta_n\leq1$,
$\zeta_n=1$ on $B(0,\rho_n)$,
$\supp\zeta_n\subset B(0,2\rho_n)$,
$|\nabla\zeta_n|\lesssim\rho_n^{-1}$, and
$|\Delta\zeta_n|\lesssim\rho_n^{-2}$. Set
$w_n:=\zeta_nu(s_n)$ and
$A_n:=B(0,2\rho_n)\setminus B(0,\rho_n)$. Then
\eqref{eq:no-ghost-global-tail} gives
\begin{align}\label{eq:no-ghost-global-transition-small}
        \lim_{n\to \infty} \left[\bar E(u(s_n);A_n)
        +
        \|u(s_n)\|_{L^{p+1}(A_n)}^{p+1}\right]
        =0.
\end{align}
Moreover,
\begin{align}
        \|u(s_n)-w_n\|_{\dot H^1}
        &\lesssim
        \|\nabla u(s_n)\|_{L^2(\R^d\setminus B(0,\rho_n))}
        +
        \|u(s_n)\|_{L^{p+1}(A_n)}
        \|\nabla\zeta_n\|_{L^d},
\end{align}
and hence
\begin{align}\label{eq:no-ghost-global-good-decomp}
        \lim_{n\to \infty}\left[\|u(s_n)-w_n\|_{\dot H^1}
        +
        \|u(s_n)-w_n\|_{L^{p+1}}\right]
        =0.
\end{align}
In particular, $\{w_n\}_{n\in\N}$ is bounded in
$\dot H^1(\R^d)$. Arguing as in the estimates leading up to
\eqref{eq:no-ghost-finite-residual-estimate}, with the present cutoff, gives
\begin{align}\label{eq:no-ghost-global-residual-estimate}
        \|\Delta w_n+|w_n|^{p-1}w_n\|_{\dot H^{-1}}
        &\lesssim
        \rho_n\|\partial_tu(s_n)\|_{L^2}
        +
         \|\nabla u(s_n)\|_{L^2(A_n)}
        \notag\\
        &\quad+
        \|u(s_n)\|_{L^{p+1}(A_n)}
        +
        \|u(s_n)\|_{L^{p+1}(A_n)}^p.
\end{align}
By \eqref{eq:no-ghost-global-good-time},
\begin{align}
       \lim_{n\to \infty} \rho_n^2\|\partial_tu(s_n)\|_{L^2}^2
        =
        \lim_{n\to \infty}\alpha_n^2s_n\|\partial_tu(s_n)\|_{L^2}^2
        =0.
\end{align}
Combining this with \eqref{eq:no-ghost-global-transition-small} and
\eqref{eq:no-ghost-global-residual-estimate}, we obtain
\begin{align}\label{eq:no-ghost-global-PS}
        \lim_{n\to \infty}\|\Delta w_n+|w_n|^{p-1}w_n\|_{\dot H^{-1}}
        =0.
\end{align}
After passing to a subsequence, Lemma~\ref{lem:whole-space-PS} and
\eqref{eq:no-ghost-global-good-decomp} give a finite sum
of stationary solutions whose difference from $u(s_n)$ converges to zero
in $\dot H^1(\R^d)$.

Set $b(t):=\|u(t)\|_{L^{p+1}}^{p+1}$. The function $b$ is continuous and bounded. Its set of limit points as $t\to T_+$ is the interval
\[
        [
        \liminf_{t\to T_+}b(t),
        \limsup_{t\to T_+}b(t)
        ].
\]
Set $M_0:=\sup_{0\leq t<T_+}\bar E(u(t)).$ If $T_+=\infty$, then $\|\nabla v_n\|_{L^2}^2\leq M_0$. If $T_+<\infty$, weak lower semicontinuity gives $\bar E(u^*)\leq M_0$, and hence
\begin{align}
        \|\nabla v_n\|_{L^2}^2
        =
        \|\nabla u(t_n)-\nabla u^*\|_{L^2}^2
        \leq
        2\bar E(u(t_n))+2\bar E(u^*)
        \leq
        4M_0.
\end{align}
Therefore, by \eqref{eq:arbitrary-stationary-pythagorean},
\begin{align}
        \sum_{j=1}^J
        \|\nabla W^j\|_{L^2}^2
        \leq
        4M_0
\end{align}
for every sequence of times. Fix $C>4M_0$. Let $\beta$ be a limit point of $b(t)$ as $t\to T_+$ and choose a sequence $t_n\to T_+$ such that $b(t_n)\to\beta$. If $T_+=\infty$, Lemma~\ref{lem:arbitrary-stationary-profiles} gives $\beta\in\mathscr M_C$. If $T_+<\infty$, then \eqref{eq:no-ghost-finite-BL-splitting} with $\tau_n=t_n$ gives $\beta-\|u^*\|_{L^{p+1}}^{p+1}\in\mathscr M_C$. Thus, the set of limit points is contained in $\mathscr M_C$ when $T_+=\infty$, and in $\|u^*\|_{L^{p+1}}^{p+1}+\mathscr M_C$ when $T_+<\infty$. Since in either case the set of limit points is a nonempty interval contained in a countable set, it consists of one point. Therefore, $b(t)$ has a limit as $t\to T_+$.

By Sobolev's inequality and the uniform $\dot H^1$ bound, $E(u(t))$ is
bounded below. Hence, by \eqref{eq:energy-identity}, it has a finite
limit as $t\to T_+$. When $T_+=\infty$, set $b_\infty:=\lim_{t\to\infty}b(t)$ and $E_\infty:=\lim_{t\to\infty}E(u(t))$. Along the Palais-Smale sequence constructed above, \eqref{eq:no-ghost-global-good-decomp} and  Lemma~\ref{lem:whole-space-PS} give
\begin{align}
        E(u(s_n))
        =
        \frac1d\|u(s_n)\|_{L^{p+1}}^{p+1}
        +
        o_n(1).
\end{align}
When $T_+<\infty$, set $b_{T_+}:=\lim_{t\to T_+}b(t)$, $E_{T_+}:=\lim_{t\to T_+}E(u(t))$, and $v_n:=u(s_n)-u^*$. By \eqref{eq:no-ghost-finite-good-decomp} and Lemma~\ref{lem:whole-space-PS}, $v_n$ is, up to a term converging to zero in $\dot H^1(\R^d)$, a finite sum of asymptotically pairwise orthogonal stationary solutions. The asymptotic orthogonality of these stationary solutions gives
\begin{align}
        E(v_n)
        =
        \frac1d\|v_n\|_{L^{p+1}}^{p+1}
        +
        o_n(1).
\end{align}
Applying \eqref{eq:no-ghost-finite-BL-splitting} and
\eqref{eq:no-ghost-finite-Dirichlet-splitting} with $\tau_n=s_n$ gives
\begin{align}
        E(u(s_n))-E(u^*)
        =
        E(v_n)+o_n(1),\text{ and }
        \|u(s_n)\|_{L^{p+1}}^{p+1}
        -
        \|u^*\|_{L^{p+1}}^{p+1}
        =
        \|v_n\|_{L^{p+1}}^{p+1}
        +
        o_n(1).
\end{align}
Passing to the limit in the above display gives
\begin{align}\label{eq:no-ghost-energy-good-sequence}
        E_\infty
        &=
        \frac1d b_\infty
        \quad
        \textrm{if }T_+=\infty,
        \notag\\
        E_{T_+}-E(u^*)
        &=
        \frac1d
        \left(
        b_{T_+}-\|u^*\|_{L^{p+1}}^{p+1}
        \right)
        \quad
        \textrm{if }T_+<\infty.
\end{align}
Passing to a subsequence realizing the $\limsup$, we may assume
that
\begin{align}
        \lim_{n\to\infty}\|\nabla r_n\|_{L^2}^2
        =
        \limsup_{n\to\infty}\|\nabla r_n\|_{L^2}^2
        =:\ga,
\end{align}
for some $\ga\in [0,\infty).$ If
$T_+=\infty$, then Lemma~\ref{lem:arbitrary-stationary-profiles} gives
\begin{align}\label{eq:no-ghost-energy-global-arbitrary}
        E_\infty
        =
        \frac1d
        \sum_{j=1}^J\|\nabla W^j\|_{L^2}^2
        +
        \frac12\gamma,\qquad 
        b_\infty=
        \sum_{j=1}^J\|\nabla W^j\|_{L^2}^2.
\end{align}
If $T_+<\infty$, set $v_n:=u(t_n)-u^*$. Applying
\eqref{eq:no-ghost-finite-BL-splitting} and
\eqref{eq:no-ghost-finite-Dirichlet-splitting} with $\tau_n=t_n$, we obtain
\begin{align}
        E_{T_+}-E(u^*)
        =
        \lim_{n\to\infty}
        \left[
        \frac12\|\nabla v_n\|_{L^2}^2
        -
        \frac1{p+1}\|v_n\|_{L^{p+1}}^{p+1}
        \right],\qquad
        b_{T_+}-\|u^*\|_{L^{p+1}}^{p+1}
        =
        \lim_{n\to\infty}
        \|v_n\|_{L^{p+1}}^{p+1}.
\end{align}
Using \eqref{eq:arbitrary-stationary-pythagorean},
$\|r_n\|_{L^{p+1}}\to0$, and $\|\nabla W^j\|_{L^2}^2=\|W^j\|_{L^{p+1}}^{p+1}$ (as the $W^j$ solve \eqref{eq:yamabe}), we get
\begin{align}\label{eq:no-ghost-energy-finite-arbitrary}
        E_{T_+}-E(u^*)
        =
        \frac1d
        \sum_{j=1}^J\|\nabla W^j\|_{L^2}^2
        +
        \frac12\gamma,
        \quad
        b_{T_+}-\|u^*\|_{L^{p+1}}^{p+1}
        =
        \sum_{j=1}^J\|\nabla W^j\|_{L^2}^2.
\end{align}
Comparing \eqref{eq:no-ghost-energy-good-sequence} with
\eqref{eq:no-ghost-energy-global-arbitrary} and
\eqref{eq:no-ghost-energy-finite-arbitrary}, we obtain $\ga=0$. Thus
$\limsup_{n\to\infty}\|\nabla r_n\|_{L^2}^2=0$, which proves the lemma.
\end{proof}
We now treat $d\geq5$. By
\eqref{eq:no-ghost-residual-identity}, it suffices to show that the
left-hand side tends to zero. In the global case, we prove
$\|\partial_tu(t)\|_{\dot H^{-1}}\to0$ as $t\to\infty$; we treat the
finite-time case separately.
% The argument is different when $d\geq 5$. The aim here is to estimate $\|\partial_t u\|_{\dot{H}^{-1}}$ since the energy of $r_n$ is roughly the same as that of the integral of the product of $\partial_t u$ and $r_n.$ Thus to show that the energy of $r_n$ vanishes, we will aim to show that the $\dot{H}^{-1}$ norm of $\partial_t u$ vanishes when $d\geq 5.$
\begin{lem}\label{lem:arbitrary-time-no-ghost-high-dimensional}Assume that $d\geq5$. Let $t_n\to T_+$ and let $r_n$ be the remainder in Lemma~\ref{lem:arbitrary-stationary-profiles}. Then $\|\nabla r_n\|_{L^2}\to 0$ as $n\to \infty.$
\end{lem}
\begin{proof}
Set $S_n:=\sum_{j=1}^J W^j[a_n^j,\lambda_n^j]$. By part
$\mathrm{(b)}$ of Proposition~\ref{prop:profile-decomp},
\begin{align}\label{eq:no-ghost-gradient-pairing}
        \int_{\R^d}\nabla u(t_n)\cdot\nabla r_n\,\ud x
        =
        \|\nabla r_n\|_{L^2}^2+o_n(1).
\end{align}
Indeed, when $T_+=\infty$, this follows by writing $u(t_n)=S_n+r_n$ in the above integral. When $T_+<\infty$, every profile converges weakly to zero in $\dot H^1(\R^d)$ by \eqref{eq:arbitrary-profile-finite-localization}. Since
$u(t_n)-u^*\rightharpoonup0$ weakly in $\dot H^1(\R^d)$, it follows that $r_n\rightharpoonup0$ weakly in $\dot H^1(\R^d)$, and hence $\int_{\R^d}\nabla u^*\cdot\nabla r_n\,\ud x=o_n(1)$.

Since $\lim_{n\to\infty}\|r_n\|_{L^{p+1}}=0$ and
$\{u(t_n)\}_{n\in\N}$ is bounded in $L^{p+1}(\R^d)$, we have
\begin{align}
        \lim_{n\to\infty}
        \int_{\R^d}
        |u(t_n)|^{p-1}u(t_n)r_n\,\ud x
        =
        0.
\end{align}
Therefore, using \eqref{eqn:NLH} and
\eqref{eq:no-ghost-gradient-pairing}, we obtain
\begin{align}\label{eq:no-ghost-residual-identity}
        &-\int_{\R^d}\nabla u(t_n)\cdot\nabla r_n\,\ud x
        +
        \int_{\R^d}|u(t_n)|^{p-1}u(t_n)r_n\,\ud x
        =
        -\|\nabla r_n\|_{L^2}^2+o_n(1).
\end{align}
We now claim that, for every $0<t<T<T_+$,
\begin{align}\label{eq:no-ghost-tail-action}
        \int_0^{T-t}
        \|e^{s\Delta}\partial_tu(t)\|_{L^2}^2\,\ud s
        \lesssim
        \int_t^T\|\partial_tu(s)\|_{L^2}^2\,\ud s.
\end{align}
Here the implicit constant depends only on $d$ and
$\sup_{0\leq\tau<T_+}\|\nabla u(\tau)\|_{L^2}$. Fix $h>0$ such that $T+h<T_+$, and set
$z_h(s):=h^{-1}[u(s+h)-u(s)]$. We first prove
\eqref{eq:no-ghost-tail-action} with $z_h(t)$ in place of
$\partial_tu(t)$ and then let $h\to0$. This avoids differentiating
\eqref{eqn:NLH} in time. Using \eqref{eqn:NLH}, we obtain
\begin{align}\label{eq:no-ghost-difference-quotient}
        \partial_sz_h
        =
        \Delta z_h+V_hz_h,
        \quad
        V_h(s,x)
        :=
        p\int_0^1
        |(1-\theta)u(s,x)+\theta u(s+h,x)|^{p-1}\,\ud\theta.
\end{align}
Since $(p-1)d/2=p+1$, Minkowski's inequality and Sobolev's inequality
give
\begin{align}
        \|V_h(s)\|_{L^{d/2}}
        \lesssim
        \|u(s)\|_{L^{p+1}}^{p-1}
        +
        \|u(s+h)\|_{L^{p+1}}^{p-1}
        \lesssim
        \sup_{0\leq\tau<T_+}
        \|\nabla u(\tau)\|_{L^2}^{p-1}
        \lesssim
        1.
\end{align}
Therefore,
\begin{align}\label{eq:no-ghost-potential-uniform}
        \sup_{0<h<T_+-T}
        \sup_{0<s<T}
        \|V_h(s)\|_{L^{d/2}}
        \lesssim
        1.
\end{align}
For a distribution $f\in\dot H^{-2}(\R^d)$, we use the Fourier
characterization
\begin{align}
        \|f\|_{\dot H^{-2}}^2
        :=
        \int_{\R^d}
        \frac{|\widehat f(\xi)|^2}{|\xi|^4}\,\ud\xi.
\end{align}
We next estimate the product $V_hg$ in $\dot H^{-2}(\R^d)$ for an arbitrary function $g\in L^2(\R^d).$ Let
$\varphi\in C_c^\infty(\R^d)$. Since $\frac{2}{d}+\frac12+\frac{d-4}{2d}=1,$ H\"older's inequality and the Sobolev embedding
$\dot H^2(\R^d)\subset L^{2d/(d-4)}(\R^d)$ give
\begin{align}
        \left|
        \int_{\R^d}V_hg\varphi\,\ud x
        \right|
        &\leq
        \|V_h\|_{L^{d/2}}
        \|g\|_{L^2}
        \|\varphi\|_{L^{2d/(d-4)}}
        \lesssim
        \|V_h\|_{L^{d/2}}
        \|g\|_{L^2}
        \|\Delta\varphi\|_{L^2}.
\end{align}
Taking the supremum over all
$\varphi\in C_c^\infty(\R^d)$ satisfying
$\|\Delta\varphi\|_{L^2}\leq1$, we obtain
\begin{align}\label{eq:no-ghost-potential-estimate}
        \|V_hg\|_{\dot H^{-2}}
        \lesssim
        \|V_h\|_{L^{d/2}}
        \|g\|_{L^2}.
\end{align}
Let $G\in L^2((0,T-t);\dot H^{-2}(\R^d))$, set
\begin{align}
        \cT G(s)
        :=
        \int_0^s
        e^{(s-\sigma)\Delta}G(\sigma)\,\ud\sigma.
\end{align}
We next prove
\begin{align}\label{eq:no-ghost-volterra}
        \|\cT G\|_{L_s^2(0,T-t;L_x^2)}
        \leq
        \|G\|_{L_s^2(0,T-t;\dot H_x^{-2})}.
\end{align}
By density, it suffices to prove this estimate for smooth $G$. Taking the
Fourier transform in the spatial variable gives
\begin{align}
        \widehat{\cT G}(s,\xi)
        =
        \int_0^s
        e^{-(s-\sigma)|\xi|^2}
        \widehat G(\sigma,\xi)\,\ud\sigma.
\end{align}
Extending $\widehat G(\cdot,\xi)$ by zero outside $(0,T-t)$ and applying
Young's inequality in time with
$k_\xi(s):=e^{-s|\xi|^2}\mathbf{1}_{\{s>0\}}$ gives
\begin{align}
        \|\widehat{\cT G}(\cdot,\xi)\|_{L_s^2(0,T-t)}
        &\leq
        \|k_\xi\|_{L^1(\R)}
        \|\widehat G(\cdot,\xi)\|_{L_s^2(0,T-t)}
        =
        \frac{1}{|\xi|^2}
        \|\widehat G(\cdot,\xi)\|_{L_s^2(0,T-t)}.
\end{align}
Therefore, by Plancherel's theorem and Tonelli's theorem,
\begin{align}
        \|\cT G\|_{L_s^2(0,T-t;L_x^2)}^2
        &=
        \int_{\R^d}
        \|\widehat{\cT G}(\cdot,\xi)\|_{L_s^2(0,T-t)}^2
        \,\ud\xi
        \notag\\
        &\leq
        \int_{\R^d}
        \frac{1}{|\xi|^4}
        \|\widehat G(\cdot,\xi)\|_{L_s^2(0,T-t)}^2
        \,\ud\xi
        \notag\\
        &=
        \|G\|_{L_s^2(0,T-t;\dot H_x^{-2})}^2.
\end{align}
Taking square roots on both sides implies \eqref{eq:no-ghost-volterra}. By \eqref{eq:lwp-positive-time-H2}, $z_h\in L_s^2(t,T;L_x^2)$, and hence
$V_hz_h\in L_s^2(t,T;\dot H_x^{-2})$ by \eqref{eq:no-ghost-potential-estimate}. Thus Duhamel's formula in the sense of distributions gives, for every $0\leq s\leq T-t$,
\begin{align}
        e^{s\Delta}z_h(t)
        =
        z_h(t+s)
        -
        \int_0^s
        e^{(s-\sigma)\Delta}
        V_h(t+\sigma)z_h(t+\sigma)\,\ud\sigma.
\end{align}
It follows from \eqref{eq:no-ghost-potential-estimate} and
\eqref{eq:no-ghost-volterra} that
\begin{align}
        \|e^{s\Delta}z_h(t)\|_{L_s^2(0,T-t;L_x^2)}
        \lesssim
        \|z_h\|_{L_s^2(t,T;L_x^2)}.
\end{align}
By \eqref{eq:lwp-positive-time-H2},
$u(s+h)-u(s)=\int_s^{s+h}\partial_\tau u(\tau)\,\ud\tau$ in
$L^2(\R^d)$. Consequently,
\begin{align}
        \lim_{h\to0}
        \|z_h-\partial_tu\|_{L_s^2(t,T;L_x^2)}
        =
        0.
\end{align}
Moreover, \eqref{eq:lwp-positive-time-H2} gives
$z_h(t)\to \partial_tu(t)$ in $L^2(\R^d)$. Passing to the limit as $h\to0$ proves \eqref{eq:no-ghost-tail-action}.

We now prove that $r_n\to0$ in $\dot H^1(\R^d)$ as $n\to \infty.$ Assume first that $T_+=\infty$. By \eqref{eqn:NLH}, H\"older's inequality, and Sobolev's inequality, $\partial_tu(t)\in\dot H^{-1}(\R^d)$ for every $t>0$. For every $f\in\dot H^{-1}(\R^d)$, Plancherel's theorem gives
\begin{align}\label{eq:no-ghost-free-identity}
        \int_0^\infty
        \|e^{s\Delta}f\|_{L^2}^2\,\ud s
        =
        \int_0^\infty
        \int_{\R^d}
        e^{-2s|\xi|^2}
        |\widehat f(\xi)|^2\,\ud\xi\ud s
        =
        \frac12
        \int_{\R^d}
        \frac{|\widehat f(\xi)|^2}{|\xi|^2}\,\ud\xi
        =
        \frac12\|f\|_{\dot H^{-1}}^2.
\end{align}
Letting $T\to\infty$ in \eqref{eq:no-ghost-tail-action}, we obtain
\begin{align}\label{eq:no-ghost-global-residual}
        \|\partial_tu(t)\|_{\dot H^{-1}}^2
        \lesssim
        \int_t^\infty\|\partial_tu(s)\|_{L^2}^2\,\ud s
        \to0
        \qquad
        \textrm{as }t\to\infty.
\end{align}
By the definition of the $\dot H^{-1}$-norm and the uniform boundedness
of $r_n$ in $\dot H^1(\R^d)$,
\begin{align}
        &\left|
        -\int_{\R^d}\nabla u(t_n)\cdot\nabla r_n\,\ud x
        +
        \int_{\R^d}|u(t_n)|^{p-1}u(t_n)r_n\,\ud x
        \right|
        \leq
        \|\partial_tu(t_n)\|_{\dot H^{-1}}
        \|\nabla r_n\|_{L^2}.
\end{align}
Thus, \eqref{eq:no-ghost-global-residual} and
\eqref{eq:no-ghost-residual-identity} give
$\lim_{n\to\infty}\|\nabla r_n\|_{L^2}=0$.

Assume now that $T_+<\infty$, and set
$\rho_n:=\sqrt{T_+-t_n}$. For every
$f\in\dot H^{-1}(\R^d)$,
\begin{align}
        \int_0^{\rho_n^2}
        \|e^{s\Delta}f\|_{L^2}^2\,\ud s
        &=
        \int_{\R^d}
        \frac{1-e^{-2\rho_n^2|\xi|^2}}{2|\xi|^2}
        |\widehat f(\xi)|^2\,\ud\xi
        \gtrsim
        \|P_{>\rho_n^{-1}}f\|_{\dot H^{-1}}^2.
\end{align}
Letting $T\to T_+$ in \eqref{eq:no-ghost-tail-action}, we deduce that
\begin{align}\label{eq:no-ghost-finite-high-residual}
        \|P_{>\rho_n^{-1}}\partial_tu(t_n)\|_{\dot H^{-1}}^2
        \lesssim
        \int_{t_n}^{T_+}\|\partial_tu(s)\|_{L^2}^2\,\ud s
        \to0.
\end{align}
For every fixed $\alpha>0$,
\eqref{eq:arbitrary-profile-finite-localization} and a change of
variables give
\begin{align}
        \bar E\left(
        W^j[a_n^j,\lambda_n^j];
        \R^d\setminus B(x^{i_j},\alpha\rho_n)
        \right)
        \to0
\end{align}
as $n\to \infty$ for every $1\leq j\leq J$. Hence
\eqref{eq:arbitrary-profile-finite-concentration}, the triangle
inequality, and
$r_n=v_n-\sum_{j=1}^J W^j[a_n^j,\lambda_n^j]$ give
\begin{align}
        \bar E\left(
        r_n;
        \R^d\setminus\cup_{i=1}^K B(x^i,\alpha\rho_n)
        \right)
        \to0,
\end{align}
as $n\to \infty$. Together with $\|r_n\|_{L^{p+1}}\to0$ as $n\to \infty$, a diagonal argument gives a
sequence $\alpha_n>0$ such that $\lim_{n\to\infty}\alpha_n=0$ and
\begin{align}\label{eq:no-ghost-finite-localization}
        \lim_{n\to\infty}
        \left[
        \bar E\left(
        r_n;
        \R^d\setminus
        \cup_{i=1}^K B(x^i,\alpha_n\rho_n)
        \right)
        +
        \|r_n\|_{L^{p+1}
        (\R^d\setminus
        \cup_{i=1}^K B(x^i,\alpha_n\rho_n))}^{p+1}
        \right]
        =
        0.
\end{align}
Choose $\chi_n\in C_c^\infty(\R^d)$ such that $\chi_n=1$ on
$\cup_{i=1}^K B(x^i,\alpha_n\rho_n)$, $ \supp\chi_n\subset \cup_{i=1}^K B(x^i,2\alpha_n\rho_n),$ $0\leq \chi_n\leq 1$ and $\|\nabla\chi_n\|_{L^d}\lesssim1$. Set $f_n:=\chi_nr_n$. Since
\begin{align}
        \nabla(r_n-f_n)
        =
        (1-\chi_n)\nabla r_n-r_n\nabla\chi_n
\end{align}
and
\begin{align}
        \|r_n\nabla\chi_n\|_{L^2}
        \leq
        \|r_n\|_{L^{p+1}}\|\nabla\chi_n\|_{L^d}
        \to0,
\end{align}
\eqref{eq:no-ghost-finite-localization} gives
\begin{align}\label{eq:no-ghost-finite-cutoff}
        \lim_{n\to\infty}\|r_n-f_n\|_{\dot H^1}
        =0.
\end{align}
For all sufficiently large $n$, the balls
$B(x^i,2\alpha_n\rho_n)$ are pairwise disjoint. Applying
Poincar\'e's inequality on each ball and using the uniform boundedness
of $\|\nabla f_n\|_{L^2}$ gives 
\begin{align}
        \|f_n\|_{L^2}
        \lesssim
        \alpha_n\rho_n\|\nabla f_n\|_{L^2}
        \lesssim
        \alpha_n\rho_n.
\end{align}
Consequently, Bernstein's inequality gives
\begin{align}\label{eq:no-ghost-finite-low-remainder}
\|\nabla P_{\leq\rho_n^{-1}}f_n\|_{L^2}
\lesssim
\rho_n^{-1}\|f_n\|_{L^2}
\lesssim
\alpha_n
\to0,
\end{align}
as $n\to \infty.$ Next, \eqref{eqn:NLH}, H\"older's inequality, and Sobolev's inequality give
\begin{align}\label{eq:no-ghost-uniform-negative-bound}
        \sup_{0<t<T_+}
        \|\partial_tu(t)\|_{\dot H^{-1}}
        &\lesssim
        \sup_{0<t<T_+}
        \left(
        \|\nabla u(t)\|_{L^2}
        +
        \|u(t)\|_{L^{p+1}}^p
        \right)
        \notag\\
        &\lesssim
        \sup_{0<t<T_+}
        \left(
        \|\nabla u(t)\|_{L^2}
        +
        \|\nabla u(t)\|_{L^2}^p
        \right)
        <\infty.
\end{align}
For $g\in\dot H^1(\R^d)$, set
\begin{align}
        \cR_n(g)
        &:=
        -\int_{\R^d}\nabla u(t_n)\cdot\nabla g\,\ud x
        +
        \int_{\R^d}|u(t_n)|^{p-1}u(t_n)g\,\ud x.
\end{align}
By \eqref{eqn:NLH}, $\cR_n(g)=\partial_tu(t_n)[g]$. Since $r_n=f_n+(r_n-f_n)$ and
$f_n=P_{>\rho_n^{-1}}f_n+P_{\leq\rho_n^{-1}}f_n$, we have
\begin{align}
        |\cR_n(r_n)| &=\left|
        -\int_{\R^d}\nabla u(t_n)\cdot\nabla r_n\,\ud x
        +
        \int_{\R^d}|u(t_n)|^{p-1}u(t_n)r_n\,\ud x
        \right|
        \notag\\
        &\leq
        |\cR_n(f_n)|
        +
        |\cR_n(r_n-f_n)|
        \notag\\
        &\leq
        |\cR_n(P_{>\rho_n^{-1}}f_n)|
        +
        |\cR_n(P_{\leq\rho_n^{-1}}f_n)|
        +
        |\cR_n(r_n-f_n)|
        \notag\\
        &\leq \|P_{>\rho_n^{-1}}\partial_tu(t_n)\|_{\dot H^{-1}}
        \|\nabla f_n\|_{L^2}
        +
        \|\partial_tu(t_n)\|_{\dot H^{-1}}
        \|\nabla P_{\leq\rho_n^{-1}}f_n\|_{L^2}\\
        &\quad +
        \|\partial_tu(t_n)\|_{\dot H^{-1}}
        \|\nabla(r_n-f_n)\|_{L^2}.
\end{align}
In the first term, we used
\begin{align}
        \cR_n(P_{>\rho_n^{-1}}f_n)
        =
        (P_{>\rho_n^{-1}}\partial_tu(t_n))[f_n].
\end{align}
It follows from \eqref{eq:no-ghost-finite-high-residual},
\eqref{eq:no-ghost-finite-cutoff},
\eqref{eq:no-ghost-finite-low-remainder}, and
\eqref{eq:no-ghost-uniform-negative-bound} that $\cR_n(r_n)\to0$ as $n\to \infty$. Combining this with
\eqref{eq:no-ghost-residual-identity} gives
$\|\nabla r_n\|_{L^2}\to0$ as $n\to\infty$.
\end{proof}
\begin{proof}[Proof of Lemma~\ref{lem:arbitrary-time-no-ghost}]
The conclusion now follows from  Lemma~\ref{lem:arbitrary-time-no-ghost-high-dimensional} when $d\geq 5$ and
Lemma~\ref{lem:arbitrary-time-no-ghost-low-dimensional} when $d\in \{3,4\}$.
\end{proof}
We next show that, when $T_+<\infty$, each point in the singular set is associated to at least one profile in the profile decomposition obtained in Lemma~\ref{lem:arbitrary-stationary-profiles}.
\begin{lem}
\label{lem:arbitrary-singular-coverage}
Assume that $T_+<\infty$. Let $t_n\to T_+$ and suppose that
\begin{align}\label{eq:arbitrary-coverage-decomp}
        u(t_n)
        =
        u^*
        +
        \sum_{j=1}^J W^j[a_n^j,\lambda_n^j]
        +
        o_{\dot H^1}(1),
\end{align}
where $W^1,\ldots,W^J$ are nonzero stationary solutions and, for every
$1\leq j\leq J$, there exists $x^{i_j}\in\calS$ such that the
parameters satisfy \eqref{eq:arbitrary-profile-finite-localization}. Then, for every
$x^i\in\calS$, there exists $j\in\{1,\ldots,J\}$ such that
$x^{i_j}=x^i$.
\end{lem}

\begin{proof}
Fix $x^i\in\calS$ and suppose that no profile is associated to $x^i$.
Since $\calS$ is finite, we can choose a decreasing sequence $r_m>0$
such that
\begin{align}\label{eq:arbitrary-coverage-radii}
        \lim_{m\to\infty}r_m
        =
        0,
        \qquad
        0<4r_m\leq R_0,
        \qquad
        \overline{B(x^i,4r_m)}
        \cap
        (\calS\setminus\{x^i\})
        =
        \emptyset
\end{align}
for every $m$. By the absolute continuity of the $\dot H^1$ and
$L^{p+1}$ norms of $u^*$,
\begin{align}\label{eq:arbitrary-coverage-regular-small}
        \lim_{m\to\infty}
        [
        \bar E(u^*;B(x^i,4r_m))
        +
        \|u^*\|_{L^{p+1}(B(x^i,4r_m))}^{p+1}
        ]
        =
        0.
\end{align}

Fix $m$. Suppose first that
$\calS\setminus\{x^i\}\neq\emptyset$, and set
\[
        \delta_{i,m}
        :=
        \dist(
        \overline{B(x^i,4r_m)},
        \calS\setminus\{x^i\}
        ).
\]
By \eqref{eq:arbitrary-coverage-radii}, we have $\delta_{i,m}>0$. Since no
profile is associated to $x^i$, we have
$x^{i_j}\in\calS\setminus\{x^i\}$ for every $1\leq j\leq J$.
Therefore, \eqref{eq:arbitrary-profile-finite-localization} gives
\[
        \dist(a_n^j,B(x^i,4r_m))
        \geq
        \frac{\delta_{i,m}}{2}
\]
for every $1\leq j\leq J$ and all sufficiently large $n$. Moreover,
$\lim_{n\to\infty}\lambda_n^j=0$. Hence, after changing variables,
\begin{align}
        \bar E(
        W^j[a_n^j,\lambda_n^j];
        B(x^i,4r_m)
        )
        &\leq
        \int_{|y|\geq \delta_{i,m}/(2\lambda_n^j)}
        |\nabla W^j(y)|^2\,\ud y,
        \notag\\
        \|W^j[a_n^j,\lambda_n^j]\|_
        {L^{p+1}(B(x^i,4r_m))}^{p+1}
        &\leq
        \int_{|y|\geq \delta_{i,m}/(2\lambda_n^j)}
        |W^j(y)|^{p+1}\,\ud y.
\end{align}
Since $W^j\in\dot H^1(\R^d)\cap L^{p+1}(\R^d)$, both terms converge to
zero as $n\to\infty$.

If $\calS\setminus\{x^i\}=\emptyset$, then every profile is associated
to $x^i$ by \eqref{eq:arbitrary-profile-finite-localization}. Thus, the
contradiction assumption implies that $J=0$. In either case, since the
number of profiles is finite, \eqref{eq:arbitrary-coverage-decomp} gives,
for every fixed $m$,
\begin{align}\label{eq:arbitrary-coverage-defect-small}
        \lim_{n\to\infty}
        [
        \bar E(u(t_n)-u^*;B(x^i,4r_m))
        +
        \|u(t_n)-u^*\|_
        {L^{p+1}(B(x^i,4r_m))}^{p+1}
        ]
        =
        0.
\end{align}

We can now choose increasing integers $n_m\in \N$ such that
\begin{align}\label{eq:arbitrary-coverage-diagonal-choice}
        &\bar E(u(t_{n_m})-u^*;B(x^i,4r_m))
        +
        \|u(t_{n_m})-u^*\|_
        {L^{p+1}(B(x^i,4r_m))}^{p+1}
        \leq
        \frac1m,\quad 
        0<T_+-t_{n_m}
        \leq
        \frac{r_m^2}{m}.\qquad 
\end{align}
By \eqref{eq:arbitrary-coverage-regular-small},
\eqref{eq:arbitrary-coverage-diagonal-choice}, and the inequalities
$|\nabla u|^2\leq2|\nabla(u-u^*)|^2+2|\nabla u^*|^2$ and
$|u|^{p+1}\lesssim|u-u^*|^{p+1}+|u^*|^{p+1}$, we obtain
\begin{align}\label{eq:arbitrary-coverage-initial-small}
        \lim_{m\to\infty}
        [
        \bar E(u(t_{n_m});B(x^i,2r_m))
        +
        \|u(t_{n_m})\|_
        {L^{p+1}(B(x^i,2r_m))}^{p+1}
        ]
        =
        0.
\end{align}
Moreover,
\begin{align}\label{eq:arbitrary-coverage-time-scale}
        \lim_{m\to\infty}
        \frac{T_+-t_{n_m}}{r_m^2}
        =
        0.
\end{align}

Since $x^i\in\calS$ and $r_m\leq R_0$,
\eqref{defn:singular-points} implies that there exists $s_m\in(t_{n_m},T_+)$ such that
\begin{align}\label{eq:arbitrary-coverage-later-lower}
        \bar E(u(s_m);B(x^i,r_m))
        +
        \|u(s_m)\|_{L^{p+1}(B(x^i,r_m))}^{p+1}
        \geq
        \frac{\veps_0}{2}.
\end{align}
Since $t_{n_m},s_m\to T_+$, we apply
Lemma~\ref{lem:short-time-prop} to $u(t,x+x^i)$ with $\sigma_m:=t_{n_m}$, $\ta_m:=s_m$, $W:=0$, and radius
$r_m$. By \eqref{eq:arbitrary-coverage-time-scale},
\begin{align}
        0
        \leq
        \frac{s_m-t_{n_m}}{r_m^2}
        \leq
        \frac{T_+-t_{n_m}}{r_m^2},
        \quad
        \lim_{m\to\infty}
        \frac{s_m-t_{n_m}}{r_m^2}
        =
        0.
\end{align}
Therefore, \eqref{eq:arbitrary-coverage-initial-small} and
Lemma~\ref{lem:short-time-prop} give
\begin{align}
        \lim_{m\to\infty}
        [
        \bar E(u(s_m);B(x^i,r_m))
        +
        \|u(s_m)\|_{L^{p+1}(B(x^i,r_m))}^{p+1}
        ]
        =
        0.
\end{align}
This contradicts \eqref{eq:arbitrary-coverage-later-lower}. Hence there
exists $j\in\{1,\ldots,J\}$ such that $x^{i_j}=x^i$.
\end{proof}
The following lemma shows that if a profile decomposition consists of stationary profiles, its remainder vanishes in $\dot H^1$ and $L^{p+1}$ on $B(y_n,\beta_n\rho_n)$, and the full sequence vanishes on $B(y_n,\beta_n\rho_n)\setminus B(y_n,\alpha_n\rho_n)$, then the localized distance $\bs\de_{\gamma_0}(v_n;B(y_n,\rho_n))$ to a multi-bubble configuration tends to zero.
\begin{lem}
\label{lem:localize-strong-bubbles}
Fix $\gamma_0\in(0,\bar E_*/2)$. Let $J\geq0$, let
$W^1,\ldots,W^J$ be fixed nonzero stationary solutions, and let $(a_n^j,\lambda_n^j)$ satisfy \eqref{eq:arbitrary-profile-orthogonality} for every
$j,k\in\{1,\ldots,J\}$ with $j\neq k$. Set $Q_n^j:=W^j[a_n^j,\lambda_n^j]$. Let $e_n\in\dot H^1(\R^d)$ and set
\begin{align}
        v_n:=\sum_{j=1}^JQ_n^j+e_n.
\end{align}
Let $y_n\in\R^d$, $\rho_n,\alpha_n,\beta_n>0$,
$\alpha_n\to0$, and $\beta_n\to\infty$ as $n\to \infty.$ Suppose
\begin{align}\label{eq:localize-strong-error}
        \lim_{n\to\infty}
        \left[
        \bar E(e_n;B(y_n,\beta_n\rho_n))
        +
        \|e_n\|_{L^{p+1}(B(y_n,\beta_n\rho_n))}^{p+1}
        \right]
        =0,
\end{align}
and
\begin{align}\label{eq:localize-strong-neck}
        \lim_{n\to\infty}
        \left[
        \bar E(v_n;B(y_n,\beta_n\rho_n)\setminus
        B(y_n,\alpha_n\rho_n))
        +
        \|v_n\|_{L^{p+1}(B(y_n,\beta_n\rho_n)\setminus
        B(y_n,\alpha_n\rho_n))}^{p+1}
        \right]
        =0.
\end{align}
Then
\begin{align}\label{eq:localize-strong-conclusion}
        \lim_{n\to \infty}\bs\de_{\gamma_0}(v_n;B(y_n,\rho_n))=0.
\end{align}
\end{lem}
\begin{proof}
It suffices to prove that every subsequence has a further subsequence
along which \eqref{eq:localize-strong-conclusion} holds.
After passing to a further subsequence, for every $1\leq j\leq J$ the limit
\begin{align}
        m_j
        :=
        \lim_{n\to\infty}
        \|Q_n^j\|_{L^{p+1}(B(y_n,2\rho_n))}^{p+1}
\end{align}
exists. Set $\calJ:=\{j\in\{1,\ldots,J\}:m_j>0\}$. For all sufficiently large $n\in \N$, choose
$\phi_n\in C_c^\infty(\R^d)$ such that $0\leq\phi_n\leq1$,
$\phi_n=1$ on
$B(y_n,2\rho_n)\setminus B(y_n,2\alpha_n\rho_n)$, $ \supp\phi_n \subset B(y_n,3\rho_n)\setminus B(y_n,\alpha_n\rho_n)$, and $\|\nabla\phi_n\|_{L^d}\lesssim1$. Since
$\lim_{n\to\infty}\beta_n=\infty$, the support of $\phi_n$ is contained
in
$B(y_n,\beta_n\rho_n)\setminus B(y_n,\alpha_n\rho_n)$ for all
sufficiently large $n$.  Set
\begin{align}
        S_n:=\sum_{j=1}^JQ_n^j.
\end{align}
By \eqref{eqn:weighted-product-orth} and
\eqref{eq:weighted-dec-pointwise-direct} with zero remainder,
\begin{align}\label{eq:localize-profile-Lp-decoupling}
        \int_{\R^d}|S_n|^{p+1}\phi_n^2\,\ud x
        =
        \sum_{j=1}^J
        \int_{\R^d}|Q_n^j|^{p+1}\phi_n^2\,\ud x
        +
        o_n(1).
\end{align}
Since $J$ is fixed and
$\|Q_n^j\|_{L^{p+1}}=\|W^j\|_{L^{p+1}}$, the sequence
$\{S_n\}_{n\in\N}$ is bounded in $L^{p+1}(\R^d)$. Therefore,
H\"older's inequality and \eqref{eq:localize-strong-error} give
\begin{align}
        \left|
        \int_{\R^d}
        \left(
        |v_n|^{p+1}-|S_n|^{p+1}
        \right)
        \phi_n^2\,\ud x
        \right|
        &\lesssim
        \|S_n\|_{L^{p+1}}^p
        \|e_n\|_{L^{p+1}(\supp\phi_n)}+
        \|e_n\|_{L^{p+1}(\supp\phi_n)}^{p+1}
        \to0.
\end{align}
On the other hand, \eqref{eq:localize-strong-neck} and the support
property of $\phi_n$ imply
\begin{align}
        \int_{\R^d}|v_n|^{p+1}\phi_n^2\,\ud x
        \to0.
\end{align}
Combining the preceding three displays and using that $\phi_n=1$ on
$B(y_n,2\rho_n)\setminus B(y_n,2\alpha_n\rho_n)$ gives
\begin{align}\label{eq:localize-profile-neck}
        \lim_{n\to\infty}
        \sum_{j=1}^J
        \|Q_n^j\|_{L^{p+1}
        (B(y_n,2\rho_n)\setminus
        B(y_n,2\alpha_n\rho_n))}^{p+1}
        =
        0.
\end{align}
Consequently, for every $j\in\calJ$,
\begin{align}
        \lim_{n\to\infty}
        \|Q_n^j\|_{L^{p+1}
        (B(y_n,2\alpha_n\rho_n))}^{p+1}
        =
        m_j>0.
\end{align}
Writing this ball in the $(a_n^j,\lambda_n^j)$-frame and using
$W^j\in L^{p+1}(\R^d)$, we obtain
\begin{align}\label{eq:localize-strong-parameter-bound}
        |a_n^j-y_n|+\lambda_n^j
        \lesssim_j
        \alpha_n\rho_n
        \qquad
        \textrm{for every }j\in\calJ.
\end{align}
Indeed, for every $j\in\calJ$, the change of variables
$x=a_n^j+\lambda_n^jz$ gives
\begin{align}
\|Q_n^j\|_{L^{p+1}(B(y_n,2\alpha_n\rho_n))}^{p+1}
=
\int_{B((y_n-a_n^j)/\lambda_n^j,\,
2\alpha_n\rho_n/\lambda_n^j)}
|W^j(z)|^{p+1}\,\ud z.
\end{align}
Since this quantity tends to $m_j>0$ and $|W^j|^{p+1}\in L^1(\R^d)$, the absolute continuity of the integral and integrability at infinity give, respectively,
\begin{align}
\liminf_{n\to\infty}
\frac{\alpha_n\rho_n}{\lambda_n^j}>0,
\qquad
\limsup_{n\to\infty}
\frac{|a_n^j-y_n|-2\alpha_n\rho_n}{\lambda_n^j}<\infty.
\end{align}
Consequently, \eqref{eq:localize-strong-parameter-bound} holds.

We next show that the profiles with index outside $\calJ$ make no
contribution on $B(y_n,\rho_n)$. Fix $j\notin\calJ$. By the definition
of $m_j$,
\begin{align}
        \lim_{n\to\infty}
        \|Q_n^j\|_{L^{p+1}(B(y_n,2\rho_n))}
        =
        0.
\end{align}
Choose $\eta_n\in C_c^\infty(B(y_n,2\rho_n))$ such that
$\eta_n=1$ on $B(y_n,\rho_n)$ and
$|\nabla\eta_n|\lesssim\rho_n^{-1}$. Testing the stationary equation for
$Q_n^j$ against $Q_n^j\eta_n^2$ and applying Young's inequality gives
\begin{align}
        \bar E(Q_n^j;B(y_n,\rho_n))
        &\lesssim
        \|Q_n^j\|_{L^{p+1}(B(y_n,2\rho_n))}^{p+1}+
        \rho_n^{-2}
        \|Q_n^j\|_{L^2(B(y_n,2\rho_n))}^2.
\end{align}
Since $p+1=2d/(d-2)$, H\"older's inequality gives
\begin{align}
        \rho_n^{-2}
        \|Q_n^j\|_{L^2(B(y_n,2\rho_n))}^2
        \lesssim
        \|Q_n^j\|_{L^{p+1}(B(y_n,2\rho_n))}^2.
\end{align}
Therefore,
\begin{align}\label{eq:localize-strong-omitted}
        \lim_{n\to\infty}
        \left[
        \bar E(Q_n^j;B(y_n,\rho_n))
        +
        \|Q_n^j\|_{L^{p+1}(B(y_n,\rho_n))}^{p+1}
        \right]
        =
        0
\end{align}
for every $j\notin\calJ$. Since the number of profiles is finite,
\eqref{eq:localize-strong-error} and
\eqref{eq:localize-strong-omitted} imply
\begin{align}\label{eq:localize-strong-approximation}
        \lim_{n\to\infty}
        \left[\bar E(
        v_n-\sum_{j\in\calJ}Q_n^j;
        B(y_n,\rho_n)
        )
        +
        \|
        v_n-\sum_{j\in\calJ}Q_n^j
        \|_{L^{p+1}(B(y_n,\rho_n))}^{p+1}\right]
        =
        0.
\end{align}
Choose a sequence $L_n>0$ such that
\begin{align}
        \lim_{n\to\infty}L_n
        =
        \infty,
        \quad
        \lim_{n\to\infty}L_n\alpha_n
        =
        0,
        \quad
        \lim_{n\to\infty}\frac{L_n}{\beta_n}
        =
        0.
\end{align}
Set $\xi_n:=L_n\alpha_n\rho_n$ and
$\nu_n:=(\beta_n\rho_n)/L_n$. Then $\xi_n/\rho_n\to0$ and $\rho_n/\nu_n\to0$ as $n\to \infty.$ Moreover,
\[
        B(y_n,\nu_n)\setminus B(y_n,\xi_n)
        \subset
        B(y_n,\beta_n\rho_n)\setminus
        B(y_n,\alpha_n\rho_n).
\]
Thus \eqref{eq:localize-strong-neck} gives
\begin{align}\label{eq:localize-reduced-neck}
        \lim_{n\to\infty}
        \left[
        \bar E(v_n;B(y_n,\nu_n)\setminus B(y_n,\xi_n))
        +
        \|v_n\|_{L^{p+1}
        (B(y_n,\nu_n)\setminus B(y_n,\xi_n))}^{p+1}
        \right]
        =
        0.
\end{align}

If $\calJ=\emptyset$, take the zero multi-bubble configuration and set
$\vec\nu_n:=(\nu_n)$ and $\vec\xi_n:=(\xi_n)$. Then
\eqref{eq:localize-strong-approximation},
\eqref{eq:localize-reduced-neck}, and the preceding limits show directly from Definitions~\ref{def:d} and \ref{defn:delta} that
\[
        \lim_{n\to\infty}
        \bs\de_{\gamma_0}(v_n;B(y_n,\rho_n))
        =
        0.
\]
We may therefore assume that $\calJ\neq\emptyset$.

For every $j\in\calJ$, set $b_{j,n}:=a(Q_n^j)$ and
$\mu_{j,n}:=\lambda(Q_n^j)$, where the center and scale are defined using
the fixed number $\gamma_0$. By Lemma~\ref{lem:well-defn-scale-center},
\begin{align}
        \mu_{j,n}
        =
        \lambda(W^j)\lambda_n^j,
        \quad 
        |b_{j,n}-a_n^j-a(W^j)\lambda_n^j|
        \leq
        2\lambda(W^j)\lambda_n^j.
\end{align}
Together with \eqref{eq:localize-strong-parameter-bound}, this gives
\begin{align}\label{eq:localize-global-geometry}
        \lim_{n\to\infty}
        \max_{j\in\calJ}
        \frac{|b_{j,n}-y_n|+\mu_{j,n}}{\xi_n}
        =
        0.
\end{align}
In particular, $b_{j,n}\in B(y_n,\xi_n)$ for every $j\in\calJ$ and all
sufficiently large $n$. Choose $\kappa_n\in(0,1)$ such that
$\lim_{n\to\infty}\kappa_n=0$ and
\begin{align}\label{eq:localize-kappa-choice}
        \lim_{n\to\infty}
        \max_{j\in\calJ}
        \frac{\mu_{j,n}}{\kappa_n\xi_n}
        =
        0.
\end{align}
Set $\xi_{j,n}:=\kappa_n\mu_{j,n}$. After discarding finitely many terms, for every $j\in\calJ$ choose
$\nu_{j,n}\in [\kappa_n^{-1}\mu_{j,n}, 2\kappa_n^{-1}\mu_{j,n}]$ such that
\begin{align}\label{eq:localize-boundary-gap}
        |\nu_{j,n}-|b_{k,n}-b_{j,n}||
        \geq
        \kappa_n^{1/2}\mu_{j,n}
        \qquad
        \textrm{for every }k\in\calJ\setminus\{j\}.
\end{align}
Such a choice is possible because $\calJ$ is finite. Indeed, the interval
from which $\nu_{j,n}$ is chosen has length
$\kappa_n^{-1}\mu_{j,n}$, while the union of the excluded intervals has
length at most $2|\calJ|\kappa_n^{1/2}\mu_{j,n}$.

By \eqref{eq:localize-global-geometry} and 
\eqref{eq:localize-kappa-choice}, $ B(b_{j,n},\nu_{j,n}) \subset B(y_n,\xi_n)$ for every $j\in\calJ$ and all sufficiently large $n$. Moreover,
\begin{align}\label{eq:localize-component-geometry}
        \lim_{n\to\infty}
        \left[
        \frac{\mu_{j,n}}
        {\dist(b_{j,n},\partial B(y_n,\xi_n))}
        +
        \frac{\mu_{j,n}}{\nu_{j,n}}
        +
        \frac{\xi_{j,n}}{\mu_{j,n}}
        \right]
        =
        0
\end{align}
for every $j\in\calJ$. For every $j\in\calJ$, set
\begin{align}
        I_{j,n}
        :=
        \{
        k\in\calJ\setminus\{j\}:
        B(b_{k,n},\xi_{j,n})
        \subset
        B(b_{j,n},\nu_{j,n})
        \}.
\end{align}
If $k\in I_{j,n}$, then
\[
        \dist(b_{k,n},\partial B(b_{j,n},\nu_{j,n}))
        =
        |\nu_{j,n}-|b_{k,n}-b_{j,n}||.
\]
Hence \eqref{eq:localize-boundary-gap} gives
\begin{align}\label{eq:localize-component-boundary}
        \lim_{n\to\infty}
        \sum_{j\in\calJ}
        \sum_{k\in I_{j,n}}
        \frac{\xi_{j,n}}
        {\dist(b_{k,n},\partial B(b_{j,n},\nu_{j,n}))}
        =
        0.
\end{align}

Finally, \eqref{eq:arbitrary-profile-orthogonality} and
Lemma~\ref{lem:well-defn-scale-center} imply
\begin{align}\label{eq:localize-interaction}
        \lim_{n\to\infty}
        \sum_{\substack{j,k\in\calJ\\j\neq k}}
        \left[
        \frac{\mu_{j,n}}{\mu_{k,n}}
        +
        \frac{\mu_{k,n}}{\mu_{j,n}}
        +
        \frac{|b_{j,n}-b_{k,n}|}{\mu_{j,n}}
        \right]^{-\frac{d-2}{2}}
        =
        0.
\end{align}

Set
\begin{align}
        \mathbf W_n
        :=
        \sum_{j\in\calJ}Q_n^j,
        \quad 
        \vec\nu_n
        :=
        (\nu_n,(\nu_{j,n})_{j\in\calJ}),
        \quad
        \vec\xi_n
        :=
        (\xi_n,(\xi_{j,n})_{j\in\calJ}).
\end{align}
By \eqref{eq:localize-strong-approximation},
\eqref{eq:localize-reduced-neck},
\eqref{eq:localize-global-geometry},
\eqref{eq:localize-component-geometry},
\eqref{eq:localize-component-boundary}, and
\eqref{eq:localize-interaction},
\begin{align}
        \lim_{n\to\infty}
        \bfd_{\gamma_0}
        (v_n,\mathbf W_n;
        B(y_n,\rho_n);
        \vec\nu_n,\vec\xi_n)
        =
        0.
\end{align}
Taking the infimum in Definition~\ref{defn:delta} proves
\eqref{eq:localize-strong-conclusion}.
\end{proof}

\section{Proofs of the Main Results}
We now prove Theorems~\ref{thm:main1} and~\ref{thm:main}, and Corollary~\ref{cor:src}.

\begin{proof}[Proof of Theorem~\ref{thm:main1}]
Let $t_n\to T_+$ be an arbitrary sequence of times. Suppose first that $T_+<\infty$. Let $u^*$ and $\calS=\{x^1,\ldots,x^K\}$ be as in Theorem~\ref{thm:lwp}. By Lemma~\ref{lem:finite-time-no-escape}, $K\geq1$. Lemmas~\ref{lem:arbitrary-stationary-profiles} and~\ref{lem:arbitrary-time-no-ghost}, applied to $t_n$, give
\begin{align}
        u(t_n)
        =
        u^*
        +
        \sum_{j=1}^J W^j[a_n^j,\lambda_n^j]
        +
        o_{\dot H^1}(1).
\end{align}
With $i_j$ as in \eqref{eq:arbitrary-profile-finite-localization}, set $J_i:=\#\{1\leq j\leq J:i_j=i\}$. Lemma~\ref{lem:arbitrary-singular-coverage} gives $J_i\geq1$, and $J=\sum_{i=1}^KJ_i$. Relabel the corresponding profiles as $W_j^i[a_{j,n}^i,\lambda_{j,n}^i]$, $1\leq j\leq J_i$. The orthogonality of parameters in \eqref{eq:arbitrary-profile-orthogonality} gives \eqref{eq:ao-finite}, while
\eqref{eq:arbitrary-profile-finite-localization} gives the desired limits for the translation and scaling parameters. This proves part $\mathrm{(i)}$ of Theorem~\ref{thm:main1}.

If $T_+=\infty$, Lemmas~\ref{lem:arbitrary-stationary-profiles} and~\ref{lem:arbitrary-time-no-ghost} give
\begin{align}
        u(t_n)
        =
        \sum_{j=1}^J W^j[a_n^j,\lambda_n^j]
        +
        o_{\dot H^1}(1).
\end{align}
Set $K:=J$. The orthogonality of parameters in \eqref{eq:arbitrary-profile-orthogonality} and
\eqref{eq:arbitrary-profile-global-localization} give
\eqref{eqn:ao-global} and the desired limits for the translation and scaling parameters. This proves
part $\mathrm{(ii)}$ of Theorem~\ref{thm:main1}.
\end{proof}

\begin{proof}[Proof of Theorem~\ref{thm:main}]
Fix $0<\gamma_0<\min\{\bar E_*/2,\gamma_*\}$, where $\gamma_*$ is as in Lemma~\ref{lem:uniform-decay-stationary}. Assume first that $T_+<\infty$. Let $u^*$ and $\calS$ be as in Theorem~\ref{thm:lwp}. Fix $y\in\R^d$ and suppose that the first conclusion of part $\mathrm{(i)}$ is false. Then there are $t_n\to T_+$ and $\eta>0$ such that, with
$\rho_n:=\sqrt{T_+-t_n}$,
\begin{align}\label{eq:main-direct-finite-contradiction}
        \bs\de_{\gamma_0}(u(t_n);B(y,\rho_n))\geq\eta.
\end{align}
After passing to a subsequence,
\begin{align}\label{eq:main-direct-finite-decomp}
        u(t_n)
        =
        u^*
        +
        \sum_{j=1}^J W^j[a_n^j,\lambda_n^j]
        +o_{\dot H^1}(1),
\end{align}
and for every $1\leq j\leq J$ there is $x^{i_j}\in\calS$ such that
\begin{align}\label{eq:main-direct-finite-localization}
        \lim_{n\to \infty}\frac{|a_n^j-x^{i_j}|+\lambda_n^j}{\rho_n}=0.
\end{align}
Choose $\beta_n\to\infty$ so slowly that $\beta_n\rho_n\to0$ and
$B(y,2\beta_n\rho_n)$ contains no point of
$\calS\setminus\{y\}$. If $y\in\calS$, choose $\alpha_n\to0$ so slowly
that
\begin{align}
       \lim_{n\to \infty} \max_{j:x^{i_j}=y}
        \frac{|a_n^j-y|+\lambda_n^j}{\alpha_n\rho_n}
        =0.
\end{align}
If $y\notin\calS$, choose any sequence $\alpha_n\to0$. Set $Q_n^j:=W^j[a_n^j,\lambda_n^j]$ and $e_n:=u(t_n)-\sum_{j=1}^JQ_n^j$. Then \eqref{eq:main-direct-finite-decomp} gives
$\|e_n-u^*\|_{\dot H^1}\to0$ as $n\to \infty$.  Since $\beta_n\rho_n\to0$, the Sobolev inequality, the absolute continuity of the $\dot H^1$ and $L^{p+1}$ norms of $u^*$, and the preceding convergence give \eqref{eq:localize-strong-error}. If $x^{i_j}=y$, then
\begin{align}
\frac{\alpha_n\rho_n-|a_n^j-y|}{\lambda_n^j}
\to\infty,
\end{align}
and hence
\begin{align}
\bar E(Q_n^j;\R^d\setminus B(y,\alpha_n\rho_n))
+
\|Q_n^j\|_{L^{p+1}
(\R^d\setminus B(y,\alpha_n\rho_n))}^{p+1}
\to0.
\end{align}
If $x^{i_j}\neq y$, then
\begin{align}
\bar E(Q_n^j;B(y,\beta_n\rho_n))
+
\|Q_n^j\|_{L^{p+1}(B(y,\beta_n\rho_n))}^{p+1}
\to0.
\end{align} 
Together with
\eqref{eq:main-direct-finite-decomp} and the absolute continuity of the
$\dot H^1$ and $L^{p+1}$ norms of $u^*$, this gives
\begin{align}
        \lim_{n\to \infty} \left[\bar E(u(t_n);B(y,\beta_n\rho_n)\setminus
        B(y,\alpha_n\rho_n))
        +
        \|u(t_n)\|_{L^{p+1}(B(y,\beta_n\rho_n)\setminus
        B(y,\alpha_n\rho_n))}^{p+1}\right]
        =0.
\end{align}
Thus \eqref{eq:localize-strong-neck} holds, and Lemma~\ref{lem:localize-strong-bubbles} contradicts
\eqref{eq:main-direct-finite-contradiction}. This proves the first
conclusion in part $\mathrm{(i)}$.

Now suppose that the second conclusion of part $\mathrm{(i)}$ is false. Then there exist sequences $t_n\to T_+$, $y_n\in \R^d$, $\rho_n>0$, $R_n>0$, $\alpha_n>0$, and $\beta_n>0$ satisfying its hypotheses and $\eta>0$ such that
\begin{align}\label{eq:main-direct-local-contradiction}
\bs\de_{\gamma_0}(u(t_n);B(y_n,\rho_n))\geq\eta
\end{align}
for every $n\geq 1$. After passing to a subsequence,
Lemmas~\ref{lem:arbitrary-stationary-profiles} and
\ref{lem:arbitrary-time-no-ghost} give
\begin{align}
        u(t_n)
        &=
        \sum_{j=1}^JQ_n^j+e_n,
        \quad
        Q_n^j:=W^j[a_n^j,\lambda_n^j],
        \quad
        \|e_n-u^*\|_{\dot H^1}\to0,
\end{align}
as $n\to\infty$. The inclusion $B(y_n,R_n\rho_n)\subset B(y,\sqrt{T_+-t_n})$ and $\beta_nR_n^{-1}\to0$ give
\[
\beta_n\rho_n
\leq
\beta_nR_n^{-1}\sqrt{T_+-t_n}
\to0,
\]
as $n\to \infty.$ Hence the Sobolev inequality, the absolute continuity of the $\dot H^1$ and $L^{p+1}$ norms of $u^*$ and $\|e_n-u^*\|_{\dot H^1}\to0$ give \eqref{eq:localize-strong-error}, while the assumed annular limit is \eqref{eq:localize-strong-neck}. Lemma~\ref{lem:localize-strong-bubbles} therefore contradicts
\eqref{eq:main-direct-local-contradiction}. This proves part $\mathrm{(i)}$.

Assume now that $T_+=\infty$. Fix $y\in\R^d$ and suppose that the first
conclusion in part $\mathrm{(ii)}$ fails. Then there exist $\eta>0$ and $t_n\to\infty$ such that, with $\rho_n:=\sqrt{t_n}$,
\begin{align}\label{eqn:lower-bound}
        \bs\de_{\gamma_0}(u(t_n);B(y,\rho_n))\geq\eta.
\end{align}
After passing to a subsequence,
\begin{align}\label{eqn:preceding-limit}
        u(t_n)
        =
        \sum_{j=1}^J W^j[a_n^j,\lambda_n^j]
        +o_{\dot H^1}(1),
        \quad
        \lim_{n\to \infty}\frac{|a_n^j|+\lambda_n^j}{\rho_n}=0.
\end{align}
Set $Q_n^j:=W^j[a_n^j,\lambda_n^j]$ and $e_n:=u(t_n)-\sum_{j=1}^JQ_n^j$. If $J=0$, choose the zero multi-bubble configuration in Definition~\ref{def:d} and set $\xi_n:=\rho_n/n$ and $\nu_n:=n\rho_n$. The strong convergence of the remainder in $\dot H^1$ and the Sobolev inequality give
\begin{align}
        \bs\de_{\gamma_0}(u(t_n);B(y,\rho_n))
        &\leq
        2\bar E(u(t_n))
        +
        2\|u(t_n)\|_{L^{p+1}}^{p+1}
        +
        \frac{2}{n}
        \rightarrow0,
\end{align}
as $n\to \infty$, contradicting \eqref{eqn:lower-bound}. Thus $J\geq1$. Since $y$ is fixed and $\rho_n\to\infty$, \eqref{eqn:preceding-limit} gives
\begin{align}
\max_{1\leq j\leq J}
\frac{|a_n^j-y|+\lambda_n^j}{\rho_n}
\to0.
\end{align}
Choose $\alpha_n\to0$ so slowly that
\begin{align}
        \lim_{n\to \infty}\max_{1\leq j\leq J}
        \frac{|a_n^j-y|+\lambda_n^j}{\alpha_n\rho_n}
        =0,
\end{align}
and choose any $\beta_n\to\infty$. For every $1\leq j\leq J$,
\begin{align}
\frac{\alpha_n\rho_n-|a_n^j-y|}{\lambda_n^j}
\to \infty,
\end{align}
as $n\to \infty.$ Hence, by a change of variables and
$W^j\in\dot H^1(\R^d)\cap L^{p+1}(\R^d)$,
\begin{align}
\lim_{n\to\infty}\left[
\bar E(Q_n^j;\R^d\setminus B(y,\alpha_n\rho_n))
+
\|Q_n^j\|_{L^{p+1}(\R^d\setminus B(y,\alpha_n\rho_n))}^{p+1}
\right]
=0.
\end{align}
Together with the strong convergence of the remainder and the Sobolev inequality, this gives \eqref{eq:localize-strong-error} and \eqref{eq:localize-strong-neck} with $v_n=u(t_n)$. Hence Lemma~\ref{lem:localize-strong-bubbles} contradicts \eqref{eqn:lower-bound}. This proves the first conclusion in part $\mathrm{(ii)}$.

Finally, suppose that the second conclusion in part $\mathrm{(ii)}$
is false. Then there exist sequences $t_n$, $y_n$, $\rho_n$, $R_n$, $\alpha_n$, and $\beta_n$ satisfying its hypotheses and a constant $\eta>0$ such that
\begin{align}
\bs\de_{\gamma_0}(u(t_n);B(y_n,\rho_n))
\geq\eta.
\end{align}
After passing to a subsequence,
Lemmas~\ref{lem:arbitrary-stationary-profiles} and
\ref{lem:arbitrary-time-no-ghost} give
\begin{align}
        u(t_n)
        &=
        \sum_{j=1}^JQ_n^j+e_n,
        \quad
        Q_n^j:=W^j[a_n^j,\lambda_n^j],
        \quad
        \|e_n\|_{\dot H^1}\to0,
\end{align}
as $n\to\infty$. The convergence $\|e_n\|_{\dot H^1}\to0$ and the Sobolev inequality give \eqref{eq:localize-strong-error}, while the assumed annular limit gives \eqref{eq:localize-strong-neck}. Lemma~\ref{lem:localize-strong-bubbles} therefore gives a contradiction. This proves part $\mathrm{(ii)}$.
\end{proof}

\begin{proof}[Proof of Corollary~\ref{cor:src}]
Approximating $u_0$ in $\dot H^1$ by nonnegative smooth functions and using the maximum principle and continuous dependence, we obtain $u(t)\geq0$ for every $t\in[0,T_+)$ (cf. proof of Lemma~E.2 in~\cite{collot}). Let $W[a_n,\lambda_n]$ be any profile obtained in Theorem~\ref{thm:main1}. In the notation of \eqref{eq:arbitrary-profile-rescaled-flow}, the proof of Lemma~\ref{lem:arbitrary-stationary-profiles} gives $U_n(0)\to W$ strongly in $L^2_{\loc}(\R^d)$. Thus $W\geq 0$. Since $W$ is nonzero, elliptic regularity and the strong maximum principle give $W>0$. Corollary~8.2(b) in~\cite{CGS} therefore gives $W=\widetilde W[b,\mu]$ for some $b\in\R^d$ and $\mu>0$. Since $W[a_n,\lambda_n]=\widetilde W[a_n+\lambda_nb,\lambda_n\mu]$, every profile may be taken equal to $\widetilde W$.

By \eqref{eq:energy-identity} and \eqref{eq:tension-L2}, $E(u(t))$ has a finite limit as $t\to T_+$. If $T_+=\infty$, \eqref{eq:arbitrary-stationary-pythagorean} and \eqref{eqn:decomp-global} give
\begin{align}
\lim_{t\to\infty}E(u(t))
=
\frac{K}{d}\bar E_*.
\end{align}
Thus $K$ is independent of the chosen sequence. If $T_+<\infty$, \eqref{eq:no-ghost-finite-energy-splitting}, \eqref{eq:arbitrary-stationary-pythagorean}, and \eqref{eqn:decomp-finite} give
\begin{align}
\lim_{t\to T_+}E(u(t))
=
E(u^*)+\frac{\bar E_*}{d}\sum_{i=1}^KJ_i.
\end{align}
Thus the total number of bubbles, $\sum_{i=1}^KJ_i$, is independent of the chosen sequence.
\end{proof}
\bibliographystyle{alpha}
\bibliography{refs}
\end{document}